\documentclass[11pt]{ucthesis}

\usepackage{defs_thesis, longtable, graphicx, rotating}

\def\e{\mathbf{e}}
\def\ie{in other words,\ }
\def\eg{for example,\ }
\def\ae{as}
\def\nb{\noindent{\bf N.B.}\ }
\def\ex{\noindent{\bf Example.\ }}
\ssp

\begin{document}

\title{The Singular Theta Correspondence, Lorentzian Lattices\\and
  Borcherds--Kac--Moody Algebras}
\author{Alexander Graham Barnard}
\degreeyear{2003}
\degreesemester{Spring}
\degree{Doctor of Philosophy}
\chair{Professor Richard E. Borcherds}
\othermembers{Professor Nicolai Reshetikhin\\Professor Oliver
M. O'Reilly}
\prevdegrees{B.A. (University of Cambridge, U.K.) 1997\\
M.A. (University of Cambridge, U.K.) 2001
}
\field{Mathematics}
\campus{Berkeley}

\maketitle

\copyrightpage

\begin{abstract}
This dissertation answers some of the questions raised in Borcherds'
papers on Moonshine and Lorentzian reflection groups.  Firstly, we
show that the {\it pseudo--cusp forms} studied by Hejhal in connection
with the Riemann hypothesis can be constructed using the singular
theta correspondence.  We prove (assuming an open conjecture of
Burger, Li and Sarnak about the automorphic spectra of orthogonal
groups) that a Lorentzian reflection group with Weyl vector is
associated to a vector--valued modular form.  This result allows us to
establish a folklore conjecture that the maximal dimension of a
Lorentzian reflection group with Weyl vector is $26$.  In addition, in
the case of elementary lattices, we show that these vector--valued
forms can be obtained by inducing scalar--valued forms.  This allows
us to explain the {\it critical signatures} which occur in Borcherds'
work.  Many of the structures which occur at these critical signatures
are especially beautiful and symmetric and have appeared independently
throughout the literature.  We investigate Borcherds--Kac--Moody ({\it
BKM}) algebras with denominator formul\ae\ that are singular weight
automorphic forms.  The results of these investigations suggest that
all such BKM algebras are related to orbifold constructions of vertex
algebras and elements of the Monster finite group.  If this conjecture
were true it would give a nice, simple classification of an important
class of BKM algebras.  These BKM algebras are interesting from a
purely Lie algebraic point of view as they can be considered natural
generalizations of finite and affine Lie algebras: The Weyl groups for
finite dimensional Lie algebras are spherical reflection groups, for
affine Lie algebras are planar reflection groups and for these BKM
algebras are hyperbolic reflection groups.  Since BKM algebras also
appear in the string theory literature as algebras associated to BPS
states, we expect a classification will be of further interest.
Finally, we show how work in this dissertation combined with results
of Bruinier gives a new insight into the {\it arithmetic mirror
symmetry conjecture} of Gritsenko and Nikulin.

\end{abstract}

\begin{frontmatter}

\begin{dedication}
\begin{center}
To my grandfather ---\vskip 1cm

who was always interested in the work of others.
\end{center}
\end{dedication}

\tableofcontents

\begin{acknowledgements}
I wish to thank my advisor, Prof. Richard Borcherds, for always
providing me with interesting and challenging problems, many helpful
discussions and support throughout my years as a graduate student.  It
has been a pleasure to work with someone that has such a unique
perspective on mathematics.
\vskip 1cm
\noindent Of course, none of this would have been possible without my
parents always encouraging my curiosity and providing me with every
opportunity to discover and learn more.
\vskip 1cm
\noindent In addition, I thank Nils Scheithauer, Gerald H\"ohn, Fritz
Grunewald, Jan Bruinier, Michael Kuss, Frank Calegari and Ian
Grojnowski for many helpful conversations which have influenced this
dissertation.
\end{acknowledgements}

\end{frontmatter}

\chapter*{Notation}

\begin{longtable}[l]{ll}
$\({a \over b}\)$\quad & The Kronecker symbol; \\
$(\cdot,\cdot)$ & The bilinear form on a lattice; \\
$\<\cdot,\cdot\>$ & A bilinear form; \\
$A$ & The discriminant form $L^*/L$ of some lattice $L$; \\
$A_n$ & Elements of order $n$ in the discriminant group; \\
$A^n$ & $n$th powers of elements of $A$; \\
$A^{n*}$ & A coset of $A^n$ in $A$; \\
$b^+, b^-$ & Signatures of a lattice; \\
$b(\lambda, n)$ & A Fourier coefficient of a modular form; \\
$\BBB_n(x)$ & A Bernoulli polynomial of degree $n$; \\
$\chi_A$ & A character associated to the discriminant form $A$; \\
$c_{\lambda,n}(y)$ & A Fourier coefficient of a modular form; \\
$\CCC$ & The complex numbers; \\
$\Delta$ & A Laplacian operator; \\
$\Delta(\tau)$ & The unique normalized cusp form of weight $12$ for
  $\lieSL_2(\ZZZ)$; \\
$\Delta_{n+}(\tau)$ & A cusp form for the group $\Gamma_0(n)+$; \\
$\ee_\gamma$ & A basis for $\CCC[L^*/L]$; \\
$\e(x)$ & The function $e^{2\pi i x}$; \\
$\eta(\tau)$ & Dedekind's eta function; \\
$F_L(\tau; \gamma,s)$ & A Maass--Poincar\'e\ series; \\
$\F$ & The usual fundamental domain for $\lieSL_2(\ZZZ)$; \\
$\F_t$ & A truncated version of $\F$; \\
$\Gamma(s)$ & The Gamma function; \\
$\tilde\Gamma_0(N)$ & A congruence subgroup of $\lieMP_2(\ZZZ)$; \\
$\Gr(L)$ & The Grassmannian manifold of $L$; \\
$G_L(a,\beta,c)$ & A Gauss sum; \\
$\GGG_n$ & A Hecke triangle group; \\
$H(\delta)$ & The Heegner--divisor corresponding to $\delta$; \\
$\H$ & The upper half plane; \\
$\Im(\tau)$ & The imaginary part of $\tau$; \\
$\even_{1,1}$ & A lattice with bilinear form
  $\smallmatrix{0}{1}{1}{0}$; \\
$\even_{r,s}(A)$ & A genus of even lattices; \\
$K$ & The lattice $(L\cap z^\perp)/\ZZZ z$; \\
$K_n(y)$ & A modified Bessel function of the third kind; \\
$k^+, k^-$ & Weights of a modular form; \\
$\Lambda$ & The Leech lattice; \\
$L$ & A lattice; \\
$L^*$ & The dual of a lattice; \\
$L_0^*$ & A sublattice of $L^*$; \\
$L(n)$ & The lattice with form scaled by $n$; \\
$\lieMP_2$ & The metaplectic group --- a double cover of $\lieSL_2$; \\
$\M_{n,s}(y)$ & A modified Whittaker function; \\
$N$ & The level of a lattice or modular form; \\
$\QQQ$ & The rational numbers; \\
$\QQQ_p$ & The $p$--adic numbers; \\
$\Phi(v)$ & Singular theta transformation of some function; \\
$\psi$ & Piecewise linear part of a function; \\
$q$ & $e^{2\pi i \tau}$; \\
$\rho$ & A (generalized) Weyl vector; \\
$\rho_L$ & The Weil representation associated to $L$; \\ 
$\rho_{\alpha\beta}$ & A coefficient of the Weil representation; \\
$\Re(\tau)$ & The real part of $\tau$; \\
$\RRR$ & The real numbers; \\
$\lieSL_2$ & The special linear group; \\
$\sigma_\vecv$ & The reflection in $\vecv^\perp$; \\
$S$ & $\(\smallmatrix{0}{-1}{1}{0}, \sqrt{\tau}\)\in \lieMP_2(\ZZZ)$; \\
$s$ & A complex number giving the eigenvalue; \\
$\sgn(L)$ & The signature of the lattice $L$; \\
$\tau$ & A complex number from the upper half plane; \\
$\theta_L(\tau)$ & A theta function; \\
$\Theta_L(\tau; v)$ & A Siegel theta function associated to an
  indefinite lattice $L$; \\
$T$ & $\(\smallmatrix{1}{1}{0}{1}, 1\)\in \lieMP_2(\ZZZ)$; \\
$T_\vecx$ & An automorphism of $L\tensor\RRR$; \\
$v$ & An element of the Grassmannian; \\
$v_+$ & Projection onto the positive definite space represented by
  $v$; \\
$v_-$ & Projection onto the negative definite space represented by
  $v$; \\
$W$ & The Weyl group; \\
$\W_{n,s}(y)$ & A modified Whittaker function; \\
$\xi$ & Real analytic piece of a function; \\
$\zeta(s)$ & The Riemann zeta function; \\
$Z$ & $\(\smallmatrix{-1}{0}{0}{-1}, i\)\in \lieMP_2(\ZZZ)$; \\
$z$ & A primitive norm zero vector from $L$; \\
$z'$ & A vector in $L^*$ with $(z,z') = 1$; \\
$\ZZZ$ & The integers; \\
$\ZZZ_p$ & The $p$--adic integers. \\
\end{longtable}

\addtocounter{chapter}{-1}
\chapter{Summary of Results}

The majority of this dissertation uses Harvey and Moore's extension of
the theta correspondence \cite{borcherds_grass, harvey_moore} (which
we call the {\it singular theta correspondence}) to study Lorentzian
reflection groups with Weyl vectors.  We start by recalling some facts
about the singular theta correspondence.

Most of the definitions can be found in Chapter $1$.  Suppose that $L$
is an even, integral, unimodular lattice of signature $(m,n)$ and that
$\Gr(L)$ is the Grassmannian of maximal ($m$--dimensional) positive
definite subspaces of $L\tensor \RRR$.  Let $F(\tau)$ be a holomorphic
modular form of weight $(m-n)/2$ for $\lieSL_2(\ZZZ)$.  The Siegel
theta function $\Theta_L(\tau; v)$ is a function of both $v \in
\Gr(L)$ and $\tau$ in the upper half plane and it is invariant under
the actions of $\Aut(L)$ and $\lieSL_2(\ZZZ)$.  The function
\begin{equation}
\Phi_F(v) = \int_{\F} \overline{\Theta_L(\tau; v)} F(\tau)
dxdy/y^2, \label{intro_eqn}
\end{equation}
where $\F$ is a fundamental domain for $\lieSL_2(\ZZZ)$ acting on the
upper half plane, is an automorphic form on $\Gr(L)$ invariant under
the discrete group $\lieO_L(\ZZZ)$.  The correspondence $F(\tau)
\longrightarrow \Phi_F(v)$ is roughly the theta (or Howe)
correspondence.

Suppose that we allow $F(\tau)$ to have singularities at the cusps but
require it to be holomorphic on the upper half place.  The integral in
Equation \ref{intro_eqn} then diverges wildly.  Harvey and Moore
showed, using ideas from quantum field theory, that it is still
possible to make sense of the integral by ``regularizing'' it.  This
allowed Harvey and Moore to simplify the proofs of many results from
\cite{borcherds_products}.  Their construction was generalized by
Borcherds in \cite{borcherds_grass} to, in particular, non--unimodular
lattices.  The main difference when we work with non--unimodular
lattices is that we have to replace the modular form $F(\tau)$ with a
vector--valued modular form transforming according to the Weil
representation of the lattice $L$.

It is possible to use real analytic modular forms in the singular
theta correspondence.  Real analytic modular forms are functions with
similar transformation properties to holomorphic modular forms, but
rather than being holomorphic they are assumed to be eigenfunctions of
the Laplacian.  These are a generalization of holomorphic modular
forms since being holomorphic is essentially equivalent to being
annihilated by the Laplacian.  An advantage of real analytic modular
forms is that they exist for any choice of singularities at the cusps.
A real analytic form with given singularities at cusps can be
constructed by standard Eisenstein series techniques.

In Chapter $2$ we show that the singular theta correspondence applied
to suitable real analytic modular forms can give cusp forms for
$\lieSL_2(\ZZZ)$ with logarithmic singularities at the corners of the
fundamental domain.  These functions, known as {\it pseudo--cusp
forms}, were first found in computer experiments of Haass.  He was
searching for eigenfunction of the Laplacian and found pseudo--cusp
forms accidentally due to a slight bug in his code. Later, Stark
noticed that the eigenvalues for which these pseudo--cusp forms occur
are closely related to zeros of the Riemann zeta function.  This was
explained by Hejhal who showed they were related to the zeros of the
Dedekind zeta function for $\QQQ(\sqrt{-3})$.  Chapter $2$ provides a
new construction for these pseudo--cusp forms.  It is simple to find a
real analytic function, $F(\tau)$, whose singular theta transformation
is a modular form with logarithmic singularities; what is difficult is
to ensure it is cuspidal.  This can happen in two distinct ways:
either the singular theta transformation is automatically cuspidal or
the constant term of $F(\tau)$ vanishes.  The first of these occurs
when the eigenvalue of the Laplacian is a root of the Riemann zeta
function; the second occurs when the eigenvalue is a root of a certain
L--function.  So, cuspidal behaviour occurs when the eigenvalue is a
root of the product of the Riemann zeta function and the L--function;
this product is exactly the Dedekind zeta function for
$\QQQ(\sqrt{-3})$.  This construction generalizes easily to cusp forms
with singularities at other points in the fundamental domain and
explains the appearance of other Dedekind zeta functions.

The singular theta correspondence for holomorphic modular forms, when
the lattice has signature $(1,n)$, gives piecewise linear functions on
the hyperbolic space $\Gr(L)$.  These piecewise linear functions have
singularities orthogonal to certain vectors in $L^*$ (the vectors
being determined by the singularities of the modular form).  In rare
cases the singularities will occur along the reflection hyperplanes of
the underlying lattice $L$.  When this happens we say that the modular
form $F(\tau)$ is associated to the reflection group of $L$ (we call
$F(\tau)$ a {\it reflective} modular form).  Much evidence is
presented in \cite{borcherds_refl} supporting the conjecture that all
``nice'' Lorentzian reflection groups are associated, in this manner,
to reflective forms.  Chapter $3$ provides a proof of this conjecture
assuming an open conjecture of Burger, Li and Sarnak
\cite{sarnak_ramanujan, sarnak_ramanujan2}.  The singular theta
correspondence applied to real analytic forms does not quite give
piecewise linear functions; it gives a sum of a piecewise linear
function and a smooth function.  Therefore the concept of reflective
form still makes sense in the real analytic case.  Since real analytic
modular forms can have any desired singularity structure at the cusps
it is easy to write down a real analytic modular form associated to
the reflection group of $L$.  Specializing the eigenvalue of the real
analytic modular form to be the eigenvalue which holomorphic forms
have we obtain a form $F(\tau)$, some of whose Fourier coefficients
are holomorphic.  To complete the theorem we have to show that the
remaining Fourier coefficients vanish.  Bruinier has shown that the
holomorphic coefficients correspond to the piecewise linear function
and the non--holomorphic coefficients to the smooth function.  He also
showed that the smooth function is an eigenfunction of the Laplacian
with eigenvalue $-n$.  The Laplacian for orthogonal groups has been
extensively studied and there is a conjectural spectrum.  The
eigenvalue $-n$ does not occur in this spectrum so we deduce (assuming
the conjecture) that the smooth function is identically zero.  For
most lattices this is enough to show that the individual
non--holomorphic terms are also zero and hence that $F(\tau)$ is a
holomorphic modular form.

We now have a correspondence between Lorentzian lattices with Weyl
vectors and holomorphic vector--valued modular forms with certain
singularities at the cusps.  This allows us to deduce results about
Lorentzian lattices by studying the related modular forms.  This is
good: the theory of modular forms is further developed and more
powerful tools are available.  One theorem we can deduce a folklore
conjecture that $\even_{1,25}$ is (essentially) the unique largest
Lorentzian lattice with Weyl vector.  Under the correspondence this
conjecture is equivalent to the lowest weight cusp form for
$\lieSL_2(\ZZZ)$ having weight $12$; this is a well--known and basic
result in the theory of modular forms.  We sketch a more detailed
proof below.

Primitive roots of lattices have norms less than $2$, therefore the
associated reflective form, $F(\tau)$, has singularities no worse than
$q^{-1}$.  In particular, this means that $\Delta(\tau)F(\tau)$ is
holomorphic even at the cusps (here $\Delta(\tau)$ is Ramanujan's cusp
form for $\lieSL_2(\ZZZ)$ of weight $12$).  Forms which are
holomorphic everywhere must have positive weight, so $F(\tau)$ must
have weight at least $-12$ and therefore the lattice $L$ must have
signature at least $-24$.  This implies the folklore conjecture that
the maximal signature for a Lorentzian lattice with Weyl vector is
$-24$.  In fact, since we know the space of all vector--valued
holomorphic modular forms of weight $0$ we can deduce that all
Lorentzian lattices with Weyl vectors and signature $-24$ are closely
related to $\even_{1,25}$.

In \cite{borcherds_refl} it was observed that there is a natural way
to obtain vector--valued forms for $\lieSL_2(\ZZZ)$ by inducing
scalar--valued forms of certain congruence subgroups and that many of
the vector--valued forms associated to Lorentzian lattices occurred in
this way.  In Chapter $4$ we show that if the Lorentzian lattice $L$
is elementary then the associated vector--valued form is induced from
a scalar--valued form.  This is proved by examining in detail how the
induction process behaves on the various components of the
vector--valued form.  This result is useful since scalar--valued forms
are better understood than vector--valued ones and easier for a
computer program to handle.  One consequence of this result is an
explanation for the {\it critical signatures} observed in
\cite{borcherds_refl}: These signatures correspond to the lowest
weights of cusp forms for various congruence groups.  For example, in
level $1$ the critical signature is $-24$, the corresponding
Lorentzian lattice is $\even_{1,25}$ and the corresponding definite
lattice is the Leech lattice.  In level $2$ the critical signature is
$-16$, the Lorentzian lattice is $\even_{1,17}(2^{+10})$ and the
definite lattice is the Barnes--Wall lattice.  In level $3$ the
critical signature is $-12$, the Lorentzian lattice is
$\even_{1,13}(3^{-8})$ and the definite lattice is the Coxeter--Todd
lattice.  These three lattices are well--known and have many similar
and beautiful properties.  Other lattices which occur at the critical
signatures are slightly less well--known but have similar properties.
For example, the above three examples fit into a family: For primes
$p$ such that $p+1$ divides $24$ we obtain Lorentzian lattices which
are closely related to modular lattices (a modular lattice is one
which is similar to its dual lattice).

If the discriminant form for the Lorentzian lattice $L$ is not too
small we can deduce strong restrictions on the singularities that can
occur in the corresponding scalar--valued form.  These restrictions
allowed us to calculate the critical signatures.  However, if the
discriminant form is small the scalar--valued form has a much weaker
structure and we can not deduce results as strong as those stated
above.  Hence, there exist Lorentzian lattices with signatures below
the critical bound (although they are fairly rare).  These lattices
still have interesting properties: for example, they often have finite
co--volume reflection groups.  However, one of the most interesting
properties a Lorentzian reflection group can have is possessing a norm
$0$ Weyl vector, and none of these have one.  A large amount of
computational data and various methods for dealing with such {\it
irregular lattices} are collected at the end of Chapter $4$.

A particularly important situation where Lorentzian reflection groups
occur is as Weyl groups of Borcherds--Kac--Moody ({\it BKM}) algebras.
These algebras seem to be particularly simple when the Weyl vector has
norm $0$ and the denominator formula is a {\it singular weight}
automorphic form (actually, this latter condition implies that the
Weyl vector is of norm $0$).  One reason to restrict to this class of
BKM algebras is that the dimensions of its root spaces can be
expressed by simple closed expressions involving Fourier coefficients
of modular forms.  Examples of such BKM algebras are the fake Monster
Lie algebra \cite{borcherds_fake}, which is related to the modular
form $\Delta(\tau)$, and the Monster Lie algebra
\cite{borcherds_moonshine}, which is related to the modular form
$j(\tau)$.  Various other such algebras have been constructed, for
example, the fake Monster superalgebra \cite{nils_fakesuper} and the
Baby Monster superalgebra \cite{hohn_babysuper}.  To study such BKM
algebras, in Chapter $5$ we examine the automorphic forms which could
occur as their denominator formul\ae.  To do this we use the work of
Bruinier which shows that such automorphic forms should come from the
singular theta transformation applied to vector--valued modular forms
in the signature $(2,n)$ situation.  This means we can study the
modular forms instead of the automorphic forms.  By multiplying by a
suitable Eisenstein series we obtain a linear relation satisfied by
the coefficients of these vector--valued modular forms.  As the
coefficients that occur are $0$ or $1$, even a single linear equation
is a strong restriction.  We wrote a computer program (see Appendix
\ref{programs}) to compute the coefficients of these Eisenstein series
and then to search for solutions to the linear equation.
Surprisingly, we find very few solutions; so few that it looks like
all such BKM algebras can be constructed in a uniform way.
Specifically, it seems likely that all such BKM algebras come from the
``orbifold'' construction well--known in the theory of vertex
algebras.

These singular weight automorphic forms occur in other places in the
literature: they play a fundamental role in the {\it arithmetic mirror
symmetry conjecture} of Gritsenko and Nikulin \cite{mirror_K3,
mirror_CY}.  This conjecture relates interesting Lorentzian reflection
groups to the cusps of {\it reflective automorphic forms} in the
signature $(2,n)$ case.  From the ideas in this dissertation and the
work of Bruinier we obtain new insight into this conjecture.  We have
seen how interesting Lorentzian reflection groups are closely related
to vector--valued modular forms, and the work of Bruinier shows how
reflective automorphic forms in the signature $(2,n)$ case are
similarly related to vector--valued modular forms.  Restricting to a
cusp has a simple interpretation at the level of the vector--valued
modular forms.  Hence, we see one reason for the arithmetic mirror
symmetry conjecture to be true: both objects involved are related to
vector--valued modular forms.  Indeed, in the cases where the
correspondence in this dissertation and the correspondence in
Bruinier's work applies we are able to deduce the arithmetic mirror
symmetry conjecture; this is discussed at the end of Chapter $5$.

\chapter{Introduction}

In this chapter we recall some basic definitions.  We give references
to places in the literature where similar definitions, theorems and
calculations can be found.

\section{Lattices}

A {\it lattice} is a free $\ZZZ$--module of finite rank, equipped with
a symmetric $\ZZZ$--valued bilinear form $(*,*)$.  A lattice $L$ is
called {\it even} if the associated quadratic form,
$$q(x) = {1\over 2}(x,x),$$
takes only integral values; otherwise it is called {\it odd}.  The
property of being even or odd is known as the {\it type} of the
lattice.  Most of the lattices we deal with will be even.  If the
bilinear form is non--degenerate we say that the lattice is {\it
non--degenerate}; unless otherwise stated all lattices in this
dissertation can be assumed to be non--degenerate.  The {\it
signature} of the lattice $L$ is the signature of the vector space
$L\tensor \RRR$ equipped with the natural $\RRR$--valued bilinear
form.  The signature is denoted by $(m,n)$ where $m$ is the dimension
of the maximal positive definite subspace of $L\tensor \RRR$ and $n$
is the dimension of the maximal negative definite subspace.  If the
signature is $(1,n)$ we call the lattice {\it Lorentzian}.  We also
define
$$\sgn(L) = m-n$$
if the lattice $L$ has signature $(m,n)$.

\ex Many important lattices can be found in Chapter 4 of
\cite{conway_sloane}.  These include the Leech lattice ($\Lambda$),
Barnes--Wall lattice ($BW_{16}$), Coxeter--Todd lattice ($K_{12}$),
and the $E_8$ lattice, all of which we will see later.

If the dimension of the lattice is sufficiently small it is common to
write the matrix representing the quadratic form as a symbol for the
lattice.  For instance, the symbol $\smallmatrix{0}{1}{1}{0}$ means
the lattice $\ZZZ^2$, equipped with the bilinear form,
$$\(x,y\)^2 = \[\begin{array}{c} x\\y\end{array}\]^{T}
\[\begin{array}{cc}0&1\\1&0\end{array}\]
\[\begin{array}{c} x\\y\end{array}\] = 2xy.$$
The lattice $\smallmatrix{0}{1}{1}{0}$ is very important, and is also
known as a {\it hyperbolic plane}.

The {\it dual lattice} of $L$ is
$$L^* = \{ x\in L\tensor\QQQ : (x,y)\in\ZZZ\hbox{ for all }y\in L\}.$$
The dual lattice comes naturally equipped with a $\QQQ$--valued
bilinear form.  The {\it discriminant form} of $L$ is the finite
Abelian group $A = L^*/L$ \cite{conway_sloane,
nikulin_discriminant}.  The order of this group is the {\it
determinant} of the lattice $L$.  The {\it $p$--rank} of $L$ is the
order of the $p$--part of $A$.  A lattice $L$ is called {\it
elementary} if the discriminant form is an elementary Abelian group.

Many of the modular forms we use take values in the vector space
$\CCC[L^*/L]$.  The natural basis elements for this vector space are
denoted by $\ee_\gamma$.  The discriminant form comes naturally with a
$\QQQ/\ZZZ$--valued bilinear form.  If $L$ is an even lattice then
$q(x)$ is a well--defined $\QQQ/\ZZZ$--valued quadratic form on
$L^*/L$.  So, for an even lattice $L$, every element of $L^*/L$ has a
well defined norm in $\QQQ/2\ZZZ$ given by $x^2$ for any
representative $x$ of the coset.  If $L$ is an even lattice the
signature is determined modulo $8$ by Milgram's formula (see Appendix
4 in \cite{milnor_husemoller}):
$$\sum_{\gamma \in L^*/L} \e(\gamma^2/2) = \sqrt{|L^*/L|}\ 
\e(\sgn(L)/8),$$
where
$$\e(x) = \exp(2\pi i\cdot x).$$
The sum on the left hand side is well defined by the preceding
discussion.  The {\it level} of a lattice $L$ is the minimal positive
integer $N$ such that $N\gamma^2/2$ is integral for all $\gamma\in
L^*$.

We define a second inner product on $\CCC[L^*/L]$ that is linear in
the first argument, anti--linear in the second and such that
$$\<\ee_\beta, \ee_\gamma\> = \delta_{\beta, \gamma} = \begin{cases}
1 & \hbox{if $\beta = \gamma$,} \\
0 & \hbox{otherwise.}\end{cases}$$

For any $n\in\ZZZ$ and $A = L^*/L$ we define $A^n$ to be the $n$th
powers of elements of $A$ and $A_n$ to be the elements of order $n$.
It is clear that $A^n$ and $A_n$ are orthogonal to each other.  We
define $A^{n*}$ to be the elements $\beta \in A$ such that $(\beta,
\gamma) \equiv n\gamma^2/2\ (\mod 1)$ for all $\gamma \in A_n$.
Clearly, when $A^{n*}$ is non--empty, $A^{n*}$ is a coset of $A^n$.
See \cite{borcherds_refl} for more about these sets.

Let $L$ be a lattice and $L^*$ its dual lattice.  A {\it primitive
vector} of $L^*$ is one that is not an integer multiple of any
smaller vector.  In other words $\vecv$ is primitive if and only if
$\QQQ \vecv \cap L^* = \ZZZ \vecv$.

A {\it root} of a Lorentzian lattice $L$ is any vector in $L^*$ of
negative norm for which the reflection in the plane orthogonal to
it is an automorphism of the lattice $L$.  The {\it reflection} in a
plane orthogonal to a vector $\vecv$ is
$$\sigma_\vecv : \vecx \longmapsto \vecx - 2{(\vecv,\vecx) \over
\vecv^2} \vecv,$$
where $\vecv^2$ is shorthand for $(\vecv,\vecv)$.

Given $L$, an integral lattice, $L\tensor \ZZZ_p$ is a lattice over the
$p$--adic integers, $\ZZZ_p$.  Two lattices are said to be in the same
{\it genus} if they are equivalent over the $p$--adic integers for all
primes $p$ and have the same signature (this is the same as being
equivalent over the completion at the infinite prime).  Two lattices
$K$ and $L$ are in the same genus if and only if $K\oplus
\smallmatrix{0}{1}{1}{0}$ and $L\oplus \smallmatrix{0}{1}{1}{0}$ are
equivalent over $\ZZZ$.  For indefinite forms the genus usually
consists of only one equivalence class of lattices; however for
definite forms this is usually far from the case.

\ex The genus of the Leech lattice consists of $24$ lattices (the
Niemeier lattices \cite{niemeier}).

Knowing the genus of a lattice is equivalent to knowing the signature,
dimension, type and discriminant form, $A$, of the lattice.  We follow
Conway's notation (see Chapter 15 of \cite{conway_sloane}) and assign
a symbol $\even_{r,s}(A')$ to the genus.  The signature of the lattice
is $(r,s)$ and the symbol $\even$ means the lattice is even.  $A'$ is
a symbol that represents the $p$--adic discriminant forms of the
lattice for all primes $p$ (these are known as the $p$--{\it Jordan
components} of the lattice).  Knowing all of the Jordan components is
equivalent to knowing the discriminant form $A$.  We describe below
the possible Jordan components:
\begin{enumerate}
\renewcommand{\labelenumi}{\roman{enumi}.}
\item Let $q>1$ be a power of an odd prime $p$.  The non--trivial
Jordan components of exponent $q$ are denoted by $q^{\pm n}$ for $n\ge
1$.  The indecomposable components are $q^{\pm 1}$, generated by an
element $\gamma$ with $q\gamma = 0$, $\gamma^2/2 = a/q$ where $a$ is
an integer with $\({2a \over p}\) = \pm 1$.  These components all have
level $q$.
\item Let $q$ be a power of $2$.  The non--trivial even Jordan
components of exponent $q$ are denoted by $q^{\pm 2n}$ for $n\ge 1$.
The indecomposable components are $q^{\pm 2}$, generated by $2$
elements $\gamma$ and $\delta$ with $q\gamma = q\delta = 0$, $(\gamma,
\delta) = 1/q$ and $\gamma^2/2 = \delta^2/2 = 0$ for $q^{+2}$,
$\gamma^2/2 = \delta^2/2 = 1/q$ for $q^{-2}$.  These components all
have level $q$.
\item Let $q$ be a power of $2$.  The non--trivial odd Jordan
components of exponent $q$ are denoted by $q_t^{\pm n}$ for $n\ge 1$
and $t\in \ZZZ/8\ZZZ$.  The indecomposable components are $q^{\pm
1}_t$ where $\({t \over 2}\) = \pm 1$, generated by an element
$\gamma$ with $q\gamma = 0$, $\gamma^2/2 = t/2q$.  These components
all have level $2q$.
\end{enumerate}

The sum of two Jordan components with the same $q$ is given by
multiplying the signs, adding the exponents $n$ and adding the
subscripts $t$ (where the empty subscript counts as $t=0$).

\ex $q^{+3} q^{+5} = q^{+8}$\ because both signs are positive, so
the product is positive, the sum of $3$ and $5$ is $8$ and there are
no subscripts.

\ex $q^{+3} q^{-5} = q^{-8}$\ because one sign is positive and one
negative, so the product is negative, the sum of $3$ and $5$ is $8$
and there are no subscripts.

\ex $2^{-2} 2^{-1}_3 = 2^{+3}_3$\ because both signs are negative,
so the product is positive, the sum of $2$ and $1$ is $3$ and the
subscripts are $0$ and $3$, summing to $3$.

Not all possible Jordan symbols defined above exist, \eg $2_4^{+2}$
does not exist; some of the Jordan components are isomorphic, \eg
$2_1^{+1} \cong 2_5^{-1}$.

\ex The $A_2$ lattice (also known as the face centred cubic lattice)
is in genus $\even_{2,0}(3^{-1})$.  From this we can read that the
lattice is $2$ dimensional, even, positive definite with determinant
$3$.  The discriminant group has $3$ elements, one with norm $q(x) =
0$ and two with norm $q(x) = 1/3$.

If $L$ is a positive definite lattice then its {\it theta function}
is
$$\theta_L(\tau) = \sum_{\gamma\in L} q^{\gamma^2/2}.$$
If we regard $q$ as $e^{2\pi i \tau}$, for $\tau$ in the upper half
plane, then the theta function is a modular form of weight $\dim(L)/2$
for the group $\Gamma_0(N)$, where $N$ is the level of the lattice
$L$.  In the case of an indefinite lattice the above sum will not
converge, so we have to define the theta function in a different way.
We will see how to do this in the next section.

\ex The theta function for the Leech lattice is
$$\theta_\Lambda(\tau) = 1 + 196560q^2 + 16773120q^3 + 398034000q^4 +
\cdots,$$
it is a modular form of weight $12$ for the group $\lieSL_2(\ZZZ)$.

\section{Vector--Valued Modular Forms}

In this section we recall the definition of a vector--valued modular
form \cite{borcherds_grass}.  We will discuss both the real
analytic and holomorphic case.

When dealing with theta functions of lattices, half--integral weight
modular forms naturally occur.  This means keeping careful track of
exactly which sign the square root should take.  One of the easiest
ways to do this is to use a double cover of the modular group.  If we
wanted to deal with more general weights we could use the universal
covering group instead; this will not be needed here.
$\lieSL_2(\RRR)$ has a non--trivial double cover $\lieMP_2(\RRR)$ (the
{\it metaplectic group}).  Elements of this group can be written as
pairs
$$\(\[ \begin{array}{cc}a&b\\c&d\end{array}\] , \pm\sqrt{c\tau+d}\).$$
The matrix is in $\lieSL_2(\RRR)$ and $\tau$ is a formal complex
variable in the upper half plane $\H$ (\ie $\Im(\tau)>0$).  The multiplication
of two elements of $\lieMP_2(\RRR)$ is given by
$$(A,f(\tau))(B,g(\tau)) = (AB,f(B(\tau))g(\tau)).$$
The matrices act on the upper half plane as M\" obius transformations
$$\[ \begin{array}{cc}a&b\\c&d\end{array}\] : \tau \longmapsto
{a\tau + b\over c\tau + d}.$$
The group $\lieMP_2(\ZZZ)$ is defined to be the inverse image of
$\lieSL_2(\ZZZ)$ in $\lieMP_2(\RRR)$.  It is generated by two elements
$T$ and $S$, where
$$T = \(\[\begin{array}{rr}1&1\\0&1\end{array}\],+1\),\qquad
S = \(\[\begin{array}{rr}0&-1\\1&0\end{array}\],+\sqrt{\tau}\).$$
These satisfy the relations $S^2 = (ST)^3 = Z$, $Z^4 = 1$, where
$$Z = \(\[\begin{array}{rr}-1&0\\0&-1\end{array}\],\ i\)$$
is a generator of the centre of $\lieMP_2(\ZZZ)$.

Suppose that $\rho$ is a representation of $\lieMP_2(\ZZZ)$ on a
finite--dimensional complex vector space $V_\rho$.  Choose
$(k^+,k^-)\in {1\over 2}\ZZZ^2$.  A {\it nearly holomorphic modular
form of weight $(k^+,k^-)$ and type $\rho$} \cite{borcherds_grass}
is defined to be a holomorphic function $F$ on the upper half plane
with values in the vector space $V_\rho$ such that it transforms under
elements of $\lieMP_2(\ZZZ)$ as
\begin{equation}
F\({a\tau+b \over c\tau +d}\) = \sqrt{c\tau+d}^{\,2k^+}
\sqrt{c\bar\tau + d}^{\,2k^-} \rho\(\(\[
\begin{array}{cc}a&b\\c&d\end{array}\], \sqrt{c\tau+d}\)\)F(\tau).
\label{transformation}
\end{equation}
The function $F$ is allowed to have poles at the cusp $i\infty$.

The {\it hyperbolic Laplacian of weight $(k^+,k^-)$} on the upper half
plane \cite{maass} is
$$\Delta = -y^2\( {\d^2 \over \d x^2} + {\d^2 \over \d y^2}\) + iy \(
(k^+-k^-){\d \over \d x} + i(k^++k^-){\d \over \d y}\).$$
It will always be clear from context what the weight is, so we do not
bother to write things like $\Delta_{k^+,k^-}$.  This definition of
the Laplacian is equivalent to that of the Casimir operator on
$\lieMP_2(\RRR)$ restricted to the upper half plane.  Recall that the
upper half plane can be realized as the quotient of $\lieMP_2(\RRR)$
by the double cover of $\lieO_2(\RRR)$ (a maximal compact subgroup).
This identification allows us to transfer differential operators on
$\lieMP_2(\RRR)$ to differential operators on the upper half plane
after assuming some transformation behaviour under the double cover of
$\lieO_2(\RRR)$.  This transformation behaviour is determined by the
weight, see Equation (\ref{character}).

Choose $\lambda\in\CCC$.  A {\it real analytic modular form of weight
$(k^+,k^-)$, type $\rho$ and eigenvalue $\lambda$} is a real analytic
function $F$ on the upper half plane that transforms as shown in
Equation (\ref{transformation}) and is an eigenfunction of the
hyperbolic Laplacian with eigenvalue $\lambda$.  It is allowed to have
poles at the cusp $i\infty$.

It is easy to see that any power of $q$ is an eigenfunction of the
hyperbolic Laplacian with eigenvalue $0$, hence holomorphic modular
forms are real analytic modular forms with eigenvalue $0$.

One way to think about these functions is as functions on the group
$\lieMP_2(\RRR)$.  This is done by defining
$$G(\gamma) = F(\gamma i)\ \sqrt{ci+d}^{\ -2k^+} \sqrt{-ci+d}^{\ -2k^-},$$
where
$$\gamma = \(\[\begin{array}{cc} a&b \\ c&d\end{array}\],
\sqrt{c\tau + d}\) \in \lieMP_2(\RRR).$$
The resulting function $G$ is a vector--valued
function on $\lieMP_2(\RRR)$ that transforms by the representation
$\rho$ under the left action of the discrete subgroup
$\lieMP_2(\ZZZ)$, \ie
$$G(g\gamma) = \rho(g) G(\gamma),\quad\hbox{for $g\in \lieMP_2(\ZZZ)$,}$$
and by the character
$$\chi(\theta) = e^{i(k^+-k^-)\theta}$$
under the right action of the maximal compact subgroup, the double cover of
$\lieO_2(\RRR)$, \ie
\begin{equation}
G(\gamma g_\theta) = \chi(\theta) G(\gamma),\quad\hbox{for }g_\theta
= \(\[\begin{array}{cc}\cos\theta & -\sin\theta\\\sin\theta &
\cos\theta\end{array}\], \pm\sqrt{\cos\theta + \tau\sin\theta}\),
\label{character}
\end{equation}
where the sign is $+$ if $0\le \theta < 2\pi$ and $-$ if $2\pi \le
\theta < 4\pi$.

If $L$ is a lattice, then we define the {\it Grassmannian manifold},
$\Gr(L)$, to be the set of maximal positive definite subspaces of
$L\tensor\RRR$ \cite{borcherds_grass}.  It is a symmetric space
acted on by the orthogonal group $\lieO_L(\RRR)$ (the group of
transformation preserving the norm induced from $L$).  If $v$ is an
element of $\Gr(L)$ and $\lambda\in L\tensor\RRR$ then we denote by
$\lambda_{v^+}$ the projection of $\lambda$ onto the positive definite
space represented by $v$.  Similarly, we denote by $\lambda_{v^-}$ the
projection of $\lambda$ onto the negative definite space orthogonal to
$v$.  If the lattice is either positive or negative definite then the
Grassmannian consists of a single point; in general the Grassmannian
is a connected manifold of dimension $b^+\times b^-$.

Suppose $L$ is an even lattice.  The {\it Siegel theta function} of a
coset $L+\gamma$ of $L$ in $L^*$ is
$$\theta_{L+\gamma}(\tau;v) = \sum_{\lambda\in L+\gamma}
\e\({\tau\lambda_{v^+}^2 \over 2} + {\bar\tau\lambda_{v^-}^2 \over
2}\),$$
for $\tau$ in the upper half plane and $v\in\Gr(L)$.  Combining these
for all elements of $L^*/L$ gives a $\CCC[L^*/L]$--valued function,
the {\it Siegel $\Theta$--function} of $L$
$$\Theta_L(\tau;v) = \sum_{\gamma\in L^*/L} \ee_\gamma
\theta_{L+\gamma}(\tau;v).$$
If the lattice $L$ has signature $(b^+,b^-)$ then the function
$\Theta_L(\tau;v)$ is a modular form of weight $(b^+/2,b^-/2)$ and
type $\rho_L$, where $\rho_L$ is the {\it Weil representation} of the
group $\lieMP_2(\ZZZ)$ on the vector space $\CCC[L^*/L]$.  In terms of
generating elements the Weil representation is given by
\begin{eqnarray*}
\rho_L(T):\ee_\gamma &\longmapsto& \e\({\textstyle {1\over 2}}
(\gamma,\gamma)\) \ee_\gamma, \\
\rho_L(S):\ee_\gamma &\longmapsto& {\e(-\sgn(L)/8) \over \sqrt{|L^*/L|}}
\sum_{\delta\in L^*/L}\e(-(\gamma,\delta))\ee_\delta.
\end{eqnarray*}
See \cite{borcherds_grass} for proofs of these properties of Siegel
theta functions.

Let $\beta,\gamma\in L^*/L$, then the {\it
$\beta,\gamma$--coefficient} of the Weil representation, denoted by
$\rho_{\beta\gamma}$, is given by
$$\rho_{\beta\gamma}(M) = \<\rho_L(M)\ee_\gamma, \ee_\beta\>,$$
where $M \in \lieMP_2(\ZZZ)$.  There is a formula for the value of the
$\beta,\gamma$--coefficient of the Weil representation for a general
element of $\lieMP_2(\ZZZ)$.

\begin{theorem}\hskip-.5ex\emph{(Shintani \cite{shintani_weil})}
\label{shintani}
Let $\beta,\gamma \in L^*/L$ and
$$M = \(\[\begin{array}{cc}a&b\\c&d\end{array}\],\sqrt{c\tau+d}\)$$
be an element of $\lieMP_2(\ZZZ)$.  The $\beta,\gamma$--coefficient of
the Weil representation is given by
$$\delta_{\beta, a\gamma}\sqrt{i}^{\ \sgn(L)(\sgn(d)-1)} \e\({1\over 2}
ab\beta^2\)$$
if $c=0$, and by
$${\sqrt{i}^{\ -\sgn(L)\sgn(c)} \over |c|^{\dim(L)/2} \sqrt{|L^*/L|}}
\sum_{r\in L/cL} \e\({a(\beta+r)^2 - 2(\beta+r, \gamma) + d\gamma^2
\over 2c}\)$$
if $c\ne 0$.
\end{theorem}
Note, in the above formul\ae\ we are taking the square root of $i$
according to the branch chosen by the element of $\lieMP_2(\ZZZ)$.

It can be helpful to know where this theta function comes from.  Let
$W$ be a symplectic vector space and $V$ an orthogonal vector space.
Then, in a natural way, the space $W\tensor V$ is a symplectic space.
There is a theta function attached to a lattice contained in the space
$W\tensor V$.  Pick $V$ to be $L\tensor\RRR$ and $W$ to be $\RRR^2$
with the symplectic form given by the matrix
$\smallmatrix{0}{1}{-1}{0}$ and let the lattice be $\ZZZ^2 \tensor L$.
The theta function obtained for this lattice (when restricted to the
upper half plane and Grassmannian, which involves quotienting out by
the maximal compact subgroup of the symplectic group of $V\tensor W$)
is the Siegel theta function.  There is a natural representation of
the double cover of the symplectic group on functions on the
symplectic space.  In this case the double cover is the metaplectic
group and the natural representation is the Weil representation
\cite{howe_theta, weil}.

\section{Maass--Poincar\' e Series}

In this section we construct some real analytic modular forms on the
upper half plane; they will have explicitly known singularities at
$i\infty$.  We will later use their images under the singular theta
correspondence (which will then have known singularities) to attempt
to build certain piecewise linear automorphic functions.

A general way to construct (real analytic) modular forms is by an
averaging process: Pick a function invariant under the stabilizer of
the cusp at infinity that is also an eigenfunction of the Laplacian
and average over the remainder of $\lieMP_2(\ZZZ)$.  Provided that
everything converges, we obtain a real analytic eigenfunction of the
Laplacian with known singularities at $i\infty$.  Calculations similar
to this can be found throughout the literature (see, for example,
\cite{bruinier_borcherds, hejhal_PSL1})

If we have a modular form $F$ with type $\rho_L$ for some even lattice
$L$ then we know that $F(\tau+1) = \rho_L(T) F(\tau)$.  Hence, the
$\ee_\gamma$--component of $F$ satisfies
$$F_\gamma(\tau+1) = \e(\gamma^2/2)F_\gamma(\tau).$$
Thus, we can expand $F_\gamma$ as a Fourier series in the
$x$--variable
$$F_\gamma(\tau) = \sum_{n=-\infty}^\infty \e((n+\gamma^2/2)x)
c_{\gamma,n+\gamma^2/2}(y), \quad\hbox{where $\tau = x+iy$}.$$

Since the Fourier series converges uniformly we may apply the
Laplacian to the above expansion and interchange it with the
summation.  Doing this we obtain a differential equation that the
terms in the Fourier expansion satisfy.

\begin{lemma} $F$ satisfies the the differential equation
$\Delta F = \lambda F$ if and only if the functions $c_{\gamma,n}(y)$
satisfy the differential equations
$$c_{\gamma,n}''(y) + \[{b^++b^- \over y}\]c_{\gamma,n}'(y) +
\[{\lambda \over y^2} + {2\pi n(b^+-b^-) \over y} -
4\pi^2n^2\]c_{\gamma,n}(y) = 0.$$
\end{lemma}
\proof A simple calculation as indicated above.\qed

In order to simplify various expressions, it is convenient to define
$$\beta^+ = {b^+ + b^- \over 2} \qquad\hbox{and}\qquad \beta^- = {b^+ -
b^- \over 2}$$
and write the eigenvalue $\lambda$ in the form
$$\lambda = s(1-s) - \beta^+(1-\beta^+).$$

In fact, the number $s(1-s)$ is more fundamental than $\lambda$ --- if
we take our eigenfunction $F$ and consider $yF$ then it is easy to see
this is again a modular form but with weight $(b^+-1, b^--1)$.  Under
the Laplacian for this weight, the function $yF$ remains an
eigenfunction, however its new eigenvalue is $s(1-s) -
{\beta^+}'(1-{\beta^+}')$.  Here ${\beta^+}'$ is the value for the new
weight.  The reason for this is that the function induced on the group
$\lieMP_2(\RRR)$ is unchanged if we multiply by powers of $y$.  Its
eigenvalue under the Casimir element for $\lieMP_2(\RRR)$ is $s(1-s)$.
The extra factors of $\beta^+(1-\beta^+)$ come from the fact that we
regard our function on $\lieMP_2(\RRR)$ as a form on the upper half
plane with a certain weight.  Because of this we may also refer to
$s(1-s)$ as the eigenvalue.

\begin{lemma} When $n\ne 0$ a basis for the solutions to the
differential equation is given by
\begin{eqnarray*}
\M_{n,s}(y) & = & y^{-\beta^+}M_{\sgn(n)\beta^-,\ s-1/2}(4\pi|n|y) \\
\W_{n,s}(y) & = & y^{-\beta^+}W_{\sgn(n)\beta^-,\ s-1/2}(4\pi|n|y),
\end{eqnarray*}
where the functions $M$, $W$ are Whittaker's functions
\cite{abramowitz, erdelyi_I}.
\end{lemma}
\proof Let $d_n(y) = y^{\beta^+}c_{\gamma, n}(y)$.  We see that
$d_n(y)$ satisfies
$$d_n''(y) + \[{s(1-s) \over y^2} + {4\pi n\beta^- \over y} -4\pi^2 n^2
\]d_n(y) = 0.$$
If we substitute $z=4\pi |n|y$ we obtain the usual Whittaker
differential equation
$$d_n''(z) + \[{s(1-s) \over z^2} + {\sgn(n)\beta^- \over z} - {1\over
4}\]d_n(z) = 0.$$
The result now follows from the definition of Whittaker's functions.\qed

The solutions for $n=0$ are easily found to be (when $s\ne 1/2$)
$$y^{s-\beta^+}\quad\hbox{and}\quad y^{1-s-\beta^+}$$
as the differential equation is homogeneous.  Combining all of this we
see that the real analytic modular forms have a Fourier expansion of
the form (when $s\ne 1/2$)
$$F_\gamma(\tau) = A y^{s-\beta^+} + B y^{1-s-\beta^+} +
\sum_{n\in \ZZZ+\gamma^2/2} \[a_n \e(nx) \M_{n,s}(y) + b_n \e(nx)
\W_{n,s}(y)\].$$
When $s= 1/2$ we need to introduce a log term: $\log(y) y^{1/2 -
\beta^+}$.

The asymptotics for $\M$ and $\W$ are easily found from the
asymptotics for the Whittaker functions \cite{abramowitz}.

\begin{lemma}\label{asymptotics} The asymptotics of the functions $\W$
and $\M$ are:
\begin{enumerate}\parskip 0pt
\item As $y\rightarrow \infty$
\begin{enumerate}\parskip 0pt
\renewcommand{\labelenumii}{\roman{enumii}.}
\item $\W_{n,s}(y)$ decreases exponentially,
\item $\M_{n,s}(y)$ grows like $e^{2\pi |n| y}
y^{-\beta^+-\sgn(n)\beta^-}$.
\end{enumerate}
\item As $y\rightarrow 0$
\begin{enumerate}\parskip 0pt
\renewcommand{\labelenumii}{\roman{enumii}.}
\item $\W_{n,s}(y)$ is $O(y^{1-t-\beta^+})$,
\item $\M_{n,s}(y)$ is $O(y^{t-\beta^+})$.
\end{enumerate}
\end{enumerate}
\end{lemma}

Define the {\it Petersson slash operator} for the lattice $L$ to act
on vector--valued functions on the upper half plane by
$$(f|\gamma)(\tau) = \rho_L(\gamma)^{-1}f(\gamma\tau)
\sqrt{c\tau+d}^{\ 2b^+}\sqrt{c\bar\tau + d}^{\ 2b^-},$$
where $\gamma\in \lieMP_2(\ZZZ)$ and the lattice $L$ has signature
$(b^+,b^-)$.  As before, to ease notation, we are assuming that the
lattice will be obvious from the context.

Note that a function $F$ is a vector--valued modular form for the
lattice $L$ with type $\rho_L$ if and only if $F|\gamma = F$ for all
$\gamma\in \lieMP_2(\ZZZ)$.

\begin{lemma} The slash operator commutes with the Laplacian.
\end{lemma}
\proof $\lieMP_2(\ZZZ)$ acting on the upper half plane preserves the
metric which, in turn, determines the Laplacian.  Alternatively, one
can check this for the generators $S$ and $T$.\qed

Define the Maass--Poincar\'e series by
$$F_L(\tau; \delta, s) = {1 \over 2}\sum_{\gamma\in \<T\>\backslash
\lieMP_2(\ZZZ)} \left.\(\e(x\delta^2/2) \M_{\delta^2/2,s}(y)
\ee_\delta\)\right|\gamma.$$
It is easy to see that the only $\gamma$ in the above sum preserving
the cusp $i\infty$ are
$$\(\[\begin{array}{cc}1&0\\0&1\end{array}\],\pm 1\), \qquad
\(\[\begin{array}{cc}-1&0\\0&-1\end{array}\],\pm i\).$$
Looking at how these act on the summand, we can change the
definition of the Maass--Poincar\'e series to
$$F_L(\tau; \delta, s) = \sum_{\gamma\in \<T,Z\>\backslash
\lieMP_2(\ZZZ)} \left.\(\e(x\delta^2/2) \M_{\delta^2/2,s}(y)
\ee_\delta + \e(x\delta^2/2) \M_{\delta^2/2,s}(y)
\ee_{-\delta}\)\right|\gamma.$$

\begin{lemma}\label{convergence} The series converges uniformly for
$\Re(s) > 1+\beta^+$.  Consequently, the resulting function $F_L$ is
a real analytic modular form for the lattice $L$ with type $\rho_L$
and eigenvalue
$$\lambda = s(s-1)- \beta^+(\beta^+-1).$$
\end{lemma}
\proof By the asymptotics for $\M$ (Lemma \ref{asymptotics}), we can
bound the sum by a normal Eisenstein series with exponent $s-\beta^+$.
These are known to converge uniformly for $\Re(s-\beta^+) > 1$.  The
result now follows.\qed

\begin{lemma}\label{singularities}
$$F_L(\tau; \delta,s) - \e(x\delta^2/2)\M_{\delta^2/2,s}(y)\ee_\delta 
- \e(x\delta^2/2)\M_{\delta^2/2,s}(y)\ee_{-\delta}$$
tends rapidly to $0$ as $\tau\rightarrow i\infty$.
\end{lemma}
\proof This follows from the properties of the Eisenstein series used
in Lemma \ref{convergence}. \qed

\begin{theorem}\label{maas} The function $F_L(\tau; \delta, s)$, defined for
$\Re(s)>1+\beta^+$, is a real analytic modular form of weight
$(b^+,b^-)$, type $\rho_L$ and eigenvalue $s(s-1) -
\beta^+(\beta^+-1)$.  Its only singularity is at the cusp $i\infty$
and is of the form
$$\e(x\delta^2/2)\M_{\delta^2/2,s}(y)\ee_\delta +
\e(x\delta^2/2)\M_{\delta^2/2,s}(y)\ee_{-\delta}.$$
\end{theorem}
\proof This follows from Lemmas \ref{convergence} and
\ref{singularities}. \qed

These series (and simpler versions) occur in \cite{hejhal_PSL1,
hejhal_PSL2} and in these books it is shown that they can be
meromorphically continued to other values of $s$.  Fischer's book
\cite{fischer_selberg} contains a very clear account of these
functions and their uses.

\begin{theorem}\hskip-.5ex\emph{(Hejhal \cite{hejhal_PSL2})}
The function $F_L(\tau; \delta, s)$ can be analytically continued to a
meromorphic function of $s\in \CCC$.  The poles occur only in
the critical strip $0\le \Re(s)\le 1$.
\end{theorem}

\section{The Singular Theta Correspondence}

Now that we have a stock of real analytic modular forms with known
singularities, we study the functions that are obtained by applying
the singular theta correspondence.  In the holomorphic case, Borcherds
showed \cite{borcherds_grass} that we obtain functions with
singularities in the finite plane.  The singularities occur along
subspaces orthogonal to certain negative norm vectors.  In the real
analytic case we will see that this still happens but that it is also
possible to get singularities that are parallel to certain positive
norm vectors.

\begin{figure}[p]
\begin{centre}
\includegraphics[scale=1]{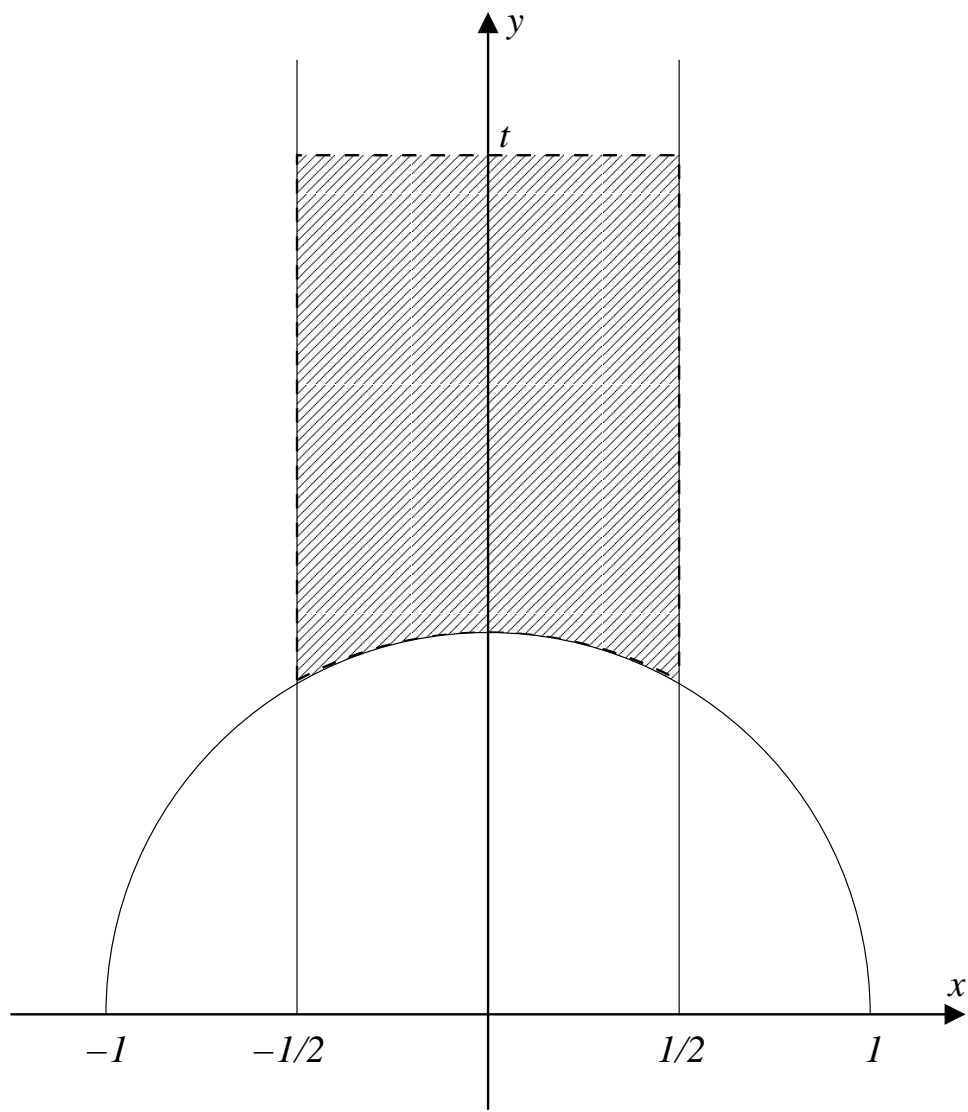}
\end{centre}
\caption{The truncated domain $\F_t$.}\label{trunc_domain}
\end{figure}

Firstly, we recall what the singular theta correspondence is (see
\cite{borcherds_grass, harvey_moore} for more details).  Suppose we have a
vector--valued modular form $F$ and an even lattice $L$ of signature
$(b^+,b^-)$.  The Siegel theta function $\Theta_L(\tau;v)$ of the
lattice $L$ is a modular form of weight $(b^+/2,b^-/2)$ and type
$\rho_L$.  Assume that $F$ has weight $(-b^-/2,-b^+/2)$ and type
$\rho_L$, then the product
$$F(\tau)\overline{\Theta_L(\tau;v)}$$
is a modular form of weight $0$.  If this product is of sufficiently
rapid decay at $i\infty$ (which occurs if $F$ is a cusp form) we can
take the integral
$$\Phi_L(v) = \int_\F F(\tau)\overline{\Theta_L(\tau;v)} {dxdy \over
y^2},$$
where $\F$ is the usual fundamental domain for $\lieSL_2(\ZZZ)$.  This
gives us a function $\Phi_L$ on $\Gr(L)$ invariant under a congruence
subgroup of $\lieO(L)$.  The map $F(\tau) \mapsto \Phi_F(v)$ is the
original theta (or Howe) correspondence.  In order to generalize this
construction to include non--cusp forms we have to find some
distributional extension of the map $F(\tau) \mapsto \Phi_F(v)$.  This
was done in \cite{harvey_moore} and \cite{borcherds_grass}.  The idea
is to truncate the integration domain in such a way that most of the
wildly non--convergent terms vanish.  The remaining non--convergent
terms are of polynomial growth and can be dealt with easily.  The
truncated domains are
$$\F_t = \{ \tau\in \F : \Im(\tau) \le t\},$$
see Figure \ref{trunc_domain} for a picture.
The {\it regularized} value of this integral is defined as the value
at $r=0$ of the analytic continuation of
$$\lim_{t\rightarrow \infty} \int_{\F_t}
F(\tau)\overline{\Theta_L(\tau;v)} {dxdy \over y^{2+r}}.$$
The above integral converges for $r$ sufficiently large 
\cite{borcherds_grass}.  This new, more general, map $F(\tau) \mapsto
\Phi_F(v)$ is the {\it singular theta correspondence}.

By performing the integrals that occur in the singular theta
correspondence it is easy to work out what kind of singularities the
function $\Phi_F(v)$ has.  The singularities occur on
sub-Grassmannians of the form $\gamma^\perp$, for $\gamma \in L^*$,
$\gamma^2<0$ where there is a non--zero coefficient corresponding to
$\gamma$ in $F$.  We now examine where the singularities are when we
apply the singular theta correspondence to the Maass--Poincar\'e functions.

\eject
\begin{lemma} If $F(\tau)$ is a real analytic modular form for the
lattice $L$ with singularity of the form
$$\e(x\delta^2/2)\M_{\delta^2/2,s}(y)\ee_\delta,$$
then the singularities of the regularized integral
$$\int_\F F(\tau)\overline{\Theta_L(\tau;v)} dx dy/y^2$$
occur only for $v\in\Gr(L)$ for which there is a vector in $L+\delta$
of norm $\delta^2$ which is in $v^+ \sqcup v^-$.
\end{lemma}
\proof Because we are looking for singularities of an integral we
can replace the integration domain $\F$ by $\{\tau\in\F : \Im(\tau)
\ge 1\}$ since the integral over a compact region will not affect
the singular behaviour of the integral.  By assumption, $F$ has only
one bad term in its Fourier expansion and we assume that this term is
of the form
$$e^{2\pi inx} \(y^m + O(1)\).$$
We can deduce the general case from this by writing the function $F$
as a sum of functions of this form.

Inserting the sum defining $\Theta_L$ and the expansion of $F$ we obtain
$$\int_{x=-1/2}^{1/2}\int_{y=1}^{\infty}\sum_{\lambda\in L+\gamma}
e^{\pi ix(\lambda^2-2n)} e^{-\pi y\(\lambda_{v^+}^2 - \lambda_{v^-}^2\)}
\( y^m +O(1) \) dxdy/y^2.$$
Due to the way the integral is regularized the $x$-integral will
remove all the terms except those with $\lambda^2 - 2n = 0$.  The
singularities therefore occur only where there is a vector
$\lambda \in L+\gamma$ with norm $2n$.  This leaves
$$\int_{y=1}^\infty \sum_{\lambda\in L+\gamma} e^{\pi y\(2|n|-
\lambda_{v^+}^2+\lambda_{v^-}^2\)}\(y^m +O(1) \) dy/y^2.$$
It is easy to see that the exponent is either $2\pi y\lambda_{v^-}^2$
or $-2\pi y\lambda_{v^+}^2$, depending on whether $n$ is positive or
negative.  So, if the exponent is non--zero then it is negative and
hence the integrand is rapidly decreasing as $y \rightarrow \infty$.
Thus the integral will converge and there will be no singularity.  So,
singularities occur only when $\lambda$ is such that the exponent is
zero --- this is exactly when $\lambda\in v^+ \sqcup v^-$.\qed

\begin{figure}[p]
\begin{centre}
\includegraphics[scale=.5]{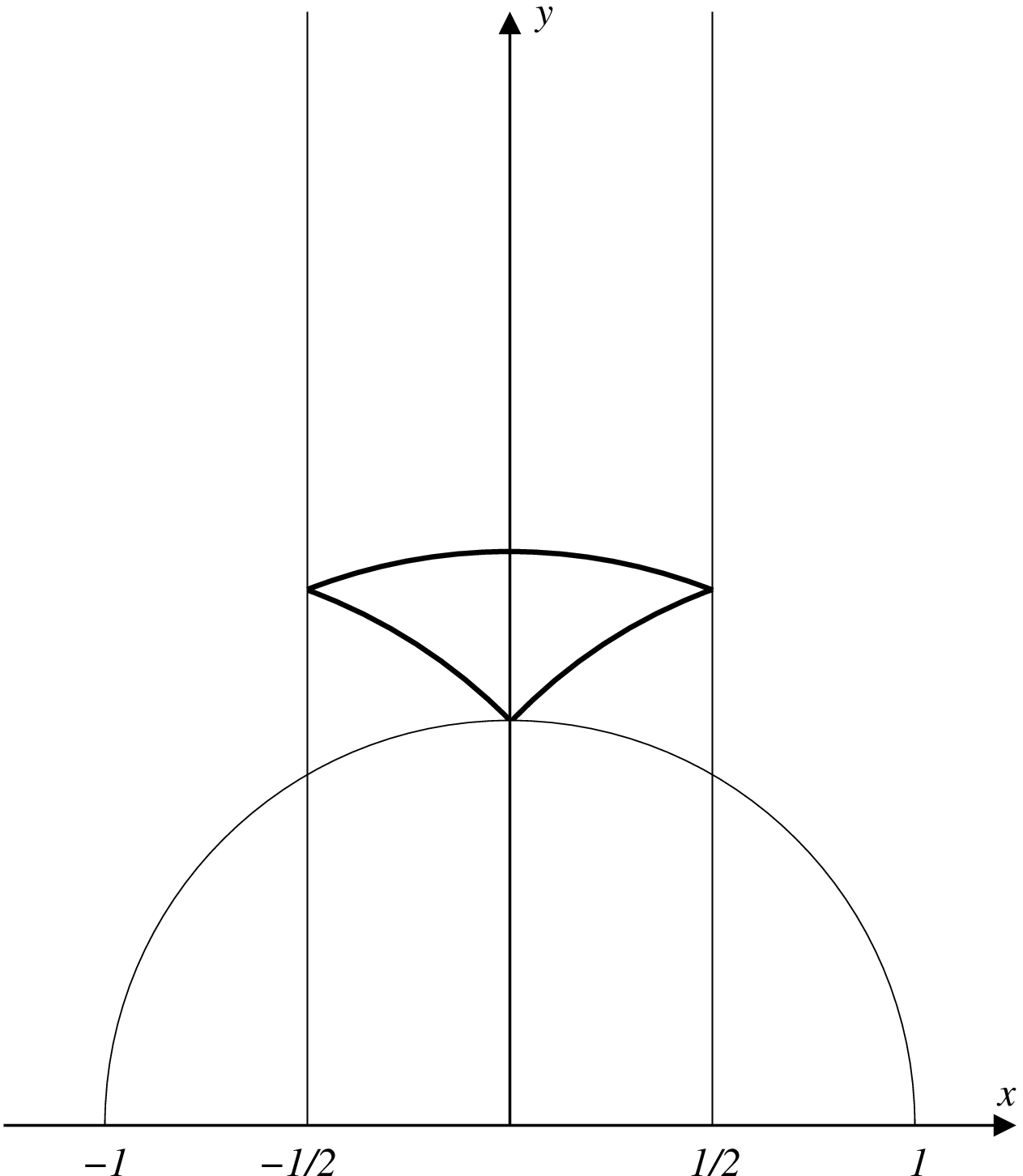}\qquad
\includegraphics[scale=.5]{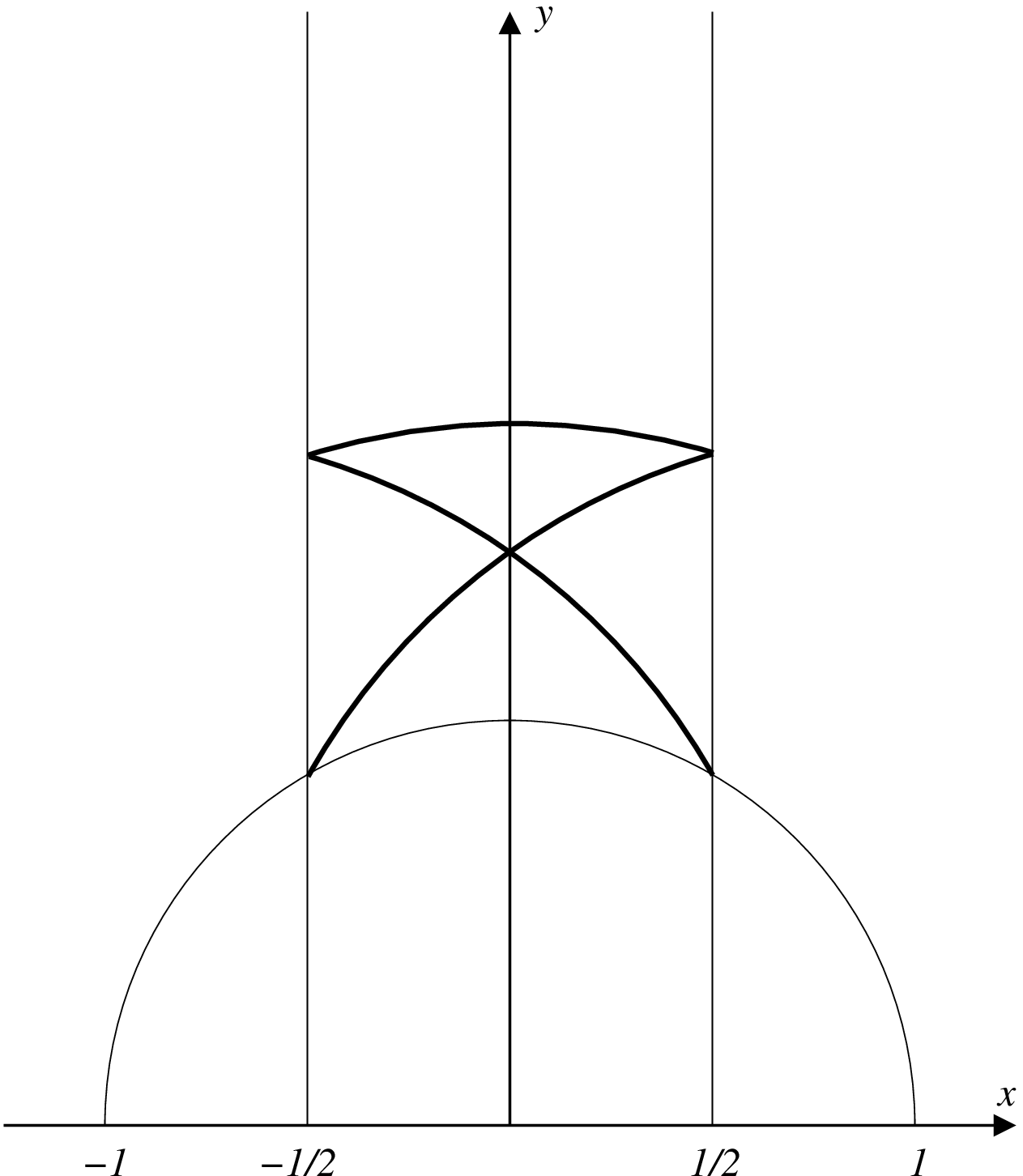}
\end{centre}
\vskip .5cm
\begin{centre}
\includegraphics[scale=.5]{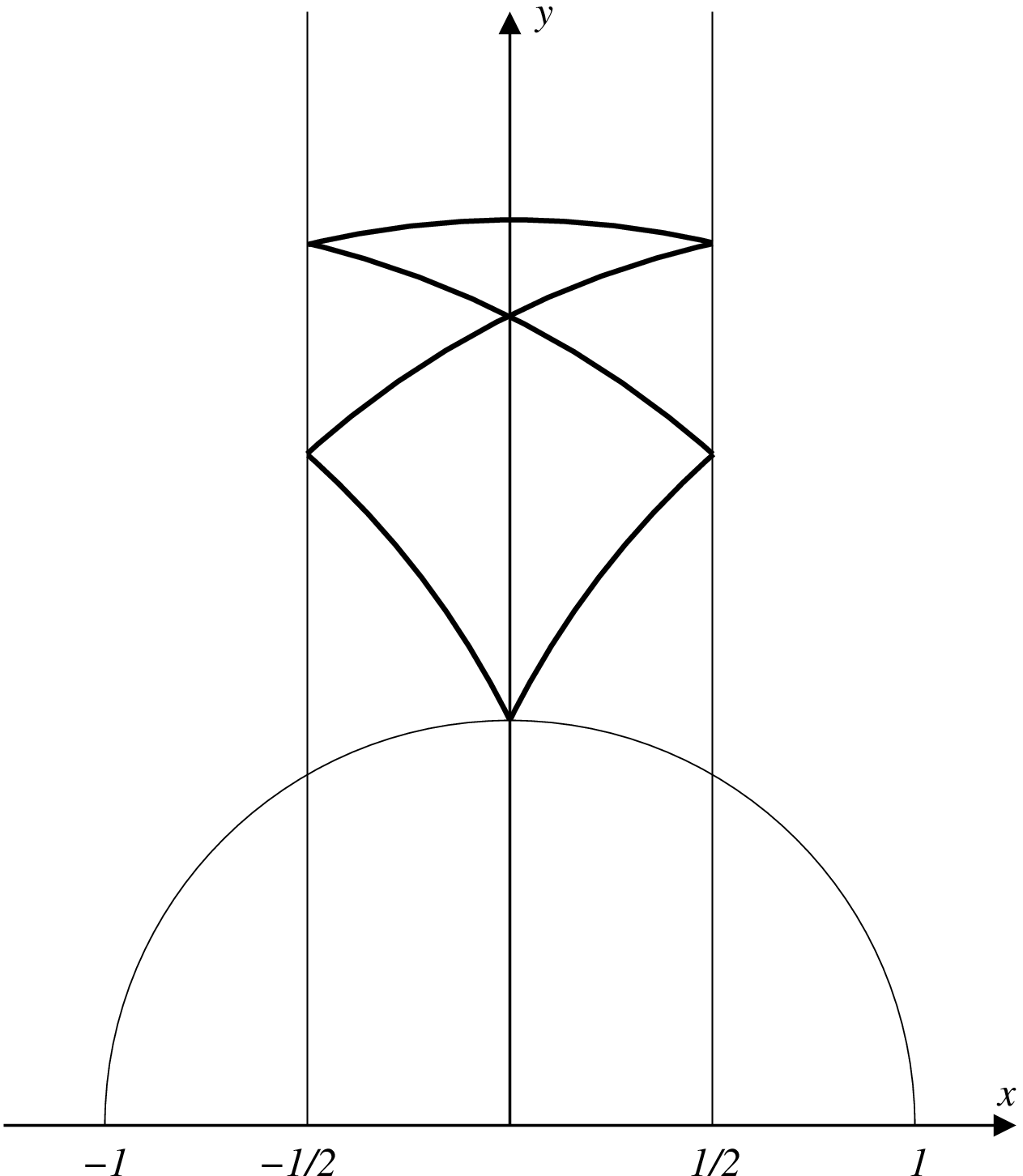}\qquad
\includegraphics[scale=.5]{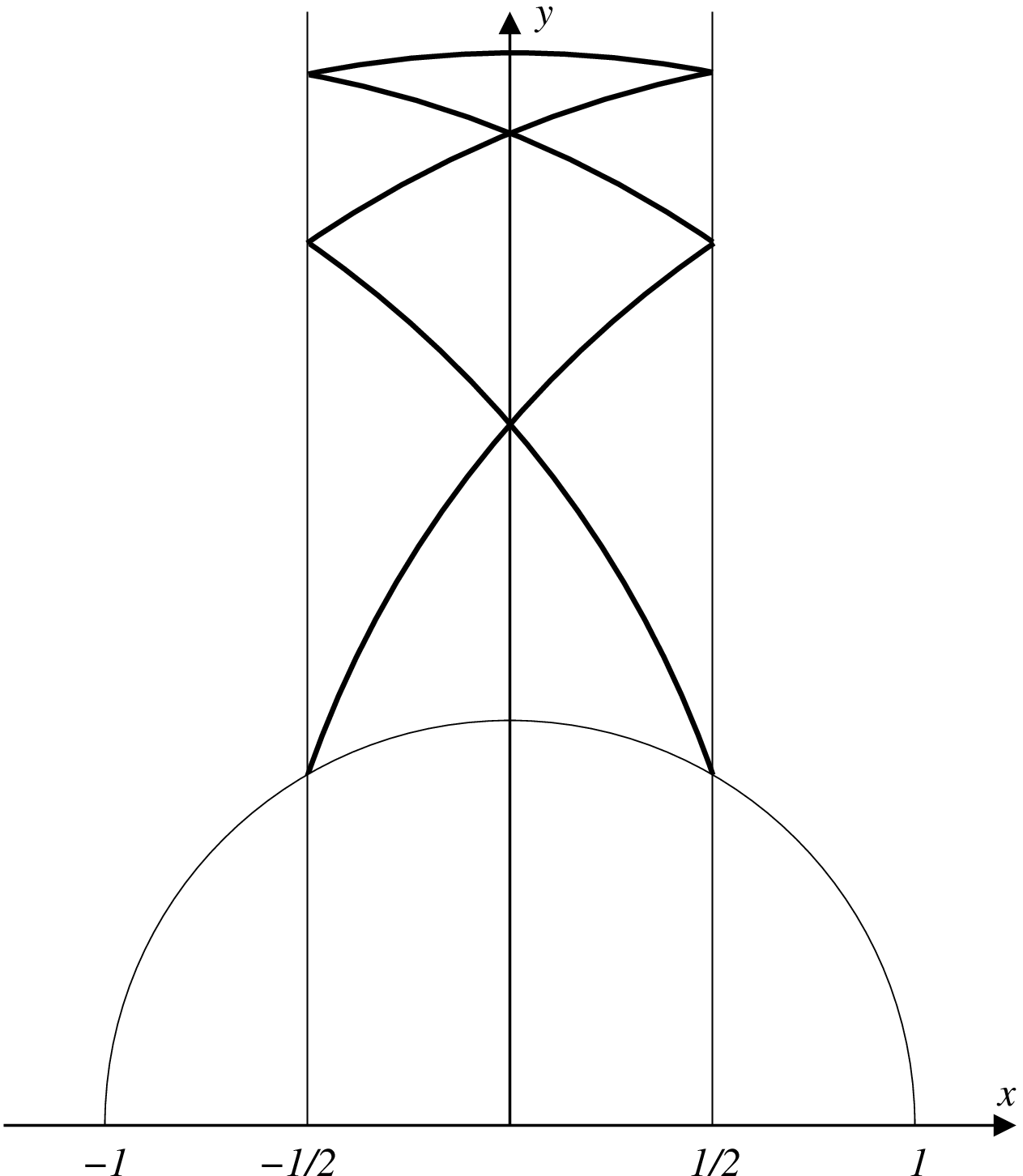}
\end{centre}
\caption{Singularities from real analytic forms.}\label{new_sings}
\end{figure}

Pictures of where these singularities occur in the fundamental domain
$\F$ for $\lieSL_2(\ZZZ)$ are given in Figure \ref{new_sings}.  The form of
the singularity that occurs along these sub--Grassmannians can be
computed using methods in Borcherds' paper \cite{borcherds_grass}.
Note, we have not yet found an interesting use for the new types of
singularity available when dealing with real analytic forms.  However,
the subspaces along which these singularities occur are certainly
interesting.  They appear, for instance, in \cite{zagier_zeta} as
curves associated to real quadratic fields of positive discriminant.

\chapter{Pseudo-cusp Forms}

\section{Introduction}

In his 1977 doctoral dissertation, Hartmut Haas investigated real
analytic modular forms for $\lieSL_2(\ZZZ)$.  In particular, he ran
several computer experiments to find the eigenvalues these real
analytic forms had under the invariant Laplacian for the upper half
plane.  His results were never published; however Harold Stark
observed that some of the numbers that occurring in his tables looked
like zeros of the Riemann zeta function.  With a little library
research, Dennis Hejhal noticed that many of the numbers from the
tables were coming from zeros of $\zeta(s) L(s,\chi_{-3})$, where
$$\chi_{-3} (n) = \( {-3 \over n}\).$$
Of course, $\zeta(s) L(s,\chi_{-3})$ is exactly the Dedekind zeta
function of $\QQQ(\sqrt{-3})$.

If it were true that each zero of the Riemann zeta function was
associated to a real analytic cusp form then one would be able to
prove the Riemann hypothesis (this follows from the fact that
self--adjoint operators have real eigenvalues).  This method for
proving the Riemann hypothesis had been independently suggested by
Polya and Hilbert (although it is fairly clear that they were not
expecting a group as simple as $\lieSL_2(\RRR)$).

Hopes of proving the Riemann hypothesis using these cusp forms were
dashed by Hejhal who noticed that these eigenfunctions were not really
eigenfunctions.  Due to a small oversight in Haas' computations the
``eigenfunctions'' could have logarithmic singularities at the point
$$\rho = {1\over 2} + {\sqrt{3} \over 2} i$$
and its conjugates under $\lieSL_2(\ZZZ)$.  So, unfortunately, they
were not really solutions to the spectral problem and therefore did
not prove the Riemann hypothesis.  The functions corresponding to
these spurious eigenvalues are called {\it pseudo--cusp forms} and
they have the following properties:

\begin{definition}\hskip-.5ex\emph{(Hejhal \cite{hejhal_zeta})}
Given $\lambda\in \CCC$.  We say that $f(z)$ is a
pseudo-cusp form [corresponding to $\lambda$] if and only if
\begin{enumerate}\parskip 0pt
\item $f(z)$ is automorphic with respect to $\Gamma$;
\item $f(z)$ is smooth at all points other than conjugates of $\rho$;
\item $\Delta f = -\lambda f$ at all points other than conjugates of
$\rho$;
\item $f(z)$ is cuspidal;
\item $f(z)$ has a logarithmic singularity at conjugates of $\rho$,
\ie
$$f(z) = \alpha \ln|z-\rho| + O(1)$$
near $z=\rho$ for some $\alpha \ne 0$.
\end{enumerate}
\end{definition}

Hejhal identified the spurious numbers in Haas' tables and gave the
following criterion for the existence of pseudo--cusp forms:

\begin{theorem}\hskip-.5ex\emph{(Hejhal \cite{hejhal_zeta})}
\label{hejhal_thm}
Suppose that $\lambda = s(1-s)$ with $\Re(s)
\ge {1\over 2}$.  There exists a pseudo-cusp form corresponding to
$\lambda$ if and only if
\begin{enumerate}\parskip 0pt
\item $f(\rho) =0$ for every [even] cusp form with eigenvalue
$\lambda$;
\item $\zeta(s) L(s,\chi_{-3}) = 0$.
\end{enumerate}
\end{theorem}

In this section we will see that the existence of pseudo--cusp forms
can be explained using the singular theta correspondence.  In
particular, we will see that pseudo--cusp forms occurring at zeros of
$\zeta(s)$ exist due to properties of the singular theta
correspondence.  The ones occurring at the zeros of $L(s,\chi_{-3})$
are a little more tricky --- they exist due to the vanishing of a
certain coefficient of an Eisenstein series.

\section{The Set Up}

After the calculation in this section were performed the author
discovered similar calculations in the work of Bruinier (in
particular, \cite{bruinier_borcherds}).  As the notation used is
slightly different and converting would probably introduce numerous
errors we choose to repeat our calculations here.  This also allows
our exposition to focus on the case at hand, rather than using the
more general formul\ae\ that occur in \cite{bruinier_borcherds}.

Let $L$ be the lattice $\ZZZ(2) \oplus \smallmatrix{0}{1}{1}{0}$.  We
represent an element of the lattice $L$ in coordinates as $(x,m,n)$,
where the squared norm is given by
$$(x,m,n)^2 = 2x^2 + 2mn.$$
It is easy to see that $|L^*/L| = 2$.  Denote the elements of the
discriminant form by $\ee_0$, $\ee_1$ which have norms $0$, $1/2$
respectively.  In coordinates we may assume
\begin{eqnarray*}
\ee_0 &=& (0,0,0), \\
\ee_1 &=& (1/2,0,0).
\end{eqnarray*}
Pick vectors $z=(0,0,1)$ and $z'=(0,1,0)$ and define $K = (L \cap
z^\perp) / \ZZZ z$.  $K$ is therefore isomorphic to $\ZZZ(2)$.

The Grassmannian for the lattice $L$ is the set of maximal (\ie $2$
dimensional) positive definite subspaces of $L\tensor\RRR$.

\begin{lemma}\label{grass} The Grassmannian for $L$ is isomorphic to
the upper half plane.
\end{lemma}
\proof An element of the Grassmannian corresponds to a $2$ dimensional
positive definite subspace $V$ of $L\tensor \RRR$.  Pick an
orthonormal basis $\{x,y\}$ for $V$ and define $z = x+iy \in
L\tensor\CCC$; it is clear that $z$ has norm $0$.  Similarly, any
norm $0$ vector in $L\tensor \CCC$ corresponds to an element of the
Grassmannian.  The plane represented is clearly unchanged under scalar
multiplication and complex conjugation on $L\tensor \CCC$.  So we can
uniquely represent an element of the Grassmannian by a vector $(\tau,
1, *) \in L\tensor \CCC$ with $\Im(\tau)>0$ and $*$ chosen to give the
vector norm $0$.  This gives the isomorphism between $\Gr(L)$ and the
upper half plane.\qed

\begin{lemma} There is an $\lieSL_2(\ZZZ)$ action on $L$
which descends to the usual M\"obius action on the upper half plane.
\end{lemma}
\proof Represent an element $(x,m,n)$ of $L$ as a matrix
$$\[ \begin{array}{cc} x & m \\ n & -x\end{array} \].$$
Note that the determinant of the matrix is $-x^2-mn = -{1\over 2}
(x,m,n)^2$.  Thus $L$ is isometric to the traceless $2\times 2$
integer matrices.  There is an $\lieSL_2(\ZZZ)$ action on the matrices
(by conjugation) which, therefore, gives an action on $L$ and hence on
$\Gr(L)$.  Examining how the two generators for $\lieSL_2(\ZZZ)$ act we
see that this is the usual M\" obius action on the upper half
plane. \qed

\begin{lemma}\label{points} The points congruent to $\rho$ under
$\lieSL_2(\ZZZ)$ are exactly the elements of $\Gr(L)$ orthogonal to
norm $-3/2$ vectors.
\end{lemma}
\proof Let $(a,b,c)$ be a norm $-3/2$ vector, \ie $2a^2 + 2bc =
-3/2$.  We know that elements of $\Gr(L)$ can be represented as
$(\tau, 1, -\tau^2) \in L\tensor \CCC$.  This will be orthogonal to
$(a,b,c)$ if and only if
$$2a\tau -b\tau^2 + c = 0.$$
The solutions to this are of the form
$$\tau = {2a \pm \sqrt{3}\ i \over 2b}.$$
Since $a,b\in \ZZZ$ the only solutions lying in the fundamental
domain $\F$ are $\rho$ and $\rho-1$.\qed

As the lattice $L$ has signature $(2,1)$ we know that $\beta^+ =
-3/4$, $\beta^- = 1/4$.  Let $\gamma$ be a vector of norm $-3/2$.
Consider the Maass--Poincar\'e function $F_L(\tau; \gamma, s)$.  Using
Lemma \ref{asymptotics} we see the asymptotics of $F_L$ as $y
\rightarrow \infty$
$$F_L(\tau; \lambda, s)\quad\sim\quad\exp(3\pi y /2) (y + O(1)).$$
From the singularity results of Borcherds \cite{borcherds_grass} we
see that the singular theta transformation of $F_L$ has logarithmic
singularities at points orthogonal to norm $-3/2$ vectors.  By the
previous lemma these are points congruent to $\rho$ under the action
of $\lieSL_2(\ZZZ)$.  The Fourier expansion is easily read from
\cite{borcherds_grass}
\begin{eqnarray*}
\Phi_L(u+iv) &=& v \Phi_K(u+iv) + 2v \sum_{\lambda \in \ZZZ(2)^*}
\sum_{n>0} e^{2\pi i n(\lambda,u)} \\
&\times& \int_{y>0} c_{\lambda, \lambda^2/2} (y) \exp(-\pi n^2 v^2 /
2 y - \pi y \lambda^2) y^{-5/2} dy.
\end{eqnarray*}
Note, we may have to regularize the $y$--integral by dividing by
$y^r$, analytically continuing and taking the limit as $r\rightarrow 0$.

\begin{lemma} For $\lambda\ne 0$, the integral
$$\int_{y>0} c_{\lambda, \lambda^2/2} (y) v \exp(-\pi n^2 v^2 / 2y - \pi
y \lambda^2) y^{-5/2} dy$$
decreases rapidly as $v \rightarrow \infty$.
\end{lemma}
\proof It is clear that the integral from $0$ to $1$ decreases rapidly
in $v$.  It suffices to consider only the integral over $y>1$.  Note
that
$$\pi n^2 v^2 / 2y + \pi \lambda^2 y \ge \sqrt{2}\pi n |\lambda| v,$$
so the integral is dominated by
$$v \exp(-\sqrt{2}\pi n |\lambda| v) \int_{y>1} c_{\lambda,
\lambda^2/2} (y) y^{-5/2} dy.$$
The integral converges since $c(y)$ is rapidly decreasing as
$y\rightarrow \infty$.  The required rapid decrease with $v$ is now
obvious.\qed

In particular, the above lemma shows that most of the terms in the
Fourier expansion are rapidly decreasing: the only terms that may not
be rapidly decreasing come from $\lambda=0$ and $\Phi_K$.

\begin{lemma}\label{zeroterm} Assume $s\ne 1/4$ and $s\ne 3/4$.  The
$\lambda = 0$ term
$$\sum_{n>0} \int_{y>0} c_{0,0}(y) \exp(-\pi n^2 v^2/2y) y^{-5/2}dy$$
is (modulo a constant factor)
$$(\pi/2)^{1/4-s} v^{1/2-2s} \zeta (2s - 1/2) \Gamma(s - 1/4).$$
\end{lemma}
\proof To simplify the calculation we ignore the sum (and remember to
put it back in at the end!).  We need to regularize the integral, so
consider instead
$$\int_{y>0} c_{0,0}(y) \exp(-\pi n^2 v^2/2y) y^{-r-5/2}dy.$$
Note that $c_{0,0}(y)$ is (modulo a constant factor) $y^{7/4 - s}$.
We calculate:
\begin{eqnarray*}
&& \int_{y>0} y^{-3/4 -r -s} \exp(-\pi n^2 v^2/2y)dy \\
&=& (\pi n^2 v^2/2)^{1/4 - r-s} \int_{y>0} y^{-5/4+r+s} e^{-y} dy \\
&=& (\pi n^2 v^2/2)^{1/4 -r-s}\Gamma(-1/4+r+s).
\end{eqnarray*}
Now, summing over $n>0$ gives
$$(\pi v^2/2)^{1/4-r-s} \zeta(2r+2s-1/2)\Gamma(-1/4+r+s).$$
Assume that this is non--singular at $r=0$.  The value at $r=0$ of the
analytic continuation is therefore
$$(\pi/2)^{1/4-s} v^{1/2-2s} \zeta(2s-1/2) \Gamma(-1/4+s).$$
The function $\pi^{-s/2} \Gamma(s/2) \zeta(s)$ is well known to have a
meromorphic continuation to all of $\CCC$ with simple poles only at
$s=0$ and $s=1$ (see, for example, \cite{bump_automorphic}).
Therefore the above calculation is valid except when
$2s-1/2 = 0$ or $2s-1/2 = 1$.\qed

In order to attempt to create a cusp form we will assume that $2s-1/2$
is a root of the zeta function.  Suppose this root is $\omega$, then
simple calculation shows that $3/4 - 4s(1-s) = \omega(\omega - 1)$.
Note that neither $0$ nor $1$ are roots of the Riemann zeta function,
therefore there will be no problems in using Lemma \ref{zeroterm}.

\begin{lemma} For $\lambda\ne 0$ the integral in the Fourier expansion
is an eigenfunction of the Laplacian
$$v^2\( {\d^2 \over \d u^2} + {\d^2 \over \d v^2} \)$$
with eigenvalue $6 - 4\lambda = 3/4 - 4s(1-s) = \omega(\omega-1)$.
\end{lemma}
\proof Notice that the following differential operators
$$v^2\({\d^2 \over \d u^2} + {\d^2 \over \d v^2}\), \quad
4y^2\( {\d^2 \over \d y^2} + {1 \over 2y}{\d \over \d y} + {\pi\lambda^2
\over 2y} - \pi^2\lambda^4\)$$
act identically on
$$v\exp(2\pi in(\lambda,u))\exp(-\pi n^2v^2/ 2y -\pi y\lambda^2).$$
We calculate by applying the $(u,v)$--Laplacian, using the above
identification to change it to an $(x,y)$--Laplacian.  Then,
integrating by parts to move it to act on the $c(y)$ term we see it
basically becomes the differential equation that $c(y)$ satisfies.

\begin{eqnarray*}
& & v^2(\d_u^2 + \d_v^2) \int_{y>0}v\exp(2\pi in(\lambda,u))\exp(-\pi n^2v^2/
2y-\pi y\lambda^2)c_{\lambda,\lambda^2/2}(y)y^{-5/2} dy \\
& = & \int_{y>0}\[4y^2(\d_y^2  + \d_y /2y + \pi\lambda^2/2y - \pi^2\lambda^4)
v\exp(2\pi in(\lambda,u))\right. \\
& & \left. \times \exp(-\pi n^2v^2/2y-\pi
y\lambda^2)\]c_{\lambda,\lambda^2/2}(y) y^{-5/2} dy \\
& = & 4\int_{y>0}v\exp(2\pi in(\lambda,u))\exp(-\pi n^2v^2/
2y-\pi y\lambda^2) \\
& & \times (\d_y^2  - \d_y /2y + 1/2y^2 + \pi\lambda^2/2y -
\pi^2\lambda^4) c_{\lambda,\lambda^2/2}(y)y^{-1/2} dy \\
& = & 4\int_{y>0}v\exp(2\pi in(\lambda,u))\exp(-\pi n^2v^2/
2y-\pi y\lambda^2)y^{-1/2} \\
& & \times (\d_y^2  - 3\d_y /2y + 3/2y^2 + \pi\lambda^2/2y -
\pi^2\lambda^4) c_{\lambda,\lambda^2/2}(y) dy \\
& = & 4\int_{y>0}v\exp(2\pi in(\lambda,u))\exp(-\pi n^2v^2/
2y-\pi y\lambda^2)y^{-1/2} \\
& & \times (-\kappa/y^2 - \pi\lambda^2/2y + \pi^2\lambda^4
+ 3/2y^2 + \pi\lambda^2/2y - \pi^2\lambda^4) c_{\lambda,\lambda^2/2}(y) dy \\
& = & 4\int_{y>0}v\exp(2\pi in(\lambda,u))\exp(-\pi n^2v^2/
2y-\pi y\lambda^2)y^{-1/2} \\
& & \times (-\kappa/y^2 + 3/2y^2) c_{\lambda,\lambda^2/2}(y) dy \\
& = & 4(3/2 - \kappa)\int_{y>0}v\exp(2\pi in(\lambda,u))\exp(-\pi n^2v^2/
2y-\pi y\lambda^2)c_{\lambda,\lambda^2/2}(y) y^{-5/2} dy.
\end{eqnarray*}
Hence, the integral terms are eigenfunctions with eigenvalue
$6-4\kappa$.\qed

A more conceptual reason why this should work can be found in
\cite{howe_theta}.  It is possible to show using Howe's results that
the singular theta transformation will ``almost'' satisfy all the
invariant differential operators for $\Gr(L)$ (the ``almost'' comes
from the fact that extra terms are generated by the regularization
procedure).  In fact, by working in the adelic setting we can get a
similar result for the Hecke operators.  This shows that the map from
modular forms to automorphic products commutes with an action of the
Hecke operators coprime to the level (where the action is additive on
modular forms and multiplicative on automorphic forms).  This possible
behaviour is mentioned in Question 10 of \cite{borcherds_products}.

\begin{lemma} The integral from the smaller lattice is
$$\#\{\gamma\in K+\delta : \gamma^2 = -3/2\}
{ \Gamma(2s) \Gamma(r) \Gamma(-1/4 + s - r) \over \Gamma(s-1/4)
\Gamma(1/4 + s +r)} (3\pi)^{-1/4 -r}.$$
\end{lemma}
\proof In the region of convergence we have
$$\int_{\lieSL_2(\ZZZ)/\H} \overline{\Theta_K(\tau;v)} \sum_{\gamma
\in \<S\>\backslash \lieMP_2(\ZZZ)}\left.\(\e(x\delta^2/2)
\M_{\delta^2/2,s}(y) \ee_\delta\)\right|\gamma\ {dxdy \over
y^{2+r}}.$$
We denote by $G(\tau)$ the function that is in the Maass--Poincar\'e
sum.  Using the transformation properties of $\Theta_K$ under
$\lieMP_2(\ZZZ)$ we can re-write this as
$$\int_{\lieSL_2(\ZZZ)/\H} \sum_{\gamma \in
\<S\>\backslash\lieMP_2(\ZZZ)}
\overline{\Theta_K(\gamma\tau; v)} G(\gamma\tau) {dxdy \over y^{2+r}}.$$
We would like to say that this is equal to
$$\int_{y>0} \int_{x\in \RRR/\ZZZ} \overline{\Theta_K(\tau;v)} G(\tau)
{dxdy \over y^{2+r}}.$$
However, we have to be careful due to the regularization of the two
integrals.  We copy the method used in \cite{borcherds_grass}.  In the
region $y>1$ the two integrals clearly have the same form of
divergence and hence are identically regularized.  So, if we can show
that the second integral converges absolutely in the region $y<1$ then
we can justify the exchange of sum and integral.  But, this is clearly
true from the behaviour of $\Theta$ and $G$ in the region $y<1$.  We
are allowed to replace $\Im(\tau)^{-s}$ and $\Im(\gamma\tau)^{-s}$ as
they both agree at $s=0$.  Hence, we can make the deduction about the
integrals.  Performing the $x$--integral leaves
$$\#\{\gamma\in K+\delta : \gamma^2 = -3/2\}
\int_{y>0} e^{-3\pi y/2}
\M_{-3/4, s}(y) y^{-2-r} dy.$$
In terms of the Whittaker functions this is
$$\#\{\gamma\in K+\delta : \gamma^2 = -3/2\}
\int_{y>0} y^{-5/4-r}
e^{-3\pi y/2} M_{-1/4, s-1/2}(3\pi y) dy.$$
This integral can be found in standard tables of integrals (see, for
example, \cite{gradshteyn})
$$\#\{\gamma\in K+\delta : \gamma^2 = -3/2\}
{ \Gamma(2s) \Gamma(r)
\Gamma(-1/4 + s - r) \over \Gamma(s-1/4) \Gamma(1/4 + s +r)}
(3\pi)^{-1/4 -r}.$$
This is the required result.\qed

In our case the term at the front is zero (because $K$ is positive
definite).

Note, that the above calculation can be done in general and it shows
that if a Maass--Poincar\'e series is put through the singular theta
correspondence for a lattice with no vectors corresponding to its
singularities then the resulting function will be identically zero.
This is unfortunate since if the resulting transformation were
non--zero then it would provide a method for constructing cusp forms
with eigenvalues corresponding to the zeros of the Riemann zeta
function.  This would then give a proof of the Riemann hypothesis.

Putting together everything from above we see that we have constructed
a function on the upper half plane with logarithmic singularities at
the conjugates of $\rho$.  This function is an eigenfunction for the
Laplacian with eigenvalue $\omega(\omega-1)$ where $\omega$ is a zero
of the Riemann zeta function.  Finally, the function is a Maass cusp
form.  In other words, we have constructed a Hejhal pseudo--cusp
form.  To summarize

\begin{theorem}\label{pseudo_Riemann} Let $\omega$ be a zero of the
Riemann zeta function.  Then there is a pseudo--cusp form with
eigenvalue $\omega(\omega-1)$.
\end{theorem}

\section{The Constant Term}

It is possible to get pseudo--cusp forms in other ways.  If the
``constant'' term of the Maass--Poincar\'e series vanishes then the
singular theta correspondence applied to the Maass--Poincar\'e series
will be an eigenfunction of the Laplacian (recall that the problems
in the previous section were all due to the constant term).  So,
we should look carefully at the constant term.

Assume that the only norm zero vector in $L^*/L$ is the zero vector,
the constant term of the Maass--Poincar\'e series is then
\cite{bruinier_borcherds}:
$${4^{1-\sgn(L)/2} \pi^{1+s-\sgn(L)/2} |m|^{s-\sgn(L)/2} \over (2s-1)
\Gamma(s+\sgn(L)/2) \Gamma(s-\sgn(L)/2)} \sum_{c\in \ZZZ-\{0\}}
|c|^{1-2s} H_c(\beta,m,0,0),$$
where $H_c(\beta,m,\gamma,n)$ is the {\it generalized Kloosterman sum}
$$H_c(\beta,m,\gamma,n) = {\exp(-\pi i \sgn(c) k/2) \over |c|}
\sum_{(c,d)=1} \rho^{-1}_{\beta\gamma} \(\[\begin{array}{cc}a&b \\
c&d\end{array}\]\) \e\({ma + nd \over c}\).$$
The matrix $\smallmatrix{a}{b}{c}{d}$ appearing in the above sum is
any element from $\lieSL_2(\ZZZ)$ having the required values for $c$
and $d$ (the Kloosterman sum is independent of the choice of $a$ and
$b$).  From the definition of $H_c$ it is clear that
$$H_c(\beta,m,0,0) = -H_{-c}^*(0,0,\beta,-m).$$
Substituting this into the formula for the constant term we see that
the constant term for the Maass--Poincar\'e series is the same as the
$(\beta,-m)$--coefficient of the Eisenstein series defined in
\cite{bruinier_kuss}.

Working through the details of Theorem 4.6 of \cite{bruinier_kuss} we
see that the constant term is therefore proportional to $L(\omega,
\chi_{-3})$.  Hence we deduce:

\begin{theorem} Let $\omega$ be a zero of $L(s,\chi_{-3})$, then there
is a pseudo--cusp form with eigenvalue $\omega(\omega-1)$.
\end{theorem}

Combining this with Theorem \ref{pseudo_Riemann} we get:

\begin{theorem} Let $\omega$ be a zero of the Dedekind eta function
for $\QQQ(\sqrt{-3})$, then there is a pseudo--cusp form with
eigenvalue $\omega(\omega-1)$.
\end{theorem}

Comparing this result with that of Hejhal (Theorem \ref{hejhal_thm}) 
we obtain:

\begin{corollary} Let $\omega$ be a zero of the Dedekind zeta function
for $\QQQ(\sqrt{-3})$ and $f(\tau)$ an even Maass cusp form with
eigenvalue $\omega(\omega-1)$.  Then $f(\rho) = 0$.
\end{corollary}

Of course, this is not surprising because it is almost certainly true
that there are no Maass cusp form with these eigenvalues!

\section{Hecke Triangle Groups}

We can use a similar construction to obtain pseudo--cusp forms for the
arithmetic Hecke triangle groups which have been studied by Hejhal
\cite{hejhal_hecke}.  We briefly discuss how to do this --- the
details are exactly the same as the previous section, the reader can
fill them in if they so desire.

The {\it Hecke triangle group} $\GGG_n$ is generated by the
transformations
\begin{equation}
\tau \longmapsto -{1\over \tau} \qquad\hbox{and}\qquad \tau
\longmapsto \tau+ 2\cos\({\pi \over n}\).
\label{heckeactions}
\end{equation}
So, $\GGG_3$ is $\liePSL_2(\ZZZ)$ (\ie $\lieSL_2(\ZZZ)$ acting on the
upper half plane by M\"obius transformations).  The other arithmetic
triangle groups are $\GGG_0$, $\GGG_4$, $\GGG_6$ and $\GGG_\infty$.
The cases of $4$ and $6$ are investigated in \cite{hejhal_hecke}.  In
this paper numerical evidence shows that there are pseudo--cusp forms
corresponding to zeros of $L(\omega, \chi_4)$ (for $\GGG_4$) and
$L(\omega, \chi_{-3})$ (for $\GGG_6$).  These pseudo--cusp forms can
be obtained by a singular theta correspondence in the same way as
before.

Let $L$ be the lattice $\ZZZ(2N) \oplus \smallmatrix{0}{1}{1}{0}$, for
$N = 1,2,3$ ($N=1$ corresponds to $\GGG_3$, $N=2$ to $\GGG_4$ and
$N=3$ to $\GGG_6$).

\begin{lemma} The Grassmannian $\Gr(L)$ is isomorphic to the upper
half plane.
\end{lemma}
\proof Essentially the same as Lemma \ref{grass} except this time we
represent elements of the Grassmannian in the form
$$(\tau/\sqrt{N} , 1 , -\tau^2),$$
where $\tau$ is in the upper half plane.\qed

\begin{lemma} The action of the Hecke triangle group on the upper half
plane comes from automorphisms of the lattice $L$.
\end{lemma}
\proof The automorphism
$$(x,m,n) \longrightarrow (x,n,m)$$
gives the action
$$\tau \longmapsto -{1\over \tau}.$$
The automorphism
$$(x,m,n) \longrightarrow (x+m,m,*)$$
(where $*$ is chosen to preserve the norm) gives the action
$$\tau \longmapsto \tau + \sqrt{N}.$$
For the cases $N=1,2,3$ this is the second transformation in
Equation (\ref{heckeactions}).\qed

Let $\delta$ be the vector $(1/2,1,-1) \in L^*$ which is of norm $N/2
-2$.  This is a representative of the unique vector of order exactly
$2$ in the discriminant group.

\begin{lemma} The points in the upper half plane orthogonal to
vectors of norm $N/2 - 2$ are exactly the points $\rho$ where
logarithmic singularities occur in \cite{hejhal_hecke}.
\end{lemma}
\proof The same as Lemma \ref{points}.\qed

The above construction then applies to give a pseudo--cusp form with
eigenvalue $\omega(\omega-1)$ for $\omega$ a zero of the Riemann zeta
function.  Working through the calculation of the constant term we get
the following:

\begin{theorem}\ 
\begin{enumerate}\parskip 0pt
\item Let $\omega$ be a root of the Dedekind zeta function
for $\QQQ(\sqrt{-2})$ then there is a pseudo--cusp form of eigenvalue
$\omega(\omega -1)$ for the Hecke triangle group $\GGG_4$.
\item Let $\omega$ be a root of the Dedekind zeta function
for $\QQQ(\sqrt{-3})$ then there is a pseudo--cusp form of eigenvalue
$\omega(\omega -1)$ for the Hecke triangle group $\GGG_6$.
\end{enumerate}
\end{theorem}

Note that the Dedekind zeta function for $\QQQ(\sqrt{-2})$ is
$\zeta(s)L(s,\chi_4)$ and for $\QQQ(\sqrt{-3})$ is
$\zeta(s)L(s,\chi_{-3})$.

Of course the above method could be used to create different types of
pseudo--cusp form with logarithmic singularities in different places.
Functions of this form are mentioned at the end of \cite{hejhal_zeta}.

\chapter{Converse Theorem for $O(1,n)$}

It was observed by Borcherds in \cite{borcherds_refl} that in many
cases ``interesting reflection groups of Lorentzian lattices are
controlled by certain modular forms with poles at cusps''; the
definition of interesting is vague (any reflection group that can be
associated to a modular form in a natural way should, by any sensible
definition, be included).  However, a very important property that a
reflection group can have is the existence of a Weyl vector.  Hence,
any definition of an interesting reflection group should include the
case a Weyl vector exists.  There are interesting lattices which do
not have Weyl vectors, for example, the odd unimodular lattices
$\odd_{1,n}$ for $n=20,21,22,23$ \cite{borcherds_20}.  It seems likely
that these lattices can be included in the framework presented below
provided we use theta functions of odd lattices.

The correspondence between modular forms with poles at cusps and
reflection groups uses the singular theta correspondence 
\cite{borcherds_grass}.  This associates to any vector--valued
modular form (of the correct weight and type) a piecewise linear
function on the hyperbolic space of $\RRR^{1,n}$.  The singularities
of the image are determined by the singularities of the modular form.
If these singularities occur along the reflection hyperplanes of a
reflection group then we say that this reflection group is associated
to the modular form.  It was noticed in \cite{borcherds_refl} that
many Lorentzian reflection groups were associated in this way to a
modular form.

It is natural to ask whether every nice Lorentzian reflection group is
indeed associated with a modular form (Problem 13.2 in
\cite{borcherds_refl}).  In this chapter we prove that such a
correspondence does exist provided the Lorentzian lattice has a
sufficiently large dimension and a sufficiently small $p$-rank.  Some
kind of condition on the lattice seems to be necessary due to the
annoying fact that scaling the norm on the lattice by an integer does
not change the properties of the reflection group but generally
prevents the existence of a corresponding modular form.  We will
discuss this in more detail later on in this chapter where we shall
give some arguments and examples to show that the result probably
still holds with fewer restrictions.

Similar results to this have been found in the case of lattices with
signature $(2,n)$ by Bruinier \cite{bruinier_borcherds}.  He
showed that (with technical conditions similar to those mentioned
above) any automorphic form on $\Gr(L)$ for $L$ with signature $(2,n)$
having all its zeros occurring orthogonal to certain lattice
vectors comes from the singular theta transformation of a holomorphic
modular form.  The first part of what follows can be regarded as the
equivalent theorem for lattices of signature $(1,n)$.  Bruinier
mentions in \cite{bruinier_borcherds} that it should be possible to
reduce the technical conditions.

We use the association between Lorentzian lattices and modular forms
to prove a folklore conjecture that the highest dimension in which a
Lorentzian lattice possesses a Weyl vector is $26$.  This upper bound
is sharp as there is a well known example in $26$ dimensions ---
$\even_{1,25}$.  We then show that any example in $26$ dimensions is a
sublattice of $\even_{1,25}$ with its norm scaled by some factor.

Nikulin has shown that there are only a finite number of elementary
Lorentzian lattices possessing Weyl vectors, although the best bound
on the maximal dimension was much higher than $26$.  Later on we will
use modular forms to give a new (heuristic) proof of this finiteness
result.  We hope that it is possible to make this heuristic proof
rigorous and that these methods can be used to classify all examples
of Lorentzian lattices with Weyl vectors.

\section{The Functions}

In this section we define the functions that will be used in the
construction of the map between interesting Lorentzian reflection
groups and modular forms.  The basic idea is to make the
Maass--Poincar\'e series from Theorem \ref{maas} look as holomorphic
as possible.

Consider the Maass--Poincar\'e series from Theorem \ref{maas}.  If we
specialize to $s(1-s)=0$ we get real analytic modular forms with
eigenvalues under the Laplacian equal to the eigenvalues that
holomorphic modular forms have.  So, some of the terms in the Fourier
expansion will look like holomorphic terms (\ie they will be of the
form $q^n$) and the others will contain some Whittaker function.  The
image of these Maass--Poincar\'e series under the singular theta
correspondence has known singularities and using results of Bruinier
\cite{bruinier_borcherds} we know the image can be decomposed
into a piecewise linear piece (this comes from the terms that look
holomorphic) and a real analytic piece (this comes from the terms that
look real analytic).  If we have a reflection group with a Weyl vector
then we naturally get a linear function on each Weyl chamber by taking
the inner product with the corresponding Weyl vector.  Combining these
linear functions on each Weyl chamber gives a piecewise linear
function on the Grassmannian.  We can sum up the piecewise linear
pieces that come from the singular theta transformation of the
Maass--Poincar\'e series to agree with the piecewise linear function
generated from the Weyl vector.  The remaining terms give an
automorphic form with known coefficients and known behaviour at the
cusps (it will be an ``exceptional form'').  There are conjectures
relating to the exceptional spectrum of the orthogonal groups
$\lieO(1,n)$ and we will use these conjectures to show that the
remaining terms must be zero.  From this we are able to deduce that
the original coefficients were also zero and hence the sum of
Maass--Poincar\'e series is, in fact, holomorphic.

Take $\delta\in L^*/L$, $m = \delta^2/2$, $\delta^2<0$.  Define the
function $F_\delta(\tau)$, based on the Maass--Poincar\' e series, as
$$F_\delta(\tau) = y^{\beta^+}F(\tau;\delta,1).$$
It has Fourier expansion
\begin{eqnarray*}
F_{\delta}(\tau) & = & q^m \ee_\delta + q^m \ee_{-\delta} +
\sum_{{\gamma\in L^*/L \atop n\in\ZZZ+\gamma^2/2} \atop n\ge 0}
b(\gamma,n) q^n \ee_\gamma \\
& + & \sum_{{\gamma\in L^*/L \atop n\in\ZZZ+\gamma^2/2} \atop n<0}
b(\gamma,n) e^{2\pi inx} W(4\pi ny) \ee_\gamma.
\end{eqnarray*}
Recall, the function $W$ is a Whittaker function that tends to zero
exponentially as $y$ tends to infinity.  Note that the function
$F_\delta$ is obtained from the Maass--Poincar\'e series at the
eigenvalue ``s=1''.  If the lattice $L$ has signature $(1,n)$ with
$n\ge 2$ then it is easy to check that this choice of $s$ is within
the half--plane of convergence.

The coefficients $b(\gamma,n)$ are real numbers (remark after Lemma
4.6 in \cite{bruinier_borcherds}) that are bounded for negative $n$
(Equation 6.10 in \cite{bruinier_borcherds}).

Define the {\it $\delta$--Heegner divisor} to be
$$H(\delta) = \bigcup_{{\lambda\in L^* \atop \delta\in\lambda +L} \atop
\lambda^2 = \delta^2} \lambda^\perp.$$

This is a locally finite collection of subspaces of the Grassmannian;
these are exactly the subspaces on which the function
$\Phi_{\delta}(v)$ (the singular theta transformation of $F_\delta$)
has singularities.  Near to a point $v_0$ on the Heegner divisor the
singularity is of the form
\begin{equation}
-4\sqrt{2}\pi\sum_{\lambda\in L^*\cap v_0^- \atop \lambda\ne 0}
|\lambda_{v^+}|.
\label{sing}
\end{equation}
We will abuse notation by also using $H(\delta)$ to represent the set
of vectors that represent the sub--Grassmannians.  It will be clear
from the context exactly which set we are talking about.

The Fourier expansion for $\Phi_{\delta}$ can be worked out using the
Rankin--Selberg method.  To do this we fix a primitive norm $0$ vector
$z\in L$ and $z'\in L^*$ such that $(z,z')=1$ (such a $z'$ exists
because $z$ is primitive).  Let $N$ be the minimum positive (integer)
inner product of $z$ with vectors in $L$.  The vector $z$ represents a
cusp of $\Gr(L)$.  Define the lattice $K$ to be $(L\cap z^\perp)/\ZZZ
z$.  It is clear that $K$ is a negative definite lattice.  However,
even though $K\subset L$, it is not always true that $K^* \subset L^*$.
Define $p: L\tensor\RRR \rightarrow K\tensor\RRR$ to be the orthogonal
projection.  Any element of $L^*$ will clearly project to an element
of $K^*$.  Define
$$L_0^* = \{ \lambda\in L^*: (\lambda,z)\equiv 0\ (\mod N)\}.$$
Then, providing $z_{v^+}^2$ is sufficiently small, the following sum
converges (Proposition 9.1 in \cite{bruinier_borcherds}):
\begin{eqnarray*}
\Phi_{\delta}(v) & = & {1\over \sqrt{2}\,|z_{v^+}|}\Phi^K +
4\sqrt{2}\,\pi|z_{v^+}|
\sum_{l=0}^{N-1} b(lz/N,0)\BBB_2(l/N) \\
& + & 4\sqrt{2}\,\pi |z_{v^+}| \sum_{{\lambda\in K^* \atop p(\delta)\in
\lambda+K} \atop \lambda^2 = \delta^2} \BBB_2 ((\lambda, \mu) +
(\delta, z')) \\
& + & 4\sqrt{2}\,(\pi / |z_{v^+}|)^{(n-1)/2} \sum_{\lambda\in K^*-0}
\sum_{\gamma\in L_0^*/L \atop p(\gamma)\in \lambda+K}
b(\gamma,\lambda^2/2) |\lambda|^{(n+1)/2} \\
&\times & \sum_{r\ge 1} r^{(n-3)/2} 
\e(r(\lambda,\mu)+r(\gamma,z')) K_{(n+1)/2}\(2\pi r |\lambda|/|z_{v^+}|\).
\end{eqnarray*}
$\Phi^K$ is a constant determined by performing the singular theta
lift on $F$ with respect to the negative definite lattice $K$.  The
function $K_n(x)$ is the $K$--Bessel function \cite{abramowitz,
erdelyi_I}.  

\section{Decomposition of $\Phi_{\delta}$}

Bruinier showed that the function $\Phi_{\delta}(v)$ can be decomposed
as a sum of two function $\psi(v)$ and $\xi(v)$, where $\xi(v)$ is
real analytic on the whole of $\Gr(L)$ and $\psi(v)$ is a piecewise
linear function on $\Gr(L)$.  Both of these functions are
eigenfunctions of the hyperbolic Laplacian on $\Gr(L)$ (outside of
their singular sets).  Note that the decomposition of $\Phi_\delta(v)$
is not canonical --- it depends on the choice of the vectors $z$ and
$z'$.

The {\it hyperboloid model} for $\Gr(L)$ is given by representing
every maximal positive definite space of $L\tensor\RRR$ by the unique
norm $1$ vector $v_1$ it contains having positive inner product
with $z$.  Using this model we can write down the equation for
$\xi(v)$ (Definition 9.3 in \cite{bruinier_borcherds}):
\begin{eqnarray*}
\xi_{\delta}(v_1) & = & {\Phi^K\over \sqrt{2}}\({1\over (z,v_1)} -
2(z',v_1)\) + 
{4\sqrt{2}\,\pi \over (z,v_1)}\sum_{{\lambda\in K^* \atop p(\delta)\in
\lambda+K} \atop \lambda^2/2 = m} (\lambda, v_1)^2 \\
& + & 4\sqrt{2}\,(\pi / (z,v_1))^{(n-1)/2} \sum_{\lambda\in K^*-0}
\sum_{\gamma\in L_0^*/L \atop p(\gamma) \in \lambda+K}
b(\gamma,\lambda^2/2) |\lambda|^{(n+1)/2} \\
&\times & \sum_{r\ge 1} r^{(n-3)/2}\ \e\(r{(\lambda,v_1)\over (z,v_1)} +
r(\gamma,z')\) K_{(n+1)/2} \({2\pi r |\lambda| \over(z,v_1)}\).
\end{eqnarray*}
This function is obviously real analytic due to the exponential decay
of the $K$--Bessel function and the boundedness of the coefficients
$b(\gamma,n)$ for $n<0$.

The function $\psi$ is just the difference between $\xi$ and the
function $\Phi$.  Although it does not look like it is piecewise linear
(because it has quadratic terms from the Bernoulli polynomials) it is
seen that these cancel out near to the cusp (see either
\cite{borcherds_grass} or \cite{bruinier_borcherds}).

Another model for the Grassmannian manifold is the {\it upper half
space model}.  Define $K$ to be the lattice $L\cap z^\perp \cap
z'^\perp$ (this is isomorphic to the definition of $K$ used before).
Any norm $1$ vector $v_1$ can then be represented as a vector in
$K\tensor\RRR$ and a positive multiple of the vector $z'$ (the $z$
coordinate is then chosen to give norm $1$).  This gives coordinates
for the upper half space model as $\RRR_{>0} \times (K\tensor \RRR)$.
We denote the coordinates $(y,\mu)$ where $y\in\RRR_{>0}$ and $\mu\in
K\tensor\RRR$.  In these coordinates the equation for $\xi$ becomes
(Section 9.1 of \cite{bruinier_borcherds}).
\begin{eqnarray*}
\xi_{\delta}(y,\mu) &=& {\mu^2 \over \sqrt{2}\,y} \Phi^K +
{4\sqrt{2}\,\pi \over y} \sum_{{\lambda\in K \atop p(\delta)\in\lambda+K}
\atop \lambda^2/2 = m} (\lambda,\mu)^2 \\
&+& 4\sqrt{2}\,(\pi y)^{(n-1)/2} \sum_{\lambda\in K^*-0} \sum_{\gamma \in
L_0^*/L \atop p(\gamma)\in \lambda+K} b(\gamma, \lambda^2/2)
|\lambda|^{(n+1)/2} \\
&\times& \sum_{r\ge 1} r^{(n-3)/2} e(r(\lambda,\mu)+r(\gamma,z'))
K_{(n+1)/2}(2\pi r |\lambda|y).
\end{eqnarray*}

In these coordinates the hyperbolic Laplacian is
$$\Delta = -y^2\({\d^2 \over \d y^2} + {\d^2 \over \d x_1^2} + \cdots
+ {\d^2 \over \d x_n^2}\) + (n-2)y {\d \over \d y},$$
where $\{x_i\}$ form orthogonal coordinates for the space
$K\tensor\RRR$.  The hyperbolic Laplacian is invariant under the
action of $\lieSO (1,n)$ \cite{grunewald}.

By explicit calculation (Theorem 9.7 of \cite{bruinier_borcherds})
one can show that the functions $\xi$ and $\psi$ are eigenfunctions of
the hyperbolic Laplacian with eigenvalues $-n$.

\section{Piecewise Linear Functions}

In this section we show that any piecewise linear function on
$\Gr(L)$, invariant under the automorphism group of the lattice, with
singularities lying on Heegner divisors can be obtained as the
singular theta transformation of a linear combination of
Maass--Poincar\' e series.  This can be thought of as the $\lieO
(1,n)$ version of a similar result for $\lieO (2,n)$ in
\cite{bruinier_borcherds}.  This is achieved by picking the only
possible choice of Maass--Poincar\' e series (which is determined by
the form of the singularities along the Heegner divisors).  We must
then show that the corresponding linear combination of the functions
$\xi$ is zero.  To show that $\xi$ is zero we will calculate its
behaviour at the cusps and under the automorphism group of the lattice
to show that it is either a cusp form or an exceptional form for the
Laplacian.  From some general conjectures about the exact structure of
the eigenvalues for the Laplacian of orthogonal groups we are able to
deduce that $\xi$ must be identically zero.

Let $P(v)$ be a piecewise linear function that is invariant under
automorphisms of the lattice that are trivial on the discriminant
group and having singularities given by a linear combination of
Heegner divisors
$$(P) = {1\over 2} \sum_{\delta\in L^*/L}c(\delta) H(\delta).$$
The divisor notation means that the locations and form of the
singularities of the function $P$ are of the types given by Equation
(\ref{sing}) (scaled by $c(\delta)$).

\nb In many cases the choice of Heegner--divisors is unique, but there
are situations where this may not be true.  For now we will ignore
this problem and assume that we have chosen Heegner--divisors giving
the correct singulariries for $P$.

Define a real analytic modular form by $F_P(\tau) = \sum_{\delta}
c(\delta)F_{\delta}(\tau)$.  Let $\Phi_P$ be its singular theta
transformation.

\begin{lemma} The function $P-\Phi_P$ is real analytic on $\Gr(L)$,
invariant under automorphisms of the lattice that are trivial on the
discriminant group and an eigenfunction of the hyperbolic Laplacian
with eigenvalue $-n$.
\end{lemma}
\proof The singularities of $\Phi_P$ are exactly the same as the
singularities of $P$ (by construction).  So their difference is
singularity free.  Consider $P-\psi$.  This is a piecewise linear
function with no singularities.  Hence it is a linear function.  Let
us suppose that it is of the form $(v_1,\nu)$.  Note that a direct
calculation shows that any linear function is an eigenfunction of the
Laplacian with eigenvalue $-n$.  By construction, $\Phi_P$ is
invariant under automorphisms of the lattice that are trivial on the
discriminant group.  Hence the function $P-\Phi_P$ has the required
properties.\qed

\begin{lemma} The group of automorphisms of the lattice $L$ that are
trivial on the discriminant group is a congruence subgroup of
$\Aut(L)$.
\end{lemma}
\proof Let the maximum order of an element of $L^*/L$ be $m$.  Regard
elements of $\Aut(L)$ as matrices (by choosing an integral basis for
$L$).  All automorphisms of the form
$$\sigma \equiv \id\ (\mod m)$$
act trivially on the discriminant group.  \qed

In particular we know that automorphisms of the above form are of
finite index in $\Aut(L)$.  Let us examine the behaviour of $P-\Phi_P$
at the cusps.

\begin{lemma} The endomorphism of $L\tensor\RRR$ defined by
\begin{eqnarray*}
(0,0,1) &\longmapsto& (0,0,1)\\
T_\vecx : (0,1,0) &\longmapsto& (\vecx,1,-\vecx^2/2)\\
(\lambda,0,0) &\longmapsto& (\lambda,0,-\vecx\cdot\lambda)
\end{eqnarray*}
is an element of $\lieO_L(\RRR)$.  If $L$ is even and $\vecx\in K$, or
if $L$ is odd and $\vecx\in 2K$, then $T_\vecx$ is an element of
$\Aut(L)$ preserving the discriminant group.
\end{lemma}
\proof It is easy to check that this endomorphism preserves the inner
product.  It is an automorphism because $T_\vecx \circ T_{-\vecx} =
T_0 = I$.

Suppose that $(\veck,a,b) \in L^*$.  Taking the inner product of this
with $z$ shows that $a\in\ZZZ$.  Taking the inner product with
elements of $K$ shows that $\veck\in K^*$.  We know that
$$T_\vecx: (\veck,a,b) \mapsto (\veck+a\vecx,a,b-a\vecx^2/2-\vecx\cdot
\veck).$$
The difference between the image and the original vector is thus
$$(a\vecx,0,-a\vecx^2/2-\vecx\cdot \veck).$$
This vector is in the lattice $L$.\qed

The action of this map on the upper half plane model is
$$T_\vecx: (y,\mu) \mapsto (y,\mu+\vecx).$$
As both $P$ and $\Phi_P$ are invariant under automorphisms of $L$
preserving the discriminant form, $P-\Phi_P$ is invariant too.

\begin{theorem} The function $P-\Phi_P$ decreases as $1/y$ (or faster)
at the cusp represented by $z$.
\end{theorem}
\proof Pick some $z' \in L^*$ with $(z,z')=1$ and decompose $\Phi_P$
into $\psi + \xi$.  Fix these functions and define the vector $\nu$ by
$$P(v_1) - \psi(v_1) = (v_1,\nu).$$
Such a vector exists as $P=\psi$ is a linear function.  The Fourier
expansion of $P(v_1)-\Phi_P(v_1)$ is then
\begin{eqnarray}
& & (v_1,\nu) - {\Phi^K\over \sqrt{2}}\({1\over (z,v_1)} -
2(z',v_1)\) -
{4\sqrt{2}\,\pi \over (z,v_1)}\sum_{\lambda\in (P)}
(\lambda, v_1)^2 \nonumber\\
& - & 4\sqrt{2}\,(\pi / (z,v_1))^{(n-1)/2} \sum_{\lambda\in K^*-0}
\sum_{\gamma\in L_0^*/L \atop p(\gamma) \in \lambda+K}
b(\gamma,\lambda^2/2) |\lambda|^{(n+1)/2} \nonumber\\
&\times & \sum_{r\ge 1} r^{(n-3)/2} \e\(r{(\lambda,v_1)\over (z,v_1)} +
r(\gamma,z')\) K_{(n+1)/2} (2\pi r |\lambda|/(z,v_1)).
\label{fourier1}
\end{eqnarray}
The function $P-\Phi_P$ is invariant under automorphisms of the
lattice fixing the discriminant group.  In particular, it is invariant
under $T_\vecx$ for any $x\in K$.  So, Fourier expansion
(\ref{fourier1}) is equal to
\begin{eqnarray}
& & (T_\vecx v_1,\nu) -
{\Phi^K\over \sqrt{2}}\({1\over (z,v_1)} - 2(z',T_\vecx v_1)\) -
{4\sqrt{2}\pi \over (z,v_1)}\sum_{\lambda\in (P)}
(\lambda, T_\vecx v_1)^2 \nonumber\\
& - & 4\sqrt{2}\,(\pi / (z,v_1))^{(n-1)/2} \sum_{\lambda\in K^*-0}
\sum_{\gamma\in L_0^*/L \atop p(\gamma) \in \lambda+K}
b(\gamma,\lambda^2/2) |\lambda|^{(n+1)/2} \nonumber\\
&\times & \sum_{r\ge 1} r^{(n-3)/2} \e\(r{(\lambda,T_\vecx v_1)\over
(z,v_1)} + m(\gamma,z')\) K_{(n+1)/2} (2\pi r |\lambda|/(z,v_1)).
\label{fourier2}
\end{eqnarray}
In simplifying (\ref{fourier2}) we have used the fact that $T_\vecx$
stabilizes the vector $z$.  After equating (\ref{fourier1}) and
(\ref{fourier2}) most terms cancel and we are left with (assuming that
$v_1 = (\vecp,a,b)$ and $\nu = (\vecq,c,d)$)
\begin{eqnarray}
-a(\vecq,\vecx) + ac{\vecx^2\over 2} + c(\vecp,\vecx) &=&
\sqrt{2}\,\Phi^K \( -a {\vecx^2 \over 2} - (\vecp,\vecx)\) \nonumber\\
&-& 4\sqrt{2}\,\pi \sum_{\lambda\in (P)} \[ 2(\lambda,\vecp)
(\lambda,\vecx) + a(\lambda, \vecx)^2 \].
\label{equated}
\end{eqnarray}
Regarding $\vecp$, $a$, $b$ and $\vecx$ as variables (\ref{equated})
becomes a polynomial equation and the terms can be extracted by
degree.  In particular we find that
$$\vecq = 0$$
and
$$c(\vecp,\vecx) = \sqrt{2}\Phi^K (\vecp,\vecx) - 8\sqrt{2}\,\pi
\sum_{\lambda\in (P)} (\lambda,\vecp)(\lambda,\vecx).$$
Which means that
\begin{equation}
\sum_{\lambda \in (P)} (\lambda,\vecp)(\lambda,\vecx) = -{\sqrt{2}\,
\Phi^K + c \over 8\sqrt{2}\,\pi}(\vecp,\vecx).
\label{vectorsystem}
\end{equation}
This shows that the sum on the left hand side, if $\vecp = \vecx$, is
constant on the sphere of norm $-1$ vectors.  This means that the
vectors $\lambda \in (P)$ form a {\it vector system} 
\cite{borcherds_products}.

If we substitute (\ref{vectorsystem}) back (\ref{fourier1}) and recall
that $v_1$ has norm $1$ we get
\begin{eqnarray*}
P(v_1) - \Phi_P(v_1) &=& (v_1,\nu) - {\Phi^K \over \sqrt{2}} \(
{1\over (z,v_1)} - 2(z',v_1)\) + {\sqrt{2}\,\Phi^K + c \over
2(z,v_1)} \vecp^2 + \O \\
&=& (v_1,\nu) - {\Phi^K \over \sqrt{2}\,a} \( 1-2ab-2a^2(z')^2-\vecp^2
\) + {c \over 2a}\vecp^2 + \O \\
&=& (v_1,\nu) + \sqrt{2}\,\Phi^Ka(z')^2 + {c \over 2a}\vecp^2 + \O \\
&=& bc + ad + ac(z')^2 + \sqrt{2}\,\Phi^Ka(z')^2 + {c \over 2a}\vecp^2 +
\O \\
&=& a\( d + \sqrt{2}\,\Phi^K (z')^2 + c (z')^2 \) + {c \over 2a} + \O,
\end{eqnarray*}
where $\O$ is shorthand for the rapidly decreasing terms containing
the $K$--Bessel functions.

Now, looking at the action of the Laplacian on this function we see
that there can be no terms of the form $1/a$ and hence $c=0$.

Since $\O$ decreases rapidly as one approaches the cusp and the terms
proportional to $a$ decrease as the reciprocal of a polynomial, the
function $P(v_1) - \Phi_P(v_1)$ decreases as $1/y$ (or faster).\qed

We note that in proving the above result we also showed the following:

\begin{corollary} The elements of $(P)$ form a vector system.
\end{corollary}

If we examine the Fourier expansion for $\xi$ knowing that the vectors
$\lambda \in (P)$ form a vector system we see:

\begin{corollary} $\xi$ is rapidly decreasing at the cusp represented by
$z$.
\end{corollary}

\section{The Spectrum of the Laplacian}

Set $\omega = {n-1 \over 2}$ and write the eigenvalues of the
Laplacian in the form $s^2-\omega^2$.  Because the eigenvalues for the
Laplacian are negative reals, this forces $s$ to lie in the following
subset of $\CCC$
$$i\RRR \cup [-\omega,\omega].$$

The continuous spectrum of the Laplacian lies in $i\RRR$ and it is
conjectured that the cuspidal part also lies in this region (these are
the {\it Ramanujan--Selberg conjectures}).  The other eigenvalues that
can occur are called {\it exceptional eigenvalues}.  They can come
from either theta liftings or from the residues at poles of the
analytic continuation of the Eisenstein series (see, for example,
\cite{howe_ramanujan}).  These functions tend to zero at the cusps but
are not cuspidal (because they do not tend to zero exponentially).
They correspond to the non--tempered automorphic representations (see,
for example, \cite{bump_automorphic}).  There is a precise conjecture
as to the location of all eigenvalues for the orthogonal group.

\begin{conjecture}\label{conj} For the group $\lieO (1,n)$ the
exceptional eigenvalues are given by
$$s\in\{\omega, \omega-1, \dots,-\omega\}.$$
\end{conjecture}

This conjecture is formulated in \cite{sarnak_ramanujan} and
\cite{sarnak_ramanujan2} and in these papers various computations are
performed to support the conjectures.  It is also noted in these
papers that the conjecture appears consistent with various conjectures
of Arthur \cite{arthur}.

For $\lieSO (1,2)$ (which is $\liePSL_2(\RRR)$) this says
that there are no exceptional eigenvalues other than the constant
function (which is the residue of the pole of an Eisenstein series at
$s=1$).  In particular, other than this zero eigenvalue, all other
eigenvalues are at most $-1/4$.  This is the famous {\it Selberg
$1/4$--conjecture}.

\begin{lemma} Assume Conjecture \ref{conj}.  For $n\ge 7$ the function
$P-\Phi_P$ is the zero function.
\end{lemma}
\proof The function $P(v_1) - \Phi_P(v_1)$ is an eigenfunction of the
Laplacian and so, by the conjecture, its eigenvalue, $-n$, should
occur for some $s\in i\RRR \cup \{ \omega,\omega-1,\dots,-\omega\}$.

For $n\ge 6$ it is clear that $n < \omega^2$ and so $s$ must be one of
the exceptional eigenvalues.  This requires us to solve
$$-n = -(n-1)k + k^2$$
in the integers, giving
$$n = {k(k+1) \over k-1}.$$
Thus $n$ is either $0$ or $6$.  Hence for $n\ge 7$ the eigenvalue $-n$
does not occur in the (conjectured) spectrum and so $P-\Phi_P$ is the
zero function. \qed

\begin{corollary} $P = \psi_P$ and $\xi_P = 0$.
\end{corollary}
\proof We have seen that $P = \Phi_P = \psi_P + \xi_P$.  $P-\psi_P$ is
some linear function and $\xi_P$ is a rapidly decreasing function at
the cusp $z$.  Hence, both are zero.\qed

\noindent{\bf Question.} Is there another way to show that $P-\Phi_P$
is zero without using Conjecture \ref{conj}?  In many examples the
difference is zero even in cases where $n\le 7$.  It is therefore
natural to conjecture that the above construction will always give
$P=\Phi_P$; even if this turns out to be false the above method will
give a way to construct exceptional eigenfunctions which may
themselves be interesting.

\section{The Maass--Poincar\' e Sum}

We have constructed a real analytic modular form whose theta
transformation is exactly the piecewise linear function we wanted.  We
would now like to show that this implies that the modular form was
actually a holomorphic one.  Unfortunately there are many
counterexamples to this for lattices with large $p$--rank.  However, if
the $p$--rank of the lattice is sufficiently small then we can make
this deduction.

\begin{lemma} There is no nearly holomorphic vector--valued modular
form associated to the lattice $\even_{1,25}(2)$ with singularities
corresponding to the primitive roots.
\end{lemma}
\proof The roots of $\even_{1,25}(2)$ are the vectors that are one
half of the vectors that were roots for $\even_{1,25}$.  These are
norm $-1$ vectors in $L^*$.  The modular form associated to this would
have singularities of the form $q^{-1/2}$ and so on multiplying be
$\Delta(\tau)$ we would obtain a non--singular weight zero modular
form that was not a constant function.  This is clearly impossible.
So, without some extra condition it is impossible to guarantee the
existence of a modular form with the singularity structure we would
like.\qed

It is possible to create a suitable modular form for the lattice
$\even_{1,25}(2)$; however it will have singularities corresponding to
non--primitive vectors.  This modular form is closely related to the
classical form $\Delta(2\tau)$: it is formed by inducing the
scalar--valued modular form $1/\Delta(2\tau)$ to a vector--valued form
(see Chapter $4$).  The form will have singularities of the form
$q^{-2}$ corresponding to norm $-4$ vectors.  Usually, some norm $-4$
vectors will not be roots of a lattice, but in the case of
$\even_{1,25}(2)$ all norm $-4$ vectors are roots (because they are
multiples of norm $-1$ vectors in the dual lattice).  This problem
does not only occur for imprimitive lattices; another example comes
from \cite{borcherds_refl} where the modular form
$$q^{-1} - 216 - 9126q + O(q^2)$$
almost corresponds to the lattice with genus
$\even_{1,15}(3^{-1})$ (the reason for ``almost'' is that the
vector--valued version of this modular form is zero).  What is
happening is that we should be looking at the Atkin--Lehner duals:
The lattice should be in the genus
$\even_{1,15}(3^{-15})$ and the modular form should be
$$q^{-3} - 90q^{-1} - 216 - 5904q^2 + O(q^3).$$
Usually this form would cause problems since some norm
$-6$ vectors may not be roots.  However, for the lattice
$\even_{1,15}(3^{-15})$ all norm $-6$ vectors are multiples of norm
$-2/3$ vectors from the dual lattice and hence they are roots.

There are other examples in \cite{borcherds_refl} having small
discriminant groups and almost corresponding scalar valued forms (they
seem always to be the ``dihedral'' cases).  These can all be explained
using Atkin--Lehner duals in the same way as above.

These examples make clear that problems occur in the cases when
Heegner--divisors for different vectors can coincide.  In the case of
$\even_{1,25}(2)$, the Heegner--divisors for norm $-4$ vectors and
those for cosets of norm $-1$ vectors are the same.  In the
$\even_{1,15}(3^{-1})$ case, the Heegner--divisors for the norm $-6$
vectors and those for cosets of norm $-2/3$ vectors are the same.  If
Heegner--divisors coincide it becomes harder to pick a canonical
choice for the singular coefficients of the modular form, \eg in the
$\even_{1,25}(2)$ case we have to pick the coefficients of $q^{-2}$
and $q^{-1/2}$ and only for the correct choice of these coefficients
will the real analytic terms cancel.  There are, however, many
examples to show that these harder lattices still seem to be
associated to modular forms.  In \cite{bruinier_priv}, Bruinier
expressed similar beliefs that the theorems carry over to these harder
cases but that some new idea is needed; it was suggested that this
might involve more careful study of the behaviour of newforms and
oldforms.

The fact that $\xi$ is the zero function gives strong restrictions on
the coefficients $b(\lambda, n)$ for $n<0$

\begin{theorem}\label{zerocoeff} For any $\lambda\in K^*-0$ we have
the following restriction on the coefficients $b(\lambda,n)$
$$\sum_{\gamma\in L_0^*/L \atop p(\gamma)
\in \lambda+K} {b(\gamma,\lambda^2/2) \over r^2} = 0.$$
\end{theorem}
\proof The formula for $\xi$ in the upper half plane model is
\begin{eqnarray*}
\xi(y,\mu) &=& 4\sqrt{2}\,(\pi y)^{(n-1)/2} \sum_{\lambda\in K^*-0}
\sum_{\gamma \in L_0^*/L \atop p(\gamma) \in \lambda+K} b(\gamma,
\lambda^2/2)|\lambda|^{(n+1)/2} \\
&\times& \sum_{r\ge 1} r^{(n-3)/2} \e(r(\lambda,\mu)+r(\gamma,z'))
K_{(n+1)/2}\(2\pi r |\lambda|y\).
\end{eqnarray*}
Note that we have removed the polynomial terms due to the vector
system property.  The $\Lambda$--Fourier coefficient is
$$4\sqrt{2}\,(\pi y)^{(n-1)/2} |\Lambda|^{(n+1)/2}
\sum_{r\cdot\lambda = \Lambda}\sum_{\gamma\in L_0^*/L \atop p(\gamma)
\in \lambda+K} {b(\gamma,\lambda^2/2) \over r^2}.$$
The fact this is zero implies that
$$\sum_{r\cdot\lambda=\Lambda} \sum_{\gamma\in L_0^*/L \atop p(\gamma)
\in \lambda+K} {b(\gamma,\lambda^2/2) \over r^2} = 0.$$
By an easy induction, starting with the primitive vectors, we see the
result.\qed

We would like to deduce from this result that the $b(\lambda,n)$ with
$n<0$ are all zero.  However, there are a couple of problems with
this.  In cases like $\even_{1,25}(2)$, the norms of vectors in the
lattice are $0, 8, 12, 16, \dots$.  In particular, half of the even
norms are missing.  However, the Maass--Poincar\'e series has terms of
the form $b(\lambda, n)$ where $n \equiv \lambda^2/2\ (\mod 1)$.  So,
the Maass--Poincar\'e series can have terms such as $b(\lambda, 7)$
which can not be ruled out by the above lemma as they never occur in
the Fourier expansion.  So we need some way to guarantee that enough
norms occur in the lattice $K^*$.  A separate problem is that the
lemma only guarantees that a sum is zero rather than the individual
terms.  If $N=1$ (\ie if $z' \in L$) then $L_0^*/L$ is trivial and the
sums all have only one term.  In the next section we will place
restrictions on the lattice $L$ which allow us to deduce that the
modular form is holomorphic.

\section{Well--Endowed Lattices}

We now make two assumptions about the lattice $L$ and will call these
lattices {\it well--endowed}.  Firstly, we assume that it splits a
hyperbolic plane.  In other words, we can decompose $L$ as
$$L = K \oplus \even_{1,1},$$
where $K$ is a negative definite lattice.  Note that this splitting is
usually not unique:  In the case $L = \even_{1,25}$ we can pick any
of the $24$ Niemeier lattices \cite{niemeier} for $K$.

Since we have split a hyperbolic plane we can pick $z$ and $z'$ from
the lattice $\even_{1,1}$ and hence get only single terms (not sums)
in Theorem \ref{zerocoeff}.  This will allow us to deduce that
$b(\lambda,n)=0$ for any norm occurring at such a cusp.

\begin{lemma} The genus of the lattice $L$ contains only one lattice.
\end{lemma}
\proof In the indefinite case, Theorem 19 from Chapter 15 of
\cite{conway_sloane} gives that form must $p$--adically diagonalize
with distinct powers of $p$ on the diagonal.  This is impossible due
to the $\even_{1,1}$ part.\qed

\begin{corollary} Let $K'$ be any lattice in the same genus as $K$.  Then
$$K'\oplus \even_{1,1} \cong L.$$
\end{corollary}
\proof The genus of $K'\oplus \even_{1,1}$ is the same as that of
$L$.  However, $L$ is unique within the genus.\qed

A genus of lattices determines the discriminant group $A$.  The genus
$g$ is said to {\it represent all possible norms} if the following
holds: for any $\gamma \in A$ and $n \equiv \gamma^2 \ (\mod 2)$ there
is a lattice $L$ in $g$ and a vector $\lambda \in L^*$ which is in the
coset represented by $\gamma$ and has norm $n$.

Well--known results of Siegel show that a genus represents every norm
not ruled out by local considerations.  In fact, Siegel computed the
average number of such representations \cite{siegel}.  To be able to
use Theorem \ref{zerocoeff} we want our genus to represent all
possible norms, so we now investigate conditions on $L$ for which we
can guarantee this.

\begin{lemma} Suppose that the even lattice $L$ contains $\even_{1,1}$
as a sublattice.  Then $L$ represents all possible norms.
\end{lemma}
\proof Write $L = M\oplus \even_{1,1}$ so that $L^*/L$ is identified
with $M^*/M$.  Pick any lift of $\gamma$ to a vector of $M^*$.  Because
the lattice $\even_{1,1}$ contains vectors of every even norm we can
add on one of these to adjust the norm.\qed

There are lattices that do not represent all possible norms.  An
imprimitive lattice $L(2)$ will miss certain norms.  Being primitive
is not sufficient: Let $L$ be any indefinite lattice and $M$ a
one--dimensional lattice then $L(p) \oplus M$ is a primitive lattice
that only represents squares or non--squares modulo $p$.  Another
example is $L(25) \oplus \smallmatrix{2}{0}{0}{4}$ which does not
represent any even integer divisible exactly once by $5$.  So, there
seem to be problems when the $p$--rank of the lattice is large.  It is
well--known that, when the $p$--rank is sufficiently small compared to
the dimension, the lattice splits a hyperbolic plane.

\begin{lemma} Suppose that $L$ is an even indefinite lattice such that
$$\prk(L) \le \dim(L)-3\quad \hbox{for all primes $p$}$$
then $L \cong M\oplus \even_{1,1}$ for some lattice $M$.
\end{lemma}
\proof This is Corollary 1.13.5 in \cite{nikulin_discriminant}.\qed

The third example above shows that this lemma is optimal, \ie there
exist lattices with $\prk(L) = \dim(L)-2$ for some primes $p$ not
splitting hyperbolic planes.

\begin{theorem} Suppose that $L$ is an even indefinite lattice such that
$$\prk(L) \le \dim(L)-5\quad \hbox{for all primes $p$}.$$
If we are given $r\in\QQQ$, $r<0$ and $\gamma\in L^*/L$ such that
$\gamma^2 \equiv r\ (\mod 2)$ then $\gamma$ can be lifted to a vector of
$L^*$ of norm $r$ such that there is a copy of $\even_{1,1}$ in its
orthogonal complement.
\end{theorem}
\proof We already know from the previous lemma that we can lift
$\gamma$ to a vector $\beta$ in $L^*$ of norm $r$.  Suppose that the
minimum positive inner product of $\beta$ with $L$ is $N$.  Let $M =
\beta^\perp \cap L$ and $A = \ZZZ\beta \cap L$.  $L$ can be formed by
``gluing'' together $A$ and $M$ (see Chapter 4 of
\cite{conway_sloane}).

Let $l$ be any vector of $L$ such that $(l,\beta)=N$.  By projection
onto $M\tensor\RRR$ we see that there is $m'\in M^*$ such that $m' + {N
\over \beta^2} \beta \in L$.  This gives us a gluing relation and due
to the minimality of $N$ it generates all the gluing relations.  Let
$G$ be the group of all gluing relations considered as a subgroup of
$A^*/A \oplus M^*/M$.  We have
\begin{diagram}
0 & \rTo & G & \rInto & A^*/A \oplus M^*/M & \rOnto & L^*/L & \rTo & 0.
\end{diagram}
As $G$ has rank at most $1$ we see that the $\hbox{$p$-rank}$ of $M$
is at most $1$ more than the $\hbox{$p$-rank}$ of $L$.  Hence $M$ can
be split as $M = T\oplus \even_{1,1}$ by Nikulin's result.  Gluing
this back to $A$ we have the required decomposition. \qed

\begin{corollary}\label{allposs} If $K$ is a negative definite lattice
with
$$\prk(K) \le \dim(K) - 3\quad\hbox{ for all primes $p$}$$
then the genus of $K$ represents all possible norms.
\end{corollary}
\proof Form $L = K\oplus\even_{1,1}$.  $L$ satisfies the conditions
of the previous theorem.  Using this theorem we can find copies of
vectors of all possible norms with an $\even_{1,1}$ orthogonal to
them.  Taking the orthogonal lattice to this $\even_{1,1}$ gives a
lattice in the genus of $K$ representing the required vector.\qed

As before, this bound on the $p$--ranks is optimal (as the face
centred cubic lattice shows).

So, we restrict to lattices $L$ that split a hyperbolic plane and
determine a genus representing all possible norms.  We call this
class of lattices {\it well--endowed}.

\begin{theorem} Assume Conjecture \ref{conj}.  Suppose that $L$ is a
well--endowed lattice with signature $(1,n)$, $n\ge 7$.  Suppose that
$P(v)$ is a continuous piecewise linear function on $\Gr(L)$ with
singularities along Heegner--divisors of $L$.  Then $P$ is given by
the singular theta lift of some nearly holomorphic vector--valued
modular form.
\end{theorem}
\proof We saw how to pick the candidate function lifting to $P$.
We now need to show that the real analytic terms cancel.  Since $L$
is well--endowed we get (using Theorem \ref{zerocoeff} at all the cusps)
that
$$b(\lambda, n) = 0\quad\hbox{for all possible $\lambda$ and $n$}.$$
Thus there are no real analytic terms in the function $F_P(\tau)$
so it is nearly holomorphic.\qed

This theorem can be regarded as the $(1,n)$ version of main theorem
proved by Bruinier in \cite{bruinier_borcherds}.  The restriction on
the dimension is only needed due to our use of Conjecture \ref{conj}.

It is worth commenting on the similarity between the conditions we
impose on the lattice and those that Bruinier requires.  We impose
conditions on the lattice so that we can deduce that vectors of all
possible norms occur as we vary over cusps.  Bruinier requires that
the lattice of signature $(2,n)$ splits two hyperbolic planes.  This
allows him to guarantee that vectors of all possible norms occur at a
single cusp.  Notice that splitting two hyperbolic planes is
guaranteed for signature $(2,n)$ lattice with $\prk(L) \le \dim(L)-5$
for all primes $p$ (simply apply Nikulin's result twice).

If the lattice is not well--endowed it can still be possible to
get results by using Atkin--Lehner involutions.

A lattice is called {\it elementary} if its discriminant group $L^*/L$
is elementary Abelian.  It is well known that every lattice has an
embedding (actually, a canonical embedding) into an elementary lattice
of the same dimension.

\begin{lemma} Every lattice $L$ has an embedding into an elementary
lattice.
\end{lemma}
\proof Let $A$ be a $p$--component of the discriminant group with
exponent $p^n$ for $n>1$.  The subgroup $H=p^{n-1}A$ is isotropic,
non-trivial and invariant under $\lieO(L)$.  Thus $L+H$ is an integral
lattice into which $L$ embeds.  This process decreases the size of the
discriminant group as so is clearly finite.\qed

If $m\in\NNN$ and $m$ square free then we can define the {\it
$m$--dual} of an elementary lattice $L$ to be
$$L^{*m} = \( L^* \cap {1\over m}L\)(m),$$
where $(m)$ means we scale the inner product by $m$.

\nb This operation does not always preserve the fact that the lattice
is even.

These operations are also called {\it Atkin--Lehner involutions}.
Recall that classical modular forms can be though of as functions on
the space of $2$--dimensional lattices \cite{serre}.  The action
of the $m$--dual on these lattices then gives an operation on the
space of modular forms.  This operation is the same as one described
by Atkin and Lehner in \cite{atkin_lehner} (it is the $W$ operation
associated to the matrix $\smallmatrix{0}{1}{m}{0}$).

\begin{lemma} $L^{*m}$ is an elementary lattice and has the following
properties
$$\prk(L^{*m}) = \left\{\begin{array}{ll}
\prk(L) & (p,m) \ne 1, \\
\dim(L) - \prk(L) & (p,m) = 1.\end{array}\right.$$
\end{lemma}
\proof It is easy to see that $L^{*m}\tensor \ZZZ_p = (L\tensor
\ZZZ_p)(m)$ if $p$ and $m$ are coprime and that $L^{*m}\tensor \ZZZ_p
= (L\tensor \ZZZ_p)^*(m)$ otherwise.\qed

\begin{lemma} The $*m$-operation is an involution, \ie
$$(L^{*m})^{*m} = L.$$
\end{lemma}
\proof If $(p,m)=1$ then
\begin{eqnarray*}
\(L^{*m}\)^{*m}\tensor\ZZZ_p &=& \(L^{*m}\tensor\ZZZ_p\)(m) \\
&=& (L\tensor\ZZZ_p)(m^2) \\
&=& L\tensor\ZZZ_p.
\end{eqnarray*}
If $(p,m)\ne 1$ then
\begin{eqnarray*}
\(L^{*m}\)^{*m}\tensor\ZZZ_p &=& \(L^{*m}\tensor\ZZZ_p\)^*(m) \\
&=& \((L\tensor\ZZZ_p)^*(m)\)^*(m) \\
&=& {1\over m}(L\tensor\ZZZ_p)(m^2) \\
&=& L\tensor\ZZZ_p.
\end{eqnarray*}
{}\qed

Exactly the same idea can be used to show that the $*m$--operators are
multiplicative
$$(L^{*m})^{*n} = L^{*mn}\quad\hbox{if $(m,n)=1$}.$$

Once we have embedded our lattice into an elementary lattice we can
use the Atkin--Lehner involutions to attempt to make the lattice
easier to deal with.  Recall that the dihedral cases from
\cite{borcherds_refl} could be interpreted using these Atkin--Lehner
involutions.  Presumably a more careful study of these operations
would allow our results to be extended to a more general class of
lattices.  For example, one of the problems with using Theorem
\ref{zerocoeff} was that a sum of coefficients occurs if the cusp has
level $N>1$.  If we knew more relations between these coefficients
then we could still deduce that the individual terms were zero.  When
we take Atkin--Lehner involutions of modular forms there are usually
many relations.

\chapter{Lorentzian Lattices with Weyl Vectors}

In this chapter we will use our previous results to see how a
Lorentzian lattice with a Weyl vector is associated with a nearly
holomorphic vector--valued modular form whose singularities correspond
to the primitive roots.  To do this we will construct a piecewise
linear function on the Grassmannian with singularities corresponding to
primitive roots.  The construction of the piecewise linear function
depends on the existence of the Weyl vector.  We will then use the
converse theorem of the previous section to show that this function
occurs as the singular theta lift of some nearly holomorphic modular
form.  We will use the existence of the modular form to study the
existence of such reflection groups in high dimensions.

\section{A Piecewise Linear Function}

\begin{lemma}\label{primitive} Suppose that $\vecv$ is a primitive
root of a (not necessarily even) lattice $L$.  Then $\vecv$ has
squared norm $-2/m$ for some integer $m$ and $m\vecv\in L$.
\end{lemma}
\proof As $\vecv$ is a primitive vector of $L^*$ we can find a vector
$\vecw$ in the lattice $L$ with inner product exactly $1$ with
$\vecv$.  Consider the reflection of this vector in the plane
$\vecv^\perp$.  For $\vecv$ to be a root we must have
\begin{equation}
{2\over \vecv^2}\vecv \in L.
\label{eqn}
\end{equation}
In particular, by squaring the above formula, we see $\vecv$ has
squared norm $-4/m$ for some integer $m$.  Substituting this back into
(\ref{eqn}) shows that ${m \over 2}\vecv\in L$.  The inner product of
${m \over 2}\vecv$ with $\vecw$ is both ${m \over 2}$ and an integer.
Thus, $m$ is even and we have the result.\qed

\begin{lemma} Suppose that $\vecv\in L^*/L$ has norm $-2/m$ and order
dividing $m$ then every lift of $\vecv$ to a norm $-2/m$ vector of
$L$ is a root.
\end{lemma}
\proof This follows easily from the formula for the reflection in
$\vecv$.\qed

\nb The above conditions on $\vecv$ do not quite show that the lift is
a primitive root.  It is easy to show that the lift is either a
primitive root or twice a primitive root.

Suppose that $L$ is a Lorentzian lattice.  Let $W$ be the reflection
group generated by reflections in the primitive roots of $L$.  This
group acts on $\Gr(L)$ with fundamental chamber $\F$.  Primitive roots
orthogonal to the faces of $\F$ are called {\it simple roots}.
Reflections in simple roots generate the whole reflection group.  A
{\it Weyl vector} for $W$ is a vector $\rho$ such that $(\rho,r) =
r^2/2$ for all simple roots $r$.

We will actually be able to use a slightly more general notion of a
Weyl vector.  A {\it generalized Weyl vector} for $W$ is a vector
$\rho$ such that the inner product $(\rho,r)$ depends only on the
coset of the simple root $r$ in $L^*/L$.  This is clearly a more
general notion than that of a Weyl vector; it is, however, not as
general as that in \cite{nikulin_big}.  It seems to be the case that
most interesting examples fall into this category.

Given a lattice $L$ with generalized Weyl vector $\rho$ we can
generate a piecewise linear function as follows.  Let $v_1$ be a norm
$1$ vector representing the point $v$ in $\Gr(L)$.  By action of the
Weyl group we get some norm $1$ vector $w_1$ in $\overline{\F}$ (the
non--uniqueness for elements of the boundary does not matter).  Define
$F(v_1)$ to be the inner product of $w_1$ with $\rho$.  This is
clearly piecewise linear and has its singularities on the reflection
hyperplanes.  The singularity along the hyperplane orthogonal to
the root $r$ is of the form $|(v,r)|$.

Any automorphism of the lattice that is trivial on the discriminant
form and fixes the fundamental chamber $\F$ will fix a generalized
Weyl vector. So, the piecewise linear function is invariant under a
congruence subgroup of the full automorphism group.  Using the results
of the previous chapter we can associate to this function a modular
form whose singularities encode the reflection hyperplanes of the
lattice provided that the lattice satisfies the previously discussed
conditions.

\begin{theorem} Assume Conjecture \ref{conj}.  If $L$ has a
generalized Weyl vector and $L$ is well--endowed then there is a
nearly holomorphic modular form whose singularities correspond to the
primitive roots of the reflection group of $L$.
\end{theorem}

As noted before, a sufficient condition for being well--endowed is
that $L$ has $\prk(L) \le \dim(L)-5$ for all primes $p$.  Note that
the singularities correspond to primitive roots of $L$ and so they are
of the form $q^{-1/n}$ for various integers $n$ by Lemma
\ref{primitive}.

\begin{theorem} Assume Conjecture \ref{conj}.  If $L$ is a
well--endowed Lorentzian lattice with generalized Weyl vector then
$\dim(L)\le 26$.
\end{theorem}
\proof By the previous theorem, associate to $L$ a nearly holomorphic
modular form $F(\tau)$.  Consider the form $F(\tau)\Delta(\tau)$,
where
$$\Delta(\tau) = q\prod_{i=1}^\infty \(1-q^i\)^{24} = q - 24q^2 +
252q^3 - 1472q^4 + 4830q^5 - \cdots$$
is the cusp form of weight $12$ for $\lieSL_2(\ZZZ)$.
$F(\tau)\Delta(\tau)$ has no singularities at the cusps because all
singularities of $F(\tau)$ are of the form $q^{-1/n}$.  Therefore it
is holomorphic with weight $13 - \dim(L)/2$.  It is well known
that a modular form of negative weight has singularities, so
$$13 - \dim(L)/2 \ge 0\quad\Rightarrow\quad \dim(L) \le 26.$$
This is the required result.  Note that the dimensions involved here
are all within the ranges that the theorems work.\qed

There is a famous Lorentzian lattice of dimension $26$ with a (norm
$0$) Weyl vector.  This lattice is $\even_{1,25}$.  It is closely
related to the Leech lattice: Let $\Lambda$ denote the Leech lattice
(the unique even unimodular lattice with minimal norm $4$ and
dimension 24).  Then $\even_{1,25} \cong \Lambda(-1) \oplus
\even_{1,1}$.  This lattice was shown to have a norm zero Weyl vector
by Conway \cite{conway_sloane}.  It can be shown to have a norm zero
vector by using the fact that it is associated to the modular form
$1/\Delta(\tau)$ \cite{borcherds_grass}.  We will now show that this
lattice is essentially the unique lattice with these properties.

\begin{theorem}\label{maxdim26}
$\even_{1,25}$ is the unique even, well--endowed, $26$ dimensional
lattice with a generalized Weyl vector.
\end{theorem}
\proof The only holomorphic modular forms of weight $0$ (which is the
weight of $F(\tau)\Delta(\tau)$) are constants.  Hence $F(\tau)$ is
of the form $\vecv \times 1/\Delta(\tau)$ where $\vecv$ is a constant
vector.  The singularity $q^{-1}$ corresponds to reflections in a
primative root of norm $-2$.  However, any primitive root of $L$ of
norm $-2$ is a vector of $L$ by Lemma \ref{primitive}.  Hence the only
component of $\vecv$ that is non-zero is the one corresponding to
$\ee_0$.  Now consider the behaviour of the modular form under the
element $S$ of $\lieMP_2(\ZZZ)$.  It is clear that this can only
satisfy the required transformation formula if $|L^*/L|=1$.  Hence $L$
is unimodular, so it must be $\even_{1,25}$. \qed

\section{Lifting from Scalar Forms}

Although vector--valued modular forms are theoretically convenient
they can become cumbersome if $|L^*/L|$ is large.  Borcherds observed
in \cite{borcherds_refl} that most of the vector--valued forms arising
from Lorentzian lattices were closely related to scalar--valued
modular forms on $\Gamma_0(N)$ where $N$ is the level of the lattice
$L$.  Bruinier and Bundschuh \cite{bruinier_p} have shown this
to be true in the case where $L^*/L \cong \ZZZ/p\ZZZ$.  We shall show
that it is true for any $\Aut(L^*/L)$--invariant vector--valued form
when $L^*/L$ is elementary.  The condition of being
$\Aut(L^*/L)$--invariant is clearly necessary.  There are examples in
\cite{nils_priv} of vector--valued forms on non--elementary lattices
that are not directly associated to scalar--valued modular forms.
These cases seem to be constructed from a scalar--valued modular form
and a lattice theta function, although this theta function remains a
little mysterious.

The Weil representation obviously has a commuting action of
$\Aut(L^*/L)$ on it (the Weil representation can be defined using the
inner product only).  The discriminant form can be decomposed into the
orthogonal direct sum
$$L^*/L = A_2 \oplus A_3 \oplus A_5 \oplus\cdots,$$
where the $A_p$ is an elementary Abelian $p$--group (we are assuming
that the lattice $L$ is elementary).  The automorphism group therefore
decomposes as
$$\Aut(L^*/L) = \Aut(A_2) \times \Aut(A_3) \times \Aut(A_5)
\times\cdots.$$
Therefore the fixed points of $\Aut(L^*/L)$ on the Weil representation
decompose as
$$\CCC[L^*/L]^{\Aut(L^*/L)} = \CCC[A_2]^{\Aut(A_2)} \tensor
\CCC[A_3]^{\Aut(A_3)} \tensor \CCC[A_5]^{\Aut(A_5)} \tensor \cdots.$$
We examine these orbits more closely.

\begin{lemma} Assume $p$ is odd.  Two non--zero vectors of equal norm
in $A_p$ are equivalent under $\Aut(A_p)$.
\end{lemma}
\proof Realize the discriminant form $A_p$ as the discriminant form of
a sufficiently large indefinite lattice $L$.  Given two non--zero
vectors of equal norm we can lift them to primitive vectors of small
positive norm.  We can now decompose the lattice $L$ as $(\ZZZ v\oplus
K) + g$ (where $g$ is a glue vector).  As $L$ was sufficiently large
we know that $K$ is unique \cite{nikulin_discriminant} and hence
we can write down an obvious automorphism of $L$ swapping the two
lifts.  This gives an automorphism of $A_p$ swapping the two
vectors.\qed

The previous lemma can also be proved by using Witt's extension
theorem for quadratic forms.  However, Witt's theorem does not extend
to the case $p=2$ whereas the above argument does.

For the $p=2$ case there is an extra possible orbit given by the
parity vectors.  A {\it parity vector} (see page xxxiv of
\cite{conway_sloane}) in $A_2$ is a vector $v$ such that $(v,w) \equiv
w^2\ (\mod 1)$ for all $w \in A_2$.  It is clear from the definition
that the parity vector is unique.  Other terminology for parity
vectors appears in the literature: they are also known as
characteristic vectors, canonical elements and test vectors.

\begin{lemma} Two non--zero, non--parity vectors of equal norm
in $A_2$ are equivalent under $\Aut(A_2)$.
\end{lemma}
\proof The proof of this is basically identical to the previous lemma;
the reason for the extra case is that in the classification of the
genus there are two different types of $2$--component --- the even and
odd ones.\qed

The {\it level} of a lattice $L$ is the minimal positive integer $N$
such that $N\gamma^2/2$ is integral for all $\gamma\in L^*$.  In the
case where $L$ is elementary we see that $N$ is either square--free or
twice a square--free number.  The scalar--valued forms we use will be
forms on $\tilde\Gamma_0(N)$ (the inverse image of $\Gamma_0(N)$ in
$\lieMP_2(\RRR)$.  For these values of $N$ the cusps of
$\tilde\Gamma_0(N)$ are easy to describe

\begin{lemma}\hskip-.5ex\emph{(see \cite{shimura_automorphic})}
A complete set of representatives for the cusps is given by $1/c$ for
$c\vert N$; this cusp has width $N/(c^2,N)$.
\end{lemma}

Matrices representing the cosets $\lieMP_2(\ZZZ) / \tilde\Gamma_0(N)$
are therefore given by matrices of the form
$$\(\[\begin{array}{cc} 1&b \\ c&d\end{array}\], +\sqrt{c\tau+d}\)
\quad \hbox{for $c$ dividing $N$}.$$

Let $\({a \over b}\)$ be the Kronecker symbol, an extension of the
Jacobi symbol to all integers with
\begin{eqnarray*}
\({a \over -1}\) &=&
\begin{cases}
-1 & \hbox{for $a<0$,} \\
1 & \hbox{for $a>0$.}
\end{cases} \\
\({a \over 2}\) &=&
\begin{cases}
0 & \hbox{for $a$ even,} \\
1 & \hbox{for $a\equiv\pm 1\ (\mod 8)$,} \\
-1 & \hbox{for $a\equiv\pm 3\ (\mod 8)$.}
\end{cases}
\end{eqnarray*}
We now define some characters of $\tilde\Gamma_0(N)$.  The character
$\chi_n$ is defined by
$$\chi_n\(\[\begin{array}{cc}a&b\\c&d\end{array}\],\pm \sqrt{c\tau +
d}\) = \( {d \over n}\).$$
The character $\chi_\vartheta$ is defined by
$$\chi_\vartheta\(\[\begin{array}{cc}a&b\\c&d\end{array}\],\pm
\sqrt{c\tau + d}\) = 
\begin{cases}
\pm\({c \over d}\) & \hbox{if $d\equiv 1\ (\mod 4)$,} \\
\mp i\({c \over d}\) & \hbox{if $d\equiv 3\ (\mod 4)$.} \\
\end{cases}$$
The character $\chi_A$, where $A$ is a discriminant form, is given by
$$\chi_A = \begin{cases}
\chi_\vartheta^{\sgn(A) + \({-1\over |A|}\)-1} \chi_{|A|2^{\sgn(A)}} &
\hbox{if $4\vert N$,} \\
\chi_{|A|} & \hbox{if $4\!\nmid\! N$.}\end{cases}$$

Firstly we study how to extract scalar--valued forms from
components of the vector--valued form.

\begin{theorem}\hskip-.5ex\emph{(see \cite{borcherds_refl})}
Suppose that $L^*/L$ has level $N$.  If $b$ and $c$ are divisible by
$N$ then
$$g = \(\[\begin{array}{cc}a&b\\c&d\end{array}\], \pm\sqrt{c\tau+d}\)
\in \lieMP_2(\ZZZ)$$
acts on the Weil representation by
$$g(\ee_\gamma) = \chi_{L^*/L}(g)\ee_{a\gamma},$$
where $\chi_{L^*/L}$ is the character defined above.
\end{theorem}

From this we can show

\begin{lemma} Suppose that $L^*/L$ has level $N$.  If $\gamma\in
L^*/L$ has norm zero then
$$g = \(\[\begin{array}{cc}a&b\\c&d\end{array}\], \pm\sqrt{c\tau+d}\)
\in \tilde\Gamma_0(N)$$
acts on $\ee_\gamma$ as
$$g(\ee_\gamma) = \chi_{L^*/L}(g)\ee_{a\gamma}.$$
\end{lemma}
\proof We can write $g\in \tilde\Gamma_0(N)$ in the form
$$\(\[\begin{array}{cc} 1&k\\0&1\end{array}\],+1\) \(\[\begin{array}{cc}
a'&b'\\c&d\end{array}\],\pm\sqrt{c\tau+d}\),$$
where $b'$ and $c$ are divisible by $N$.  If $\gamma$ is a norm zero
vector in the discriminant group then it is invariant under the action
of $\smallmx{1}{k}{0}{1}$.  The character $\chi_{L^*/L}$ is trivial on
$\smallmx{1}{k}{0}{1}$.  Finally, the discriminant group has exponent
dividing $N$ and $a\equiv a'\ (\mod N)$ so that $a\gamma$ and $a'\gamma$
are identical in the discriminant group.\qed

Using this and invariance under $\Aut(L^*/L)$ we get

\begin{theorem} Let $F(\tau)$ be an $\Aut(L^*/L)$--invariant
vector--valued modular form.  Let $\gamma$ be a norm zero vector in the
discriminant group.  Then the $\ee_\gamma$--component of $F(\tau)$ is a
modular form for the group $\tilde\Gamma_0(N)$ with the same weight as
$F$ and character $\chi_{L^*/L}$.
\end{theorem}
\proof We will show that, for an elementary discriminant form, the
vectors $\gamma$ and $a\gamma$ (where $a$ is coprime to $N$) are
$\Aut(L^*/L)$--equivalent.  From this the result will clearly follow
using the previous lemma.

By decomposing the vector $\gamma$ into its Jordan components we only
need to show this result for the individual components.  If $\gamma\in
A_p$ for $p>2$ then this follows from the fact that two non--zero
vectors of $A_p$ are conjugate under $\Aut(A_p)$ if and only if they
have the same norm.  Both $\gamma$ and $a\gamma$ have norm zero and
because $(a,p)=1$ they are both simultaneously zero or nonzero.  If
$\gamma\in A_2$ then $a$ is odd and so $\gamma$ and $a\gamma$ are
identical in $A_2$.\qed

Suppose that $f(\tau)$ is a scalar--valued modular form for the group
$\tilde\Gamma_0(N)$ with character $\chi_{L^*/L}$ we will now show how
to induce this form to a vector--valued form with type $\rho_{L^*/L}$.

For $g = \(\smallmx{a}{b}{c}{d},\pm\sqrt{c\tau +d}\)$ define the
{\it slash operator} $|_g$ of weight $k$ by
$$f|_g(\tau) = \(\pm\sqrt{c\tau+d}\)^{-2k} f(g\tau).$$
So, if $f$ has character $\chi$ then $f$ satisfies
$$f|_g = \chi(g) f \quad\hbox{for all $g\in\tilde\Gamma_0(N)$}.$$

\begin{lemma} If $g,h\in\lieMP_2(\ZZZ)$ then
$$f|_g|_h = f|_{gh}.$$
\end{lemma}
\proof This is a simple computation where we must keep careful track
of the minus signs.\qed

Define the vector valued function $F(\tau)$ by
$$F(\tau)\ \ = \sum_{g\in\tilde\Gamma_0(N)\backslash\lieMP_2(\ZZZ)}
f|_g(\tau)\rho_{L^*/L}(g^{-1})\ee_0.$$
\begin{lemma} The sum defining $F(\tau)$ is well defined.
\end{lemma}
\proof Suppose we replace $g$ by $hg$ with $h\in\tilde\Gamma_0(N)$.
Then the term in the sum becomes
\begin{eqnarray*}
f|_{hg}(\tau)\rho_{L^*/L}((hg)^{-1})\ee_0
&=& \chi_{L^*/L}(h) f|_g(\tau) \rho(g^{-1}) \rho(h^{-1})
\ee_0 \\
&=& \chi_{L^*/L}(h) f|_g(\tau) \rho(g^{-1})\chi_{L^*/L}(h^{-1})
\ee_0 \\
&=& f|_g(\tau) \rho_{L^*/L}(g^{-1})\ee_0.
\end{eqnarray*}
Hence the sum is independent of coset representative.\qed

\begin{theorem}\label{lift} The function $F(\tau)$ is a vector valued
form of weight $k$ and type $\rho_{L^*/L}$.
\end{theorem}
\proof Let
$$h = \(\[\begin{array}{cc}a&b\\c&d\end{array}\], \sqrt{c\tau + d}\)
\in \lieMP_2(\ZZZ).$$
A computation gives
\begin{eqnarray*}
F(h\tau) &=& \sum_g f|_g (h\tau) \rho(g^{-1})\ee_0 \\
&=& \sum_g f|_{gh^{-1}} (h\tau) \rho(hg^{-1})\ee_0 \\
&=& (c\tau+d)^k \rho(h) \sum_g f|_g(\tau) \rho(g^{-1})\ee_0 \\
&=& (c\tau+d)^k \rho(h) F(\tau),
\end{eqnarray*}
which proves the result.\qed

It is well known that this induction is not injective --- there are
some forms, even ones with singularities, that induce to the zero
vector--valued modular form.  However, we shall show that any
$\Aut(L^*/L)$--invariant vector valued modular form on an elementary
lattice is induced in this way from a scalar--valued modular form.  In
fact, we will show that the scalar--valued form can be chosen to be a
linear combination of the scalar--valued forms from the
$\ee_\gamma$--components of $F(\tau)$, for $\gamma$ of norm zero.

Think of $F(\tau)$ as a column vector.  Let $\vecv$ be a row vector
with non-zero entries corresponding to norm zero vectors.  So, $\vecv
F(\tau)$ is a $\tilde\Gamma_0(N)$ form with character $\chi_{L^*/L}$.
When we induce this we get
\begin{eqnarray*}
&&\sum_{g\in\tilde\Gamma_0(N)\backslash\lieMP_2(\ZZZ)} f|_g(\tau)
\rho(g^{-1})\ee_0\\
&=& \sum_g (c_g\tau+d_g)^{-k} f(g\tau) \rho(g^{-1}) \ee_0 \\
&=& \sum_g (c_g\tau+d_g)^{-k} \vecv F(g\tau) \rho(g^{-1}) \ee_0 \\
&=& \sum_g (c_g\tau+d_g)^{-k} \vecv (c_g\tau+d_g)^k \rho(g)
F(\tau) \rho(g^{-1}) \ee_0 \\
&=& \sum_g \vecv\rho(g)F(\tau)\rho(g^{-1})\ee_0.
\end{eqnarray*}

Looking only at the $\ee_0$--coordinate of the induction we get (writing
$\rho_{00}$ for the top left hand entry of $\rho(g)$)
$$\sum_{g\in\tilde\Gamma_0(N)\backslash\lieMP_2(\ZZZ)} \vecv
\rho(g) F(\tau) \overline{\rho_{00}} =
\mathbf{v}\( \sum_g \overline{\rho_{00}} \rho(g) \) F(\tau).
$$

We can evaluate the sum using Shintani's formula for the matrix
coefficients of the Weil representation (Theorem \ref{shintani}).
Let $\alpha$ and $\beta$ be norm zero vectors in the discriminant
group.  We shall work out the $(\alpha,\beta)$--coefficient for
the Weil representation of $\smallmatrix{1}{b}{c}{d}$.

\begin{proposition}
$$\rho\(\[\begin{array}{cc}1&b\\c&d\end{array}\]\)_{\alpha\beta} =
{1\over |L^*/L|}
\sum_{\mu\in L^*/L} \e\( -{ c\mu^2 + 2(\alpha-\beta,\mu) \over 2}
\).$$
\end{proposition}
\proof By the Shintani formula
\begin{eqnarray*}
&& {\sqrt{i}^{\ -\sgn(L)} \over c^{n/2} \sqrt{|L^*/L|}}
\sum\limits_{\gamma\in L/cL} \e\({ a(\alpha+\gamma)^2 -
2(\beta,\alpha+\gamma) + d\beta^2 \over 2c} \) \\
&=& {\sqrt{i}^{\ -\sgn(L)} \over c^{n/2} \sqrt{|L^*/L|}}
\sum\limits_{\gamma\in L/cL} \e\({ (\alpha+\gamma)^2 -
2(\beta,\alpha+\gamma) + (1+bc)\beta^2 \over 2c} \) \\
&=& {\sqrt{i}^{\ -\sgn(L)} \over c^{n/2} \sqrt{|L^*/L|}}
\sum\limits_{\gamma\in L/cL} \e\({ (\alpha+\gamma)^2 -
2(\beta,\alpha+\gamma) + \beta^2 \over 2c} \)
\e\({b\beta^2 \over 2}\) \\
\end{eqnarray*}

\begin{eqnarray*}
&=& {\sqrt{i}^{\ -\sgn(L)} \over c^{n/2} \sqrt{|L^*/L|}}
\sum\limits_{\gamma\in L/cL} \e\({ (\alpha+\gamma)^2 -
2(\beta,\alpha+\gamma) + \beta^2 \over 2c} \) \\
&=& {\sqrt{i}^{\ -\sgn(L)} \over c^{n/2} \sqrt{|L^*/L|}}
\sum\limits_{\gamma\in L/cL} \e\({ (\alpha - \beta)^2
\over 2c}\) \e\({\gamma^2 + 2(\alpha-\beta,\gamma) \over 2c} \).
\end{eqnarray*}

To evaluate this sum we use a Fourier transform trick.  Set
$$\psi_\gamma(x) = \e\({1\over 2c}[ (\gamma+x)^2 +
2(\gamma+x,\lambda) ]\)$$
and
$$\Psi(x) = \sum_{\gamma\in L/cL} \psi_\gamma(x).$$
The sum we want to evaluate is $\Psi(0)$ and this also given (by the
Poisson summation formula) as the sum of the Fourier coefficients.  For
$\mu\in L^*$ set
$$c_\mu = \int_\F \Psi(x)e^{-2\pi i(\mu,x)} dx.$$
Thus we have $G(0)$ given by
\begin{eqnarray*}
&& {1\over \sqrt{|L^*/L|}}\sum_{\mu\in L^*} c_\mu \\
&=& {1\over \sqrt{|L^*/L|}}\sum_{\mu\in L^*/L}\sum_{\gamma\in L/cL} \int_\F
\e\( {(x+\gamma)^2 + 2(x+\gamma,\lambda) \over 2c} - 
(\mu,x+\gamma) \) dx \\
&=& {1\over \sqrt{|L^*/L|}}\sum_{\mu\in L^*} \int_{c\F}
\e\( {1 \over 2c} [(x+\lambda-c\mu)^2 - (\lambda+c\mu)^2]\) dx \\
&=& {1\over \sqrt{|L^*/L|}}\sum_{\mu\in L^*/L} \(\int_{L\tensor\RRR}
e^{\pi i x^2/c}dx\) \e\( -{1 \over 2c} (\lambda+c\mu)^2\). \\
\end{eqnarray*}

This integral is easily seen to cancel the terms at the front of the
Shintani formula.  Hence the matrix coefficient is
$${1\over |L^*/L|}
\sum_{\mu\in L^*/L} \e\( -{ c\mu^2 + 2(\alpha-\beta,\mu) \over 2}
\).$$
\qed

What we actually want is $\overline{\rho_{00}} \rho(g)$ so this is
\begin{eqnarray*}
&& {1\over |L^*/L|^2}\sum\limits_{\mu\in L^*/L} \sum\limits_{\delta\in L^*/L}
\e\( - {c\mu^2 - c\delta^2 + 2(\alpha-\beta,\mu) \over 2} \) \\
&=& {1\over |L^*/L|^2}\sum\limits_{\mu\in L^*/L} \sum\limits_{\delta\in L^*/L}
\e\( {c\mu^2 + 2(\alpha-\beta-c\mu,\delta) \over 2} \) \\
&=& {1\over |L^*/L|}\sum\limits_{\mu\in L^*/L \atop c\mu = \alpha-\beta}
\e\( {c\mu^2 \over 2} \).
\end{eqnarray*}

Notice we can evaluate this sum by evaluating it on the subsets $A_p$
of the discriminant form.  We can also restrict to the case where
$\alpha$ and $\beta$ are in $A_p$.  So we try to evaluate the sum:
\begin{equation}
{1\over |A_p|}\sum\limits_{\mu\in A_p \atop c\mu = \alpha-\beta}
e^{\pi i c\mu^2}.\label{our_sum}
\end{equation}

\begin{lemma} Suppose $p$ divides $c$. Then (\ref{our_sum}) is $1$ if
$\alpha=\beta$ and $0$ otherwise.
\end{lemma}
\proof If $p$ divides $c$ then the second condition in the sum is that
$\alpha = \beta$.  So, we get zero if $\alpha \ne \beta$.  When
$\alpha = \beta$ all terms in the sum are $1$. \qed

\begin{lemma} Suppose $p$ does not divide $c$ and let $d \equiv c^{-1}
\ (\mod p)$.  Then (\ref{our_sum}) is
$${1\over |A_p|} e^{-2\pi i d(\alpha,\beta)}.$$
\end{lemma}
\proof In this case the unique $\mu$ in the sum is given by $\mu =
d(\alpha - \beta)$. \qed

Due to the $\Aut(A_p)$--invariance we can regard these matrices as
matrices on the norm zero vectors in $A_p^{\Aut(A_p)}$.  This means we
can get matrices that are either $1\times 1$ (if there are no
non--trivial norm $0$ vectors), $2\times 2$ (if there are non--trivial
norm $0$ vectors, but no parity ones) or $3\times 3$ (if there
are non--trivial norm $0$ vectors and parity vectors).

\begin{lemma} The sum in the induction, restricted to $A_p$ is a
non--zero scalar multiple of the equivalent sum for the discriminant
form $A_p$.
\end{lemma}
\proof Just do the sum.\qed

To check that the matrix is invertible, we only need to do it for
the case when the discriminant form has only one non--zero Jordan
component.  Note also that the matrix is independent of the dimension
of the lattice $L$.

Let $A$ be a discriminant form that is an elementary $p$--group (so
it is of the form $(\ZZZ/p\ZZZ)^n$ for some $n$).  Let $F(\tau)$ be an
$\Aut(A)$--invariant vector--valued form.  Let $W$ be the vector space
spanned by the $\ee_\gamma$--coefficients of $F(\tau)$ for $\gamma$ of
norm $0$.  We know that $W$ is finite dimensional.

\begin{lemma} The induction from $W$ to vector--valued forms is
injective.
\end{lemma}
\proof If a form $f(\tau)$ induced to the zero vector--valued form
then we can deduce that the Fourier expansion of $f$ at the cusp $0$
is identically zero.  Hence $f(\tau) = 0$ and the induction is
injective. \qed

Hence, we have an endomorphism of $W$ given by inducing to a
vector--valued form and then restricting to the $\ee_0$--component.

\begin{lemma}\hskip-.5ex\emph{(see \cite{borcherds_refl})}
For $g = \(\smallmx{a}{b}{c}{d},\pm\)$ in $\lieMP_2(\ZZZ)$,
$\rho_L(g)\ee_0$ is a linear combination of vectors from $A^{c*}$.
\end{lemma}

\begin{lemma}\label{0component} If the $\ee_0$--component of $F(\tau)$
is zero then $F(\tau)=0$.
\end{lemma}
\proof By examining the transformation behaviour under matrices from
the previous lemma we see that various components of the
vector--valued form are identically zero.  In the case of an
elementary lattice there are enough cusps (choices of $c$) to show
that all such components are zero.\qed

Hence the endomorphism defined above is an automorphism.

\begin{theorem}\label{lifting}  An $\Aut(L^*/L)$--invariant
vector--valued form, for $L$ elementary, is induced from a
scalar--valued form.
\end{theorem}
\proof The above work shows that some linear combination of the norm
zero components of $F(\tau)$ will induce to an
$\Aut(L^*/L)$--invariant vector--valued form with the same
$\ee_0$--component as $F(\tau)$.  Lemma \ref{0component} shows that
this must therefore agree with $F(\tau)$. \qed

\section{Regular Discriminant Forms}

We have already seen that we can use modular forms to bound the
signature from below by $-24$.  In the case where we have a nice
scalar--valued modular form associated to the lattice we can do even
better than this.  In this section we show that for certain
discriminant forms a bound sharper than $-24$ can be placed on the
signature.  In order to simplify the arguments in this section we will
only present them for the case of no odd 2--Jordan components.  The
arguments in the excluded cases are, as always, notationally slightly
more complicated but really no more difficult.

We call a discriminant form {\it regular} if each of the
$p$--components, $A_p$, of the discriminant form contains vectors of
all norms $k/p$ and a non--trivial norm zero vector.  It is only in
the case of small discriminant forms that this can fail.

\ex In the case of odd primes, the $p$--components are always regular
if the $\prk$ is at least $3$.  Sometimes they are regular if the
$\prk$ is $2$ and they are never regular if the $\prk$ is $1$.

Regular forms behave in a much nicer way than irregular forms.  For
example, the results about inducing from a scalar form are easy in the
regular case

\begin{lemma} If $f$ is a non-zero scalar--valued modular form for the
lattice $L$, $L$ regular, then its lift to a vector--valued modular
form is non-zero.
\end{lemma}
\proof Consider the components of the induced vector--valued form
having order equal to the level of $L$.  If they are all zero then
the form $f$ is identically zero at one cusp.  This means that $f$ is
the zero function. \qed

Suppose we have a Lorentzian lattice $L$ with a Weyl vector.  As we
have seen, associated to it is a vector--valued modular form $F(\tau)$
with singularities corresponding to the primitive roots of $L$.  As we
are using a Weyl vector (rather than a generalized Weyl vector) it is
clear that $F(\tau)$ will be invariant under $\Aut(L^*/L)$.  Thus, if
we also assume that the lattice is elementary then Theorem
\ref{lifting} shows:

\begin{theorem} An elementary Lorentzian lattice $L$ with Weyl vector
is associated to a scalar--valued modular form $f(\tau)$ on
$\tilde\Gamma_0(N)$ (where $N$ is the level of $L$) with weight
$\sgn(L)/2$ and character $\chi_{L^*/L}$.
\end{theorem}

We have already seen (Lemma \ref{primitive}) that primitive roots of
$L$ have norm $-2/n$ and order dividing $n$.  This means that the
singularities in the vector--valued modular form $F(\tau)$ must be of
the form $q^{-1/n}$.  The second condition from Lemma \ref{primitive}
that $n\vecv\in L$ shows that the singularity can only occur at the
cusp $N/n$ in the scalar--valued modular form $f(\tau)$.  Hence,
$f(\tau)$ has poles of order at most $1$ at all cusps.  This is enough
of a restriction on the singularities of $f(\tau)$ for us to produce a
bound for the signature of the corresponding lattice.  This bound
appears to be close to best possible --- there are usually lattices
occuring at the predicted bound and they are almost always
interesting.  The lattices not included in this calculation (those
with irregular discriminant forms) seem to be rare.

\ex The Leech, Barnes--Wall and Coxeter--Todd lattices all occur at
exactly the critical signature.  

Let $f$ be a modular form of weight $k$ for the group $\Gamma(N)$.  We
will use the Riemann--Roch theorem to give a lower bound for the
number of poles $f$ can have.  This then gives us a way to give a
minimum (negative) signature that a regular lattice can have.  For now
we assume that the level is square-free, this corresponds to having no
odd $2$--Jordan components.  Similar arguments work in the remaining
case but the formul\ae\ are much more complicated; we will discuss the
remaining case later.

\begin{lemma} The group $\Gamma(N)$, for $N>1$, has no elliptic
points.
\end{lemma}
\proof The stabilizers of the elliptic points in $\lieSL_2(\ZZZ)$ are
clearly not congruent to the identity matrix modulo $N$.\qed

\begin{lemma} The degree of the divisor of $f$ is given by
$$\deg((f)) = {k \over 12}\[\lieSL_2(\ZZZ) : \Gamma(N)\].$$
\end{lemma}
\proof This follows from the well known formula \cite{miyake}
$$g = 1 + {\mu \over 12} - {\nu_2 \over 4} - {\nu_3 \over 3} -
{\nu_\infty \over 2},$$
where $g$ is the genus, $\mu$ is the index in $\lieSL_2(\ZZZ)$,
$\nu_2$, $\nu_3$ are the number of elliptic points of order $2$, $3$,
respectively and $\nu_\infty$ is the number of cusps.  All cusps of
$\Gamma(N)$ are of width $N$, so $\nu_\infty = \mu / N$.  Finally, the
degree of the divisor of an automorphic form is given by
$$\deg((f)) = k(g-1) + {k \over 2}\sum_a \(1 - {1\over e_a}\),$$
where $a$ runs over the inequivalent elliptic points and cusps and
$e_a$ is the order of the point $a$ (either $2$, $3$ or $\infty$).
Combining these gives the result. \qed

For the modular forms constructed above we know that the only poles
occur at the cusps and have order at most $1$.  Thus the degree of the
divisor is at least minus the number of cusps for $\Gamma_0(N)$.  As
$\Gamma(N)$ is a subgroup of $\Gamma_0(N)$ we can regard the modular
form as a $\Gamma(N)$ form with the degree of the divisor multiplied
by $\[\Gamma_0(N) : \Gamma(N)\]$.

\begin{theorem} The signature of $L$ satisfies
$$\sgn(L) \ge -24 \times {\hbox{\rm cusp}(N) \over \hbox{\rm index}(N)},$$
where $\hbox{\rm cusp}(N)$ is the number of cusps of $\Gamma_0(N)$ and
$\hbox{\rm index}(N)$ is the index of $\Gamma_0(N)$ in $\lieSL_2(\ZZZ)$.
\end{theorem}

Note that this theorem shows that there are only a finite number of
such lattices with negative signature.  This fits in well with work of
Nikulin \cite{nikulin_big} which uses his study of ``narrow parts of
hyperbolic polyhedra'' to show that there are only a finite number of
such lattices (up to an obvious equivalence relation which is roughly
the same as assuming the lattice to be elementary).  Indeed, it seems
natural to conjecture that our methods of using modular forms to study
Lorentzian reflection groups should also be able to give this
finiteness result.  Unfortunately, it is the cases of irregular
discriminants which seem to cause problems.  However, when the method
does work we obtain good bounds on the dimension which gives the hope
that one could effectively classify all such lattices.

\ex If the lattice $L$ has level $1$ then the number of cusps is $1$
and the index is $1$.  So, we obtain the signature bound of $-24$
again.  This signature corresponds to the Leech lattice.  If the
lattice has level $2$ then the number of cusps is $2$ and the index is
$3$.  So, we obtain a signature bound of $-16$ which corresponds to
the Barnes--Wall lattice.  If the level is $3$ then there are $2$
cusps and index $4$.  This gives a signature bound of $-12$ which
corresponds to the Coxeter--Todd lattice.  The existence of such a
signature is noted in \cite{borcherds_refl} where it says
\begin{quotation}\small
What usually seems to happen is that for each level there is some
``critical'' signature, with the property that almost all lattices up
to that signature have non--zero reflective modular forms, but beyond
that signature there are only a few isolated examples, usually with
$\prk$ at most $2$ for some prime $p$.  For example, for level $N=1$
the critical signature is $-24$, corresponding to the Leech lattice,
the lattice $\even_{1,25}$ whose reflection group was described by
Conway, and so on, while for level $N=2$ the critical signature is
$-16$ corresponding to the Barnes--Wall lattice and so on.
\end{quotation}

It is obvious why this signature is the beginning of a flood of such
lattices: multiplying by an Eisenstein series preserves the
singularity structure of these scalar forms and so we easily get
scalar forms corresponding to lower dimensional lattices.

Notice that in the above quotation it is stated that the exceptional
lattices seem to correspond to a small $\prk$ for some prime $p$.
This is explained in our notation by the lattice being irregular.

There are lattices with irregular discriminant forms (\eg the even
sublattice of $\odd_{1,21}$ which gives rise to the maximal finite
co--volume reflection group according to \cite{esselmann}).  These,
like the regular lattices, seem to decrease in dimension as the level
increases.

In the cases where the critical signature is attainable (in particular
it should be an integer) the lattices are especially nice.  The table
below shows the lattices that occur.  All of them are strongly modular
lattices \cite{quebbemann2}.
\begin{table}[ht]
$$
\begin{array}{|c|c|c|}
\hline
\hbox{\bf Level} & \hbox{\bf Critical Signature} &
\hbox{\bf Lattice} \\
\hline\hline
1 & -24 & \hbox{Leech $\Lambda$}\\
2 & -16 & \hbox{Barnes--Wall $BW_{16}$}\\
3 & -12 & \hbox{Coxeter--Todd $K_{12}$}\\
5 & -8 & \hbox{Icosians $Q_8(1)$}\\
6 & -8 & \hbox{$D_4 \tensor A_2$}\\
7 & -6 & \hbox{$A_6^{(2)}$}\\
11 & -4 & \hbox{11--modular}\\
14 & -4 & \hbox{strongly 14--modular}\\
15 & -4 & \hbox{strongly 15--modular}\\
23 & -2 & \hbox{23--modular}\\
\hline
\end{array}
$$
\caption{Lattices at the critical signature.}
\end{table}
Note that these numbers also turn up in \cite{extremal}, Theorem 2.1;
the levels and signatures also occur in Table 3.2, \emph{loc. cit.},
of extremal modular lattices with minimum $4$.  The lattices in the
cases of prime level are mentioned at the end of \cite{quebbemann1}.
The cases $6$, $14$ and $15$ can be found in \cite{quebbemann2}.  The
modular forms that we associate to these lattices can also be seen in
these papers.  More information about these lattices can be found in
\cite{lattice_catalogue}.

All the lattices mentioned above have the additional property that
they have minimum norm $4$.  This is also obvious from the connection
with Lorentzian reflection groups because these lattices have a norm
$0$ Weyl vector which therefore has no roots orthogonal to it.
They also share the property that they have a uniform construction for
the modular form they are associated to --- they are all
$\eta$--products and have the property that they are all inverses of a
cusp form for the group $\Gamma_0(N)+$.  We are allowed to take the
inverses of these cusp forms because they have no zeros in the finite
plane (this follows from the fact that they are $\eta$--products or by
counting zeros).  In general we can not use the inverse of a cusp form
as it will not be holomorphic in the finite plane.
\begin{table}[ht]
$$
\begin{array}{|c|c|c|}
\hline
\hbox{\bf Level} & \hbox{\bf Modular Form} & \hbox{\bf Cusp Form} \\
\hline\hline
1 & \eta(\tau)^{-24} & \Delta(\tau) \\
2 & \eta(\tau)^{-8}\eta(2\tau)^{-8} & \Delta_{2+}(\tau) \\
3 & \eta(\tau)^{-6}\eta(3\tau)^{-6} & \Delta_{3+}(\tau) \\
5 & \eta(\tau)^{-4}\eta(5\tau)^{-4} & \Delta_{5+}(\tau) \\
6 & \eta(\tau)^{-2}\eta(2\tau)^{-2}\eta(3\tau)^{-2}\eta(6\tau)^{-2}
& \Delta_{6+}(\tau) \\
7 & \eta(\tau)^{-3}\eta(7\tau)^{-3} & \Delta_{7+}(\tau) \\
11 & \eta(\tau)^{-2}\eta(11\tau)^{-2} & \Delta_{11+}(\tau) \\
14 & \eta(\tau)^{-1}\eta(2\tau)^{-1}\eta(7\tau)^{-1}\eta(14\tau)^{-1}
& \Delta_{14+}(\tau) \\
15 & \eta(\tau)^{-1}\eta(3\tau)^{-1}\eta(5\tau)^{-1}\eta(15\tau)^{-1}
& \Delta_{15+}(\tau) \\
23 & \eta(\tau)^{-1}\eta(23\tau)^{-1} & \Delta_{23+}(\tau) \\
\hline
\end{array}
$$
\caption{Modular forms associated to critical elementary lattices.}
\end{table}
From the structure of these $\eta$--products it is clear that they are
all invariant under the Atkin--Lehner involutions, which at the
lattice level corresponds to the associated lattices being invariant
under the $*m$--duals, which is the definition of strongly modular.
Denominator formul\ae\ for Borcherds--Kac--Moody algebras associated
to these forms in the prime case were studied in \cite{niemann}; a
construction based on the singular theta correspondence can be found
in the first part of \cite{scheithauer}.  In \cite{scheithauer} it was
necessary to construct by hand the vector--valued forms which are
associated to the $\eta$--products, they can also be constructed by
applying the lift described in Theorem \ref{lift} to the
$\eta$--products.  Of course, this allows us to extend the
results in \cite{niemann, scheithauer} to the $N=6,14,15$ cases.

We now examine the possible regular discriminant forms with reflective
modular forms.  Many of the modular forms occuring are
$\eta$--products, many of these products also occur in the work of
Martin \cite{martin_eta2, martin_eta1}.  Most of these
$\eta$--products are shown by Martin to be associated to elements of
the Conway group $2\hbox{\rm Co}_1$ (see \cite{atlas} for more details
about $\hbox{\rm Co}_1$).  The modular forms can also be seen in
tables in \cite{kondo_leech} where connections to Conway's group are
discussed and \cite{conway_norton} where there are connections with
the Monster.

\eject
\begin{centre}
{\large\bf Level 1}
\end{centre}

The critical signature is $-24$.  Level $1$ lattices exist only in
signatures divisible by $8$.  There are $3$ possibilities.
\begin{table}[ht]
$$
\begin{array}{|c|c|}
\hline
\hbox{\bf Lattice} & \hbox{\bf Modular Form} \\
\hline\hline
\even_{1,25} & 1/\Delta(\tau)\\
\even_{1,17} & E_4(\tau)/\Delta(\tau) \\
\even_{1,9} & E_4(\tau)^2/\Delta(\tau) \\
\hline
\end{array}
$$
\caption{Modular forms for level $1$ lattices.}
\end{table}

\vfill
\begin{centre}
{\large\bf Level 2}
\end{centre}

The critical signature is $-16$.  Level $2$ lattices exist only in
signatures divisible by $4$.  There are $4$ possible signatures.
\begin{table}[ht]
$$
\begin{array}{|c|c|}
\hline
\hbox{\bf Lattice} & \hbox{\bf Modular Form} \\
\hline\hline
\even_{1,17}(2^{+even}) & 1/\Delta_{2+}(\tau)\\
\even_{1,13}(2^{-odd}) & \theta_{D_4}(\tau)/\Delta_{2+}(\tau) \\
\even_{1,9}(2^{+even}) & \theta_{D_4}(\tau)^2/\Delta_{2+}(\tau) \\
\even_{1,5}(2^{-odd}) & \theta_{D_4}(\tau)^3/\Delta_{2+}(\tau) \\
\hline
\end{array}
$$
\caption{Modular forms for level $2$ lattices.}
\end{table}

\vfill
\eject
\begin{centre}
{\large\bf Level 3}
\end{centre}

The critical signature is $-12$.  Level $3$ lattices exist only in
signatures divisible by $2$.  There are $6$ possible dimensions.
\begin{table}[ht]
$$
\begin{array}{|c|c|}
\hline
\hbox{\bf Lattice} & \hbox{\bf Modular Form} \\
\hline\hline
\even_{1,13}(3^{even}) & 1/\Delta_{3+}(\tau)\\
\even_{1,11}(3^{odd}) & \theta_{A_2}(\tau)/\Delta_{3+}(\tau) \\
\even_{1,9}(3^{even}) & \theta_{A_2}(\tau)^2/\Delta_{3+}(\tau) \\
\even_{1,7}(3^{odd}) & \theta_{A_2}(\tau)^3/\Delta_{3+}(\tau) \\
\even_{1,5}(3^{even}) & \theta_{A_2}(\tau)^4/\Delta_{3+}(\tau) \\
\even_{1,3}(3^{odd}) & \theta_{A_2}(\tau)^5/\Delta_{3+}(\tau) \\
\hline
\end{array}
$$
\caption{Modular forms for level $3$ lattices.}
\end{table}

\vfill
\begin{centre}
{\large\bf Level 4}
\end{centre}

The critical signature is $-14$.  Level $4$ lattices exist in all
signatures.
\begin{table}[ht]
$$
\begin{array}{|c|c|}
\hline
\hbox{\bf Lattice} & \hbox{\bf Modular Form} \\
\hline\hline
\even_{1,15}(2^{even}_2) & \eta(\tau)^{-12} \eta(2\tau)^2 \eta(4\tau)^{-4}\\
\even_{1,14}(2^{odd}_3) & \eta(\tau)^{-14} \eta(2\tau)^7 \eta(4\tau)^{-6}\\
\even_{1,13}(2^{even}_4) & 1/\Delta_{4+}(\tau)\\
\even_{1,12}(2^{odd}_5) & \theta_{A_1}(\tau)/\Delta_{4+}(\tau)\\
\even_{1,11}(2^{even}_6) & \theta_{A_1}(\tau)^2/\Delta_{4+}(\tau)\\
\even_{1,10}(2^{odd}_7) & \theta_{A_1}(\tau)^3/\Delta_{4+}(\tau)\\
\vdots & \vdots \\
\hline
\end{array}
$$
\caption{Modular forms for level $4$ lattices.}
\end{table}

\vfill
\eject
\begin{centre}
{\large\bf Level 5}
\end{centre}

The critical signature is $-8$.  Level $5$ lattices exist only in
signatures divisible by $4$.  There are $2$ possible dimensions.
\begin{table}[ht]
$$
\begin{array}{|c|c|}
\hline
\hbox{\bf Lattice} & \hbox{\bf Modular Form} \\
\hline\hline
\even_{1,9}(5^{+even}) & 1/\Delta_{5+}(\tau)\\
\even_{1,5}(5^{-even}) & (E_2(\tau) - 5E_2(5\tau))/\Delta_{5+}(\tau) \\
\even_{1,5}(5^{+odd}) & \eta(\tau)\eta(5\tau)^{-5} \\
\hline
\end{array}
$$
\caption{Modular forms for level $5$ lattices.}
\end{table}

\vfill
\begin{centre}
{\large\bf Level 6}
\end{centre}

The critical signature is $-8$.  Level $6$ lattices exist only in
signatures divisible by $2$.  There are $4$ possible dimensions.
\begin{table}[ht]
$$
\begin{array}{|c|c|}
\hline
\hbox{\bf Lattice} & \hbox{\bf Modular Form} \\
\hline\hline
\even_{1,9}(2^{even}3^{even}) & 1/\Delta_{6+}(\tau) \\
\even_{1,7}(2^{even}3^{odd}) & E_1(\tau,\chi_3)/\Delta_{6+}(\tau)\\
\even_{1,5}(2^{even}3^{even}) & E_1(\tau,\chi_3)^2/\Delta_{6+}(\tau)\\
\even_{1,3}(2^{even}3^{odd}) & E_1(\tau,\chi_3)^3/\Delta_{6+}(\tau)\\
\hline
\end{array}
$$
\caption{Modular forms for level $6$ lattices.}
\end{table}

\vfill
\eject
\begin{centre}
{\large\bf Level 7}
\end{centre}

The critical signature is $-6$.  Level $7$ lattices exist only in
signatures divisible by $2$.  There are $3$ possible dimensions.
\begin{table}[ht]
$$
\begin{array}{|c|c|}
\hline
\hbox{\bf Lattice} & \hbox{\bf Modular Form} \\
\hline\hline
\even_{1,7}(7^{odd}) & 1/\Delta_{7+}(\tau) \\
\even_{1,5}(7^{even}) & E_1(\tau,\chi_7)/\Delta_{7+}(\tau)\\
\even_{1,3}(7^{odd}) & E_1(\tau,\chi_7)^2/\Delta_{7+}(\tau)\\
\hline
\end{array}
$$
\caption{Modular forms for level $7$ lattices.}
\end{table}

\vfill
\begin{centre}
{\large\bf Level 10}
\end{centre}

The critical signature is $-4$.  Level $10$ lattices exist only in
signatures divisible by $4$.  There is $1$ possible dimension.
\begin{table}[ht]
$$
\begin{array}{|c|c|}
\hline
\hbox{\bf Lattice} & \hbox{\bf Modular Form} \\
\hline\hline
\even_{1,5}(2^{even}5^{odd}) & \eta(\tau)^{-1}\eta(2\tau)^{-2}
\eta(5\tau)^{-3} \eta(10\tau)^2 \\
\even_{1,5}(2^{even}5^{even}) & \hbox{see \cite{borcherds_refl}} \\
\hline
\end{array}
$$
\caption{Modular forms for level $10$ lattices.}
\end{table}

\vfill
\eject
\begin{centre}
{\large\bf Level 11}
\end{centre}

The critical signature is $-4$.  Level $10$ lattices exist only in
signatures divisible by $2$.  There are $2$ possible dimensions.
\begin{table}[ht]
$$
\begin{array}{|c|c|}
\hline
\hbox{\bf Lattice} & \hbox{\bf Modular Form} \\
\hline\hline
\even_{1,5}(11^{even}) & 1/\Delta_{11+}(\tau) \\
\even_{1,3}(11^{odd}) & E_1(\tau,\chi_{11})/\Delta_{11+}(\tau) \\
\hline
\end{array}
$$
\caption{Modular forms for level $11$ lattices.}
\end{table}

\vfill
\begin{centre}
{\large\bf Level 12}
\end{centre}

The critical signature is $-6$.  Level $12$ lattices exist in all
signatures.
\begin{table}[ht]
$$
\begin{array}{|c|c|}
\hline
\hbox{\bf Lattice} & \hbox{\bf Modular Form} \\
\hline\hline
\even_{1,7}(2^{even}_{0/4} 3^{odd}) & 1/\Delta_{6+}(\tau) \\
\even_{1,7}(2^{even}_{2/6} 3^{even}) & \eta(\tau)^{-4}\eta(2\tau)
\eta(4\tau)^{-1}\eta(6\tau)^{-1}\eta(12\tau)^{-1} = f(\tau)\\
\even_{1,6}(2^{odd}_{1/5} 3^{odd}) & \theta_{A_1}(\tau)/\Delta_{6+}(\tau) \\
\even_{1,6}(2^{odd}_{3/7} 3^{even}) & \theta_{A_1}(\tau)f(\tau) \\
\even_{1,5}(2^{even}_{2/6} 3^{odd}) & \theta_{A_1}(\tau)^2/\Delta_{6+}(\tau) \\
\even_{1,5}(2^{even}_{0/4} 3^{even}) & \theta_{A_1}(\tau)^2\ f(\tau) \\
\vdots & \vdots \\
\hline
\end{array}
$$
\caption{Modular forms for level $12$ lattices.}
\end{table}

\vfill
\eject
\begin{centre}
{\large\bf Level 14}
\end{centre}

The critical signature is $-4$.  Level $14$ lattices exist only in
signatures divisible by $2$.  There are $2$ possible dimensions.
\begin{table}[ht]
$$
\begin{array}{|c|c|}
\hline
\hbox{\bf Lattice} & \hbox{\bf Modular Form} \\
\hline\hline
\even_{1,5}(2^{even}7^{even}) & 1/\Delta_{14+}(\tau) \\
\even_{1,3}(2^{even}7^{odd}) & E_1(\tau,\chi_{7})/\Delta_{14+}(\tau) \\
\hline
\end{array}
$$
\caption{Modular forms for level $14$ lattices.}
\end{table}

\vfill
\begin{centre}
{\large\bf Level 15}
\end{centre}

The critical signature is $-4$.  Level $15$ lattices exist only in
signatures divisible by $2$.  There are $2$ possible dimensions.
\begin{table}[ht]
$$
\begin{array}{|c|c|}
\hline
\hbox{\bf Lattice} & \hbox{\bf Modular Form} \\
\hline\hline
\even_{1,5}(3^{even}5^{even}) & 1/\Delta_{15+}(\tau) \\
\even_{1,3}(3^{odd}5^{even}) & E_1(\tau,\chi_{3})/\Delta_{15+}(\tau) \\
\even_{1,3}(3^{odd}5^{odd}) & \eta(\tau)^{-2}\eta(3\tau)
\eta(5\tau)\eta(15\tau)^{-2} \\
\hline
\end{array}
$$
\caption{Modular forms for level $15$ lattices.}
\end{table}

\vfill
\eject
\begin{centre}
{\large\bf Level 20}
\end{centre}

The critical signature is $-4$.  Level $20$ lattices exist in all
signatures.
\begin{table}[ht]
$$
\begin{array}{|c|c|}
\hline
\hbox{\bf Lattice} & \hbox{\bf Modular Form} \\
\hline\hline
\even_{1,5}(2^{even}_{0/4} 5^{even}) & 1/\Delta_{20+}(\tau) \\
\even_{1,5}(2^{even}_{0/4} 5^{odd}) & \eta(\tau)^{-2} \eta(2\tau)^{-1}
\eta(5\tau)^2 \eta(10\tau)^{-3} = f(\tau)\\
\even_{1,4}(2^{odd}_{1/5} 5^{even}) & \theta_{A_1}(\tau)/\Delta_{20+}(\tau) \\
\even_{1,4}(2^{odd}_{1/5} 5^{odd}) & \theta_{A_1}(\tau)f(\tau) \\
\vdots & \vdots \\
\hline
\end{array}
$$
\caption{Modular forms for level $20$ lattices.}
\end{table}

\vfill
\begin{centre}
{\large\bf Level 23}
\end{centre}

The critical signature is $-2$.  Level $23$ lattices exist only in
signatures divisible by $2$.  There is $1$ possible dimension.
\begin{table}[ht]
$$
\begin{array}{|c|c|}
\hline
\hbox{\bf Lattice} & \hbox{\bf Modular Form} \\
\hline\hline
\even_{1,3}(23^{odd}) & 1/\Delta_{23+}(\tau) \\
\hline
\end{array}
$$
\caption{Modular forms for level $23$ lattices.}
\end{table}

\vfill
\eject
\begin{centre}
{\large\bf Level 28}
\end{centre}

The critical signature is $-2$.  Level $28$ lattices exist in all
signatures.
\begin{table}[ht]
$$
\begin{array}{|c|c|}
\hline
\hbox{\bf Lattice} & \hbox{\bf Modular Form} \\
\hline\hline
\even_{1,3}(2^{even}_{0/4} 7^{odd}) & \eta(\tau)^{-1} \eta(2\tau)
\eta(4\tau)^{-1} \eta(7\tau)^{-1} \eta(14\tau) \eta(28\tau)^{-1} \\
\even_{1,3}(2^{even}_{2/6} 7^{even}) & \eta(\tau)^{-1}
\eta(2\tau)^{-1} \eta(7\tau)^{-5} \eta(14\tau)^9 \eta(28\tau)^{-4}\\
\vdots & \vdots \\
\hline
\end{array}
$$
\caption{Modular forms for level $28$ lattices.}
\end{table}

\vfill
\begin{centre}
{\large\bf Level 30}
\end{centre}

The critical signature is $-2$.  Level $30$ lattices exist only in
signatures divisible by $2$.  There is $1$ possible dimension.
\begin{table}[ht]
$$
\begin{array}{|c|c|}
\hline
\hbox{\bf Lattice} & \hbox{\bf Modular Form} \\
\hline\hline
\even_{1,3}(2^{even}3^{odd}5^{odd}) & \eta(\tau) \eta(3\tau)^{-1}
\eta(5\tau)^{-1} \eta(6\tau) \eta(10\tau) \eta(15\tau)\\
\hline
\end{array}
$$
\caption{Modular forms for level $30$ lattices.}
\end{table}

\vfill
\eject

As $\eta$--products seem quite common we wrote a {\tt PARI} program to
search for negative weight $\eta$--products with the correct
singularities at cusps.  The program uses a formula for the order of
zeros and poles of $\eta$--products.

\begin{lemma}\hskip-.5ex\emph{(see \cite{martin_eta1})}
Fix a level $N$.  Let $t_1, \dots, t_s$ be the divisors of $N$ and
$$f(\tau) = \eta(t_1 \tau)^{r_1} \cdots \eta(t_s \tau)^{r_s}.$$
The the order of zero of $f(\tau)$ at the cusp $a/c$ is given by
$${1\over 24} {N \over \gcd(c^2, N)} \sum_{j=1}^s {\gcd(t_j,c)^2 \over
t_j} r_j$$
\end{lemma}

This formula leads to a matrix equation for the orders of poles in
terms of the exponents $r_j$.  Inverting this gives the exponents in
terms of the orders of poles.  The program then searches through all
possible orders of poles to find ones that come from $\eta$--products
with all the $r_j \in \ZZZ$.  The results of the program can be found
in Appendix \ref{eta_products}.

\section{Irregular Discriminant Forms}

In the case of irregular discriminants it is possible to obtain
results using a technique we will call {\it level lowering} (this has
nothing to do with work of Ribet \cite{ribet}).  The idea is as
follows: Take a vector--valued form corresponding to a lattice with
irregular discriminant form.  If one can find a suitable
vector--valued form transforming under the dual of the Weil
representation on the irregular part then by tensoring and taking a
trace we obtain a vector--valued form corresponding to a regular
discriminant.  This new form will satisfy the signature bounds of the
previous section and so we will get bounds for the original form.  The
``suitable'' in the description of the level lowering form has to do
with the strengths of the singularities that occur.  It is also
possible to get results in the irregular case by examining the
obstruction spaces defined in \cite{borcherds_gkz}.  If the dimension
of the obstruction space is strictly smaller than the freedom
available there will clearly exist a reflective form; if the
obstruction space is larger then we expect the only reflective form to
be the zero function.  Of course, it is possible that the obstructions
defined by the obstruction space are not independent (this is known to
happen in at least one case) but computational evidence shows that
this situation is rare and that, in most cases, the obstruction space
is so large that, even with the dependencies removed, we should expect
enough restrictions to remain.

We begin by examining the possible irregular discriminant forms that
occur when the level is a prime (as this is the simplest case).  We
firstly show that the irregular case $p^{\pm 2}$ occurs only if $p=2$:

\begin{theorem} Let $p$ be an odd prime.  If the lattice
$\even_{1,n+1}(p^{\pm 2})$ has an irregular discriminant form, then
there is no associated modular form.
\end{theorem}
\proof Let $f(\tau)$ be a scalar--valued modular form associated to
the lattice.  The signature of the lattice is divisible by
$4$. Suppose that the signature is not $-20$.  Then we can subtract a
suitable multiple of the level $1$ form $E_{12-n/2}(\tau) /
\Delta(\tau)$ to remove the $q^{-1}$ singularity of $f(\tau)$ from the
cusp at $0$.  The regular signature bound now applies showing that
$f(\tau)$ must be this level $1$ form.  

If the signature is $-20$, consider $f(\tau)$ as a modular form for
the group $\Gamma_0(2p)$.  We can subtract a suitable multiple of the
level $2$ form $\Theta_{D_4}(\tau)/\Delta(\tau)$ to remove the
$q^{-1}$ singularities of $f(\tau)$ from two of the cusps.  The
remaining form has at most $2p+2$ poles.  The usual arguments show
that this form should have at least
$${10(3p+3) \over 12}$$
poles.  Thus, $f(\tau)$ must be this level $2$ form.

It is easy to see that the vector--valued form induced from such a
form is identically zero, which shows that there was no such
$f(\tau)$.\qed

So we can restrict to the $p^{\pm 1}$ case.  The results of
\cite{borcherds_gkz} show that the space of obstructions to getting a
vector--valued form with singularities is dual to a space of cusp
forms.  Translating this to the theory of scalar--valued forms we see
that there will not exist a reflective form for an irregular
discriminant form $\even_{1,1+n}(p^{\pm 1})$ if and only if:
\begin{enumerate}
\renewcommand{\labelenumi}{\roman{enumi}.}
\item There is a cusp form of weight $2+{n \over 2}$ for the group
$\Gamma_0(p)$ with character $\chi_p$;
\item The Fourier coefficents $a_m=0$ if $\chi_p(2m) = \mp 1$;
\item The Fourier coefficient $a_1 = 0$;
\item The Fourier coefficient $a_p \ne 0$.
\end{enumerate}
By using \cite{modular_forms} we were able to check for the existence
of such forms for $p<500$.  With slightly more work it is possible to
make the above technique work for level $4$.  These computation lead
to the following conjecture:

\begin{conjecture} There are only $3$ irregular lattices of prime
level (not covered by the regular case), these are:
$$\even_{1,21}(2^{-2}),\quad \even_{1,19}(3^{+1}),\quad
\even_{1,9}(5^{-1}).$$
If we include level $4$ we get in addition:
$$\even_{1,19}(2^{+2}_6), \quad \even_{1,18}(2^{+1}_7), \quad
\even_{1,16}(2^{+3}_5).$$
\end{conjecture}

We will give below some of the cusp forms which prevent reflective
forms from existing in some cases which clarify results from
\cite{borcherds_refl}:

\begin{enumerate}
\renewcommand{\labelenumi}{\roman{enumi}.}
\item There is no reflective form for $\even_{1,15}(3^{-1})$ since
there is a cusp form for $\Gamma_0(3)$, character $\chi_3$, weight $9$
and $q$--expansion
$$q^2 - 3q^3 - 10q^5 + 45q^6 + \cdots.$$
In \cite{borcherds_refl} the form
$$q^{-1} - 216 - 9126q + \cdots$$
is mentioned as a reflective form for the lattice
$\even_{1,15}(3^{-1})$.  In spite of the existence of this
scalar--valued form there is no corresponding vector--valued form (it
is noted in \cite{borcherds_refl} that this form induces to zero).
Also, the lattice $\even_{1,15}(3^{-1})$ does not have a Weyl vector:
The norm $0$ vector found in \cite{borcherds_refl} has inner product
$0$ with some norm $2$ simple roots and $-1$ with others; this
behaviour is not allowed by our definition of a Weyl vector.

All of the ``dihedral'' cases mentioned in \cite{borcherds_refl} seem
to behave in this way.

\item There is no reflective form for $\even_{1,17}(5^{-1})$ since
there is a cusp form for $\Gamma_0(5)$, character $\chi_5$, weight $10$
and $q$--expansion
$$q^4 - q^5 - q^6 - 18q^9 + 19q^{10} + 20q^{11} + \cdots.$$
In \cite{borcherds_refl} a reflective form for this lattice is
mentioned.  This is similar to the previous case: although a
scalar--valued form exists it induces to a zero vector--valued form.
Perhaps this explains why the lattice $\even_{1,17}(5^{-1})$ is not
very nice.
\end{enumerate}

Another way to investigate the reflective forms which exist is to look
at the dimension of the obstruction space.  This dimension can be
computed using the Selberg trace formula, as was done in
\cite{borcherds_refl}.  A couple of points about the calculations in
\cite{borcherds_refl} should be made:
\begin{enumerate}
\renewcommand{\labelenumi}{\roman{enumi}.}
\item In the first formula of Corollary 7.4 the term $\mathbf{e}(k/2)$
  should be $\mathbf{e}(jk/2)$;
\item The definition of $\delta_N$ in Lemma 7.2 requires that $N$ be
  minimal (however, $\delta_\infty$ is independent of $N$).
\end{enumerate}
A {\tt PARI} program was written implement this formula to compute the
dimensions of the $\Aut(L^*/L)$--invariant vector--valued modular
forms.  The results of this program are displayed in Appendix
\ref{obstructions}.

Examining the dimensions of the obstruction spaces we obtain Table
\ref{irregulartable}, a list of irregular lattices with reflective
forms.  Note that this table is almost certainly not complete for
reasons outlined below.

\begin{table}[p]
$$
\begin{array}{|c|c|}
\hline
\hbox{\bf Level} & \hbox{\bf Lattice} \\
\hline\hline
2 & \even_{1,21}(2^{-2}) \\
\hline
3 & \even_{1,19}(3^{+1}) \\
\hline
4 & \even_{1,19}(2^{+2}_6) \\
4 & \even_{1,18}(2^{+1}_7) \\
4 & \even_{1,16}(2^{+3}_5) \\
\hline
5 & \even_{1,9}(5^{-1}) \\
\hline
6 & \even_{1,15}(2^{-2}3^{+1}) \\
6 & \even_{1,11}(2^{\star}3^{\pm 1}) \\
\hline
10 & \even_{1,9}(2^{-2}5^{+1}) \\
\hline
12 & \even_{1,14}(2^{+1}_1 3^{+1}) \\

12 & \even_{1,13}(2^{+2}_2 3^{\pm 1}) \\
12 & \even_{1,13}(2^{+2}_6 3^{+1}) \\

12 & \even_{1,12}(2^{+1}_7 3^{+1}) \\
12 & \even_{1,12}(2^{+3}_3 3^{+1}) \\

12 & \even_{1,11}(2^{+2}_2 3^{+2}) \\
12 & \even_{1,11}(2^{+2}_6 3^{+2}) \\

12 & \even_{1,10}(2^{+1}_1 3^{\pm 1}) \\
12 & \even_{1,10}(2^{+1}_7 3^{+2}) \\
12 & \even_{1,10}(2^{+3}_1 3^{\pm 1}) \\
12 & \even_{1,10}(2^{+3}_3 3^{+2}) \\
12 & \even_{1,10}(2^{+3}_5 3^{\pm 1}) \\

\hline
\end{array}
\qquad
\begin{array}{|c|c|}
\hline
\hbox{\bf Level} & \hbox{\bf Lattice} \\
\hline\hline
12 & \even_{1,9}(2^{+2}_2 3^{\pm 1}) \\
12 & \even_{1,9}(2^{+2}_6 3^{\pm 1}) \\

12 & \even_{1,8}(2^{+1}_1 3^{\star}) \\
12 & \even_{1,8}(2^{+1}_7 3^{+1}) \\
12 & \even_{1,8}(2^{+3}_1 3^{+2}) \\
12 & \even_{1,8}(2^{+3}_3 3^{\star}) \\
12 & \even_{1,8}(2^{+3}_5 3^{\star}) \\
12 & \even_{1,8}(2^{+3}_7 3^{\pm 1}) \\
\hline
14 & \even_{1,7}(2^{-2}7^{-1}) \\
\hline
15 & \even_{1,7}(3^{+1}5^{+1}) \\
15 & \even_{1,7}(3^{-1}5^{-1}) \\
\hline
20 & \even_{1,10}(2^{+1}_7 5^{-1}) \\

20 & \even_{1,8}(2^{+1}_1 5^{-1}) \\
20 & \even_{1,8}(2^{+3}_5 5^{-1}) \\

20 & \even_{1,7}(2^{+2}_2 5^{-1}) \\
20 & \even_{1,7}(2^{+2}_6 5^{\pm 1}) \\
20 & \even_{1,7}(2^{+2}_6 5^{-2}) \\

20 & \even_{1,6}(2^{+1}_7 5^{\star}) \\
20 & \even_{1,6}(2^{+3}_3 5^{\pm 1}) \\
20 & \even_{1,6}(2^{+3}_3 5^{-2}) \\
20 & \even_{1,6}(2^{+3}_1 5^{\pm 1}) \\
\hline
\end{array}
\qquad
\begin{array}{|c|c|}
\hline
\hbox{\bf Level} & \hbox{\bf Lattice} \\
\hline\hline
21 & \even_{1,5}(3^{+1}7^{-1}) \\
21 & \even_{1,5}(3^{-1}7^{+1}) \\
21 & \even_{1,3}(3^{\star}7^{\pm 1}) \\
21 & \even_{1,3}(3^{+1}7^{\star}) \\
\hline
22 & \even_{1,3}(2^{\star}11^{\star}) \\
\hline
28 & \even_{1,6}(2^{+1}_1 7^{\pm 1}) \\
28 & \even_{1,6}(2^{+3}_5 7^{-1}) \\

28 & \even_{1,5}(2^{+2}_2 7^{\pm 1}) \\
28 & \even_{1,5}(2^{+2}_6 7^{-1}) \\

28 & \even_{1,4}(2^{+1}_7 7^{-1}) \\
28 & \even_{1,4}(2^{+3}_3 7^{\pm 1}) \\
28 & \even_{1,4}(2^{+3}_7 7^{-1}) \\
\hline
30 & \even_{1,7}(2^{-2}3^{-1}5^{+1}) \\
30 & \even_{1,7}(2^{-2}3^{+1}5^{-1}) \\
\hline
35 & \even_{1,3}(5^{+1} 7^{+1}) \\
\hline
42 & \even_{1,5}(2^{-2} 3^{+1} 7^{+1}) \\
42 & \even_{1,3}(2^{-2} 3^{+2} 7^{-1}) \\
42 & \even_{1,3}(2^{-2} 3^{-2} 7^{+1}) \\
\hline
44 & \even_{1,3}(2^{+1}_1 11^{\pm 1}) \\
44 & \even_{1,3}(2^{+2}_0 11^{-1}) \\
\hline
52 & \even_{1,2}(2^{+1}_7 13^{\pm 1}) \\
\hline
\end{array}
$$
\caption{Irregular elementary lattices.}\label{irregulartable}
\end{table}

The forms in Table \ref{irregulartable} are those guaranteed to
exist because the dimension of the obstruction space is strictly
smaller than the available freedom in choosing the singular
coefficients of the modular form.  It is, however, possible that some
of the obstructions are not independent of each other.  Indeed, the
case of level $11$ lattices shows that this does occur.  The dimension
of the obstruction space for regular level $11$ lattices of signature
$-4$ is $2$.  As the freedom available is only $2$ dimensional we
would only expect the zero form to survive.  This is clearly not true
from the results of the previous section.  What is happening here is
that the $2$ dimensional obstruction space contains a cusp form which
is perpendicular to the space of reflective forms.  So, the relevant
obstructions are really only $1$ dimensional.  The scalar--valued
version of this form is the square of the cusp form of weight $2$
which exists because $\Gamma_0(11)$ has genus $1$.

Because of the above comments it would be interesting to look more
carefully at the cases where the dimension of the obstruction space is
equal to the freedom available as some of these may also give
reflective lattices.

By using the Selberg trace formula we can easily show that the
dimension of obstructions grows much faster than the available
freedom.  If the level is $N = p_1p_2\cdots p_n$ then the number of
reflective terms is at most $2^n$.  The dimension of the obstruction
space grows like
$${k-1 \over 12}{p_1+1 \over 2}{p_2+1 \over 2}\cdots {p_n+1 \over 2} -
2^n$$
where $k$ is the weight of the obstruction space.  Actually, much
better bounds than this can be made.  It is easy to see that the
obstruction space is much larger than the freedom available in most
cases.  So, even if we remove the obstructions orthogonal to the
reflective forms, we expect there to still be enough obstructions to
prevent the existence of a reflective form.

The non--existence of these modular forms is related to the fact that
it is unusual for a modular form to have large numbers of consecutive
zero Fourier coefficients near the beginning of the Fourier expansion,
see \cite{brent_zerocoeff} for a case where it is possible to show
this.  However, there do exist modular forms with almost all
coefficients equal to $0$, for example $\eta(\tau)^{26}$.  Such
modular forms are called {\it lacunary}; more details can be found in
Serre's paper \cite{serre_lacunary}.

\begin{sidewaysfigure}[p]
\includegraphics[scale=.52]{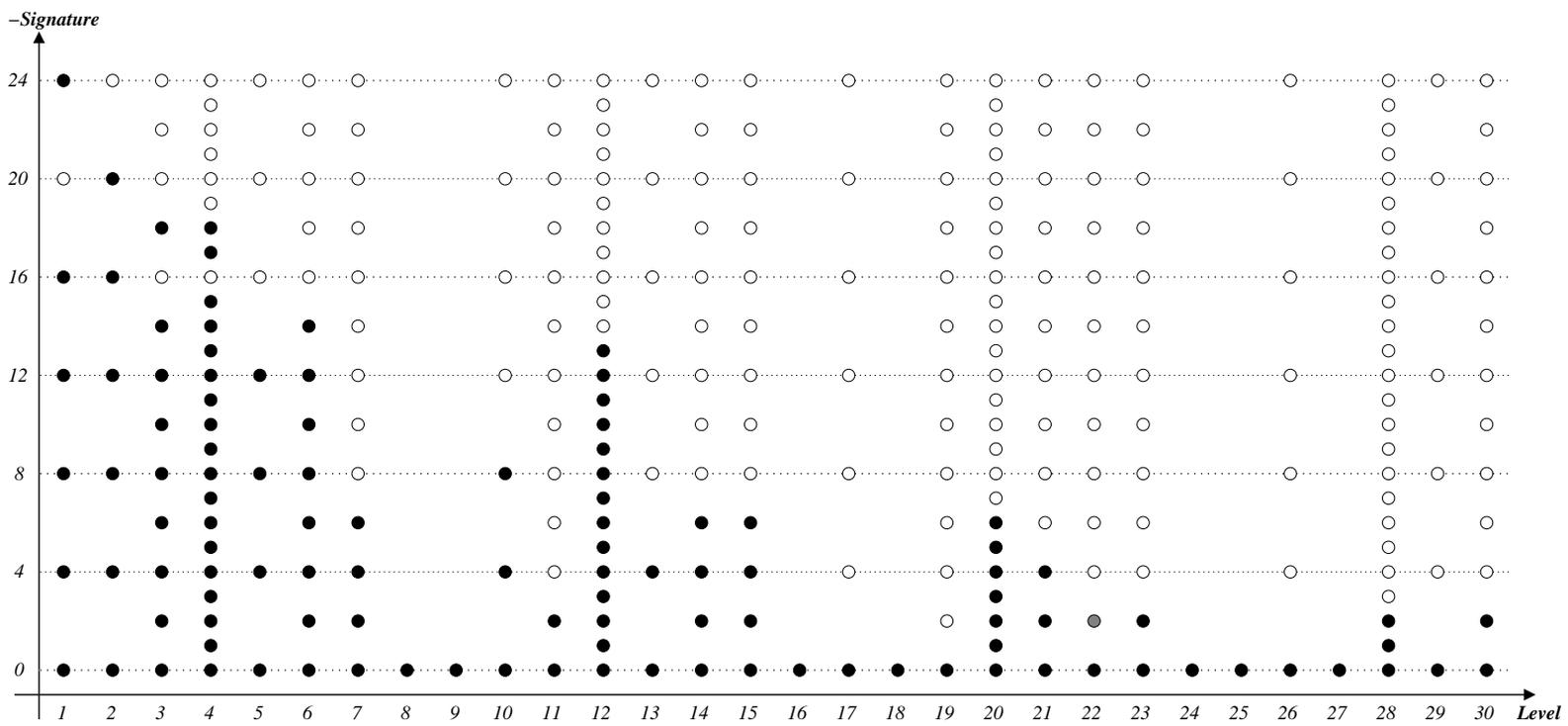}
\caption{Possible forms.}\label{possible_form}
\end{sidewaysfigure}

\begin{sidewaysfigure}[p]
\includegraphics[scale=.52]{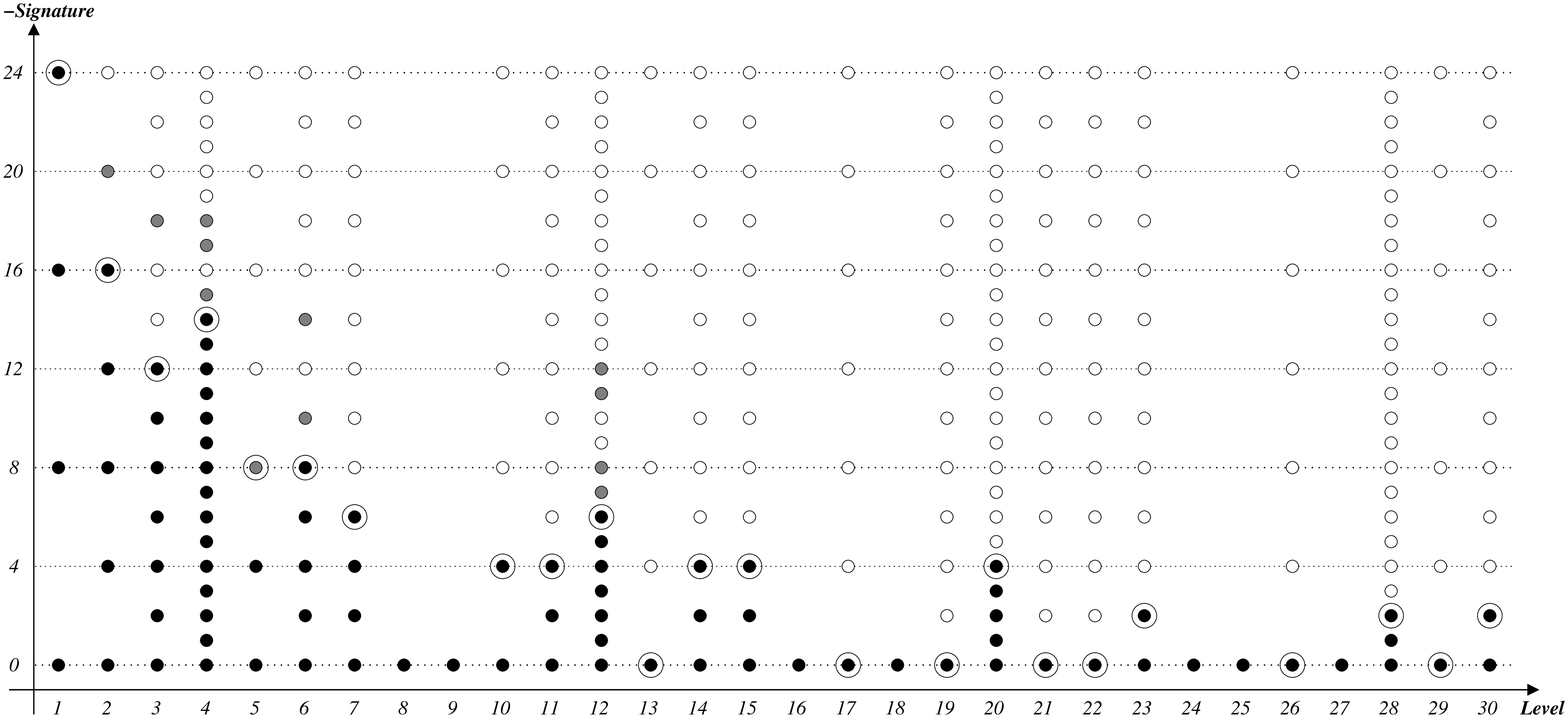}
\caption{Possible lattices.}\label{possible_lattice}
\end{sidewaysfigure}

Figure \ref{possible_form} shows graphically the dimensions and levels
in which we know there exist reflective forms.  Figure
\ref{possible_lattice} requires also that there exist a lattice
corresponding to the reflective form; the regular bound is included in
the figure.

Suppose that $L$ is a Lorentzian lattice with discriminant form $A' =
A\oplus B$ associated to a vector--valued form $F(\tau)$.  Let
$G(\tau)$ be a vector--valued form for a lattice with discriminant
$-B$.

\begin{lemma} The Weil representations on $B$ and $-B$ are dual.
\end{lemma}
\proof It is clear from the generators that changing the sign of the
inner product corresponds to taking the complex conjugate.  The result
follows since the Weil representation is unitary.\qed

We can form the vector--valued form $(F\tensor G)(\tau)$ which is a
form with values in the vector space
$$\CCC[A'] \tensor \CCC[-B] \cong \CCC[A'\oplus -B]$$
and which transforms under the representation
$$\rho_{A'\oplus -B} \cong \rho_{A'} \tensor \rho_{-B} \cong \rho_A
\tensor (\rho_B \tensor \rho_B^*).$$

\begin{lemma} Suppose that $F(\tau)$ and $G(\tau)$ are reflective
forms for lattice $\even_{1,m}(A)$ and $\even_{1,n}(B)$.  Then
$$\Delta(\tau)\cdot \( F(\tau) \tensor G(\tau) \)$$
is a reflective form for $\even_{1,m+n-25}(A\oplus B)$.
\end{lemma}
\proof The only thing to check is that the singularities are of the
correct form which is easy to see.\qed

\ex Many of the reflective forms for irregular lattices look like
$\Theta_L(\tau)/\Delta(\tau)$, for some lattice $L$.  For compound
levels the lattices which occur are often direct sums of the ones
which occur for prime levels.  This behaviour is explained by the
above lemma.

If $F(\tau)$ is a form on $A\oplus B$ and $G(\tau)$ is a form on $-B$
then we can restrict $(F\tensor G)(\tau)$ to be a form transforming
under the representation $\rho_A$ by taking the natural {\it trace}
$${\rm tr}: \rho_B \tensor \rho_B^* \longrightarrow 1.$$
We denote the resulting form by $(F\odot G)(\tau)$.

\begin{lemma} The form $(F\odot G)(\tau)$ is a vector--valued modular
form with valued in the vector space $\CCC[A]$, transforming under
the representation $\rho_A$ and with weight the sum of the weights of
$F(\tau)$ and $G(\tau)$.
\end{lemma}
\proof This is just collecting all the above facts.\qed

\nb Although it is clear that $(F\tensor G)(\tau)$ is non--zero it
seems not to be obvious that $(F\odot G)(\tau)$ is non--zero.  This,
unfortunately, causes many problems below.

We can now (almost) generalize the fact that the existence of the cusp
form $\Delta(\tau)$ bounds the signature below by $-24$.

\begin{theorem}\label{cuspformbound}
Suppose that $F(\tau)$ is a reflective form for a lattice
$\even_{1,n}(A)$.  Let $G(\tau)$ be a cusp form for $\even(-A)$ of
weight $m$.  Then, either
$$\sgn(A) \ge -2m\qquad\hbox{or}\qquad (F\odot G)(\tau) \equiv 0.$$
\end{theorem}
\proof The function $(F\odot G)(\tau)$ is non--singular and has weight
$${\sgn(A) \over 2} + m.$$
Since the function is non--singular its weight must be positive or the
function must be zero.\qed

\ex For $\even(3^{-1})$ there is a cusp form of weight $9$.  In all
other level $3$ cases there are cusp forms of weight $6$ or $7$.
Hence, the smallest signature that should occur is $-18$, which is
exactly the case.

\ex For level $2$ there is always a cusp form of weight $8$.  Hence,
the smallest signature that should occur is $-16$.  However, we know
of the lattice $\even_{1,21}(2^{-2})$.  Therefore, the reflective form
for this lattice, $F(\tau)$, and the weight $8$ cusp form, $G(\tau)$,
give an example of when $(F\odot G)(\tau) \equiv 0$.

We can now explain the method of level lowering.  As above, we have
the problem that $(F\odot G)(\tau)$ can be identically zero even if
$F(\tau)$ and $G(\tau)$ are not.  In spite of this problem the bounds
that this method produces seem to be correct.

\begin{theorem} Suppose that $F(\tau)$ is a reflective form for a lattice
$\even(A\oplus B)$ and $G(\tau)$ is a cusp form for $\even(-B)$.  Then
$$F(\tau) \odot \(G(\tau) \over \Delta(\tau)\)$$
is a reflective form on $\even(A)$.
\end{theorem}
\proof Obvious.\qed

By picking $B$ to be the irregular part corresponding to a prime the
above theorem gives a reflective form on a lattice of lower level.  We
can therefore use the bounds for lower levels to get bounds for
irregular discriminant forms.

If it were not for the possibility that the product gives zero we
would easily be able to deduce:

\begin{conjecture} The signature bound for a compound level is at least
the maximum of the bounds for its divisors.
\end{conjecture}

This conjecture seems to be consistent with the computer calculations
so far performed.  If this conjecture were true it would show that, in
order to classify elementary reflective lattices we would get a large
amount of information from the study of prime level reflective
lattices.

If $F(\tau)$ is an element of the obstruction space to some level $N$
lattice and $G(\tau)$ is a non--singular form then $(F \odot G)(\tau)$
is again an element of the obstruction space.  By using this we see:

\begin{conjecture} The signature bound for the irregular level is at
most $8$ more than the regular bound.
\end{conjecture}

The above bounds for the signatures would be enough to give a
``modular forms'' proof of Esselmann's result on cofinite Lorentzian
reflections groups.  For more details about this lattice see
\cite{borcherds_20}.

\begin{theorem}\hskip-.5ex\emph{(see \cite{esselmann})}
The second largest dimension of a Lorentzian reflection group with
Weyl vector is $22$.  In this case the lattice is (essentially)
$\even_{1,21}(2^{-2})$.  The Weyl vector has positive norm and so the
reflection group has finite co--volume.
\end{theorem}

\chapter{Singular Weight Forms}

It is well--known that any holomorphic automorphic form on
$\lieO_{2,n}(\RRR)$ either has weight $0$ (in which case it is a
constant) or weight at least $(n-2)/2$.  As objects which occur at
some theoretical minimum are normally particularly interesting we look
more at these.  The weight $(n-2)/2$ is known as the {\it singular
weight}.  If an automorphic form has singular weight $(n-2)/2$ then
all the Fourier coefficients corresponding to vectors of non--zero
norm vanish.  Conjecturally \cite{nikulin_lkma} there should only be a
finite number of singular weight forms with zeros on Heegner divisors.
In this section we study these forms in various simple cases.

The class of singular weight forms which are of particular interest
are those whose zeros (and sometimes poles) lie on sub--Grassmannians
orthogonal to primitive roots.  As we are now in signature $(2,n)$ we
can use results of Borcherds to get infinite products from the
singular weight forms.  In particular, we are especially interested in
cases when the infinite product expansion corresponds to the
denominator formula for a Borcherds--Kac--Moody ({\it BKM}) algebra.
For this to be possible we need the dimension of the simple roots to
be $1$ and so we need the automorphic form to have order $1$ zeros and
poles.  In the $(2,n)$ case, the work of Bruinier (up to some
technical conditions which we have previously mentioned) shows that
the study of these forms is basically the same as the study of certain
vector--valued modular forms (the correspondence again being by the
singular theta lift).  The singularity structure of the vector--valued
form is determined by the location of the zeros of poles of the
singular weight form.  Hence we are particularly interested in
vector--valued modular forms corresponding to signature $(2,n)$
lattice with singular coefficients corresponding to primitive roots
and occurring with coefficient $1$.  We also know the
$\ee_0$--coefficient of the vector--valued form since we know the
weight of the associated automorphic form \cite{borcherds_products,
borcherds_grass}.

\section{Eisenstein Series}

We study singular forms of this kind as follows.  In
\cite{bruinier_kuss} a vector--valued Eisenstein series which
transforms under $\rho_L^*$ is defined.  Its Fourier coefficients have
very nice properties for lattices of signature $(2,n)$, see
Proposition 14.  In particular, they are all non--positive rational
numbers.  The weight of the Eisenstein series is of weight $1+(n/2)$.
If we multiply this Eisenstein series by the vector--valued function
giving the singular weight form we, therefore, get a weight $2$
classical modular form.

\begin{lemma} Let $f(\tau)$ be a weight $2$ classical modular form,
possibly with singularities at $i\infty$.  Then the constant term in
the $q$--expansion is zero.
\end{lemma}
\proof It is easy to see that $f(\tau)d\tau$ is a $1$--form on the
curve $X(1)$.  It is well known that the integral of such a form
should be zero.  Expressed in $q$--coordinates this becomes
$$\int_{X(1)} f(q) {dq \over q} = 0.$$
Hence the residue is $0$ which, in this case, is the constant term
in the $q$--expansion of $f$.\qed

This gives us an equation relating the coefficients of the singular
terms in $F(\tau)$ to the constant term.  We can use this relationship
to severely cut down the possibilities for the singularities that can
occur.

In \cite{bruinier_kuss}, the coefficients for the vector--valued
Eisenstein series are worked out.  In particular, it is shown

\begin{theorem} Let $\gamma \in L^*$ and $n \in \ZZZ - \gamma^2/2$
with $n>0$.  The coefficient $q(\gamma,n)$ of the Eisenstein series
$E(\tau)$ of weight $k = m/2$ is equal to
$$2^{k+1} \pi^k n^{k-1} (-1)^{b^+/2} \over \sqrt{|L^*/L|} \Gamma(k)$$
times
$$\begin{cases}\displaystyle
{\sigma_{1-k}(\tilde n, \chi_{4D}) \over L(k, \chi_{4D})}
\prod_{p|2\det(S)} p^{w_p(1-2k)} N_{\gamma,n}(p^{w_p}) & 2|m,
\\ \displaystyle
{L(k-1/2,\chi_\D) \over \zeta(2k-1)} \sum_{d|f} \mu(d) \chi_\D(d)
d^{1/2-k} \sigma_{2-2k}(f/d) \prod_{p|2\det(S)} {p^{w_p(1-2k)}
N_{\gamma,n}(p^{w_p}) \over 1-p^{1-2k}} & 2\!\nmid\! m.
\end{cases}$$
\end{theorem}
Here $\sigma_k(n, \chi)$ is the twisted divisor sum
$$\sigma_k(n, \chi) = \sum_{d|n} \chi(d) d^k.$$
$N_{\gamma,n}(a)$ counts the number of vectors $v\in L/aL$ such that
the norm of $\gamma+v$ is $-2n$ modulo $a$.  $w_p$ is a number
depending on $n$ and $\gamma$ defined by
$$w_p = 1 + 2v_p(2nd_\gamma),$$
where $v_p$ is the usual $p$--adic valuation and $d_\gamma$ is the
level of $\gamma$.  $D$ and $\D$ are certain discriminants related to
the determinant of the lattice $L$.

Most of these terms are easy to compute using a computer.  A {\tt C++}
program to implement the above calculation is mentioned in
\cite{bruinier_kuss}.  Unfortunately, as we were told by Kuss
\cite{kuss_priv} the source code for this program was lost.  So, it
was necessary for us to implement a new version of this program.

\section{Lattice Gauss Sums}

One method to calculate the quantity $N_{\gamma,n}(a)$ quickly
involves evaluating Gauss sums over arbitrary lattices (see
\cite{gauss_sums} for the $\ZZZ^n$ case).  This section explains how
to evaluate such sums.

The Gauss sum we shall consider is given by
$$G_L(a,\beta,c) = {1 \over \sqrt{|L/cL|}} \sum_{r\in L/cL} \e\( -{a(r
+ \beta)^2 \over 2c}\),$$
where $a$ and $c$ are coprime integers, $\beta$ is an element of the
dual lattice and $L$ is an even lattice.

\begin{lemma} If $\beta \not\in A^{c*}$ then the Gauss sum is zero.
\end{lemma}
\proof As $\beta \not\in A^{c*}$ we can find an element $\gamma$ of
order $c$ such that
$$2c(\gamma,\beta) \not\equiv c^2 \gamma^2\ \mod 2c.$$
If, in the above sum, we send $r \mapsto r - c\gamma$ then this is
simply a rearrangement of the sum (as $c\gamma \in L$).  Now assuming
the sum is non--zero, we compute
\begin{eqnarray*}
\sum_{r\in L/cL} \e\( -{a(r + \beta)^2 \over 2c}\) &=&
\sum_{r\in L/cL} \e\( -{a(r-c\gamma + \beta)^2 \over
2c}\) \\
&=& \sum_{r\in L/cL} \e\( -{a((r + \beta)^2 + c^2\gamma^2 -
2c(\beta,\gamma)) \over 2c}\) \\
&\ne& \sum_{r\in L/cL} \e\( -{a(r + \beta)^2 \over 2c}\).
\end{eqnarray*}
For the inequality in the final line note that $2c(\gamma,\beta) -
c^2\gamma^2$ is an integer not divisible by $2c$ and $a$ and $c$ are
coprime.  Hence, the sum must be zero.\qed

We have the obvious property
$$G_L(a,\beta,c) = G_{L(a)} (1,\beta,c).$$
Thus we can restrict to the case when $a=1$.  In this case we can
convert the sum to be over $L^*/L$ by a Fourier transformation.

\begin{theorem}
$$G_L(1,\beta,c) = {\sqrt{i}^{\ -\sgn(L)} \over \sqrt{|L^* / L|}}
\sum_{\mu\in L^*/L} \e \( c\mu^2/2 + (\beta,\mu) \).$$
\end{theorem}
\proof Define
$$\psi(x) = \e\( -{(x + \beta)^2 \over 2c} \).$$
This function is clearly periodic under translations by the lattice
$cL$.  Set
$$\Psi(x) = \sum_{r \in L/cL} \psi(x+r).$$
This function is periodic under translations by $L$.  The value of the
Gauss sum is basically given by $\Psi(0)$.  As $\Psi(x)$ is a
periodic function on $L\tensor\RRR$ we can define its Fourier
coefficients: For $\mu\in L^*$ set
$$c(\mu) = \int_{L\tensor\RRR / L} \Psi(x) e^{-2\pi i(\mu,x)} dx.$$
Now we compute (by the Poisson summation formula)
\begin{eqnarray*}
\sqrt{|L^*/L|}\ \Psi(0) &=& \sum_{\mu \in L^*} c(\mu) \\
&=& \sum_{\mu\in L^*}\sum_{r\in L/cL} \int_{L\tensor\RRR / L}
\e\( -{(x+r+\beta)^2 + \over 2c} - (\mu,x+r) \) dx \\
&=& \sum_{\mu\in L^*} \int_{L\tensor\RRR / cL}
\e\( -{(x+\beta+c\mu)^2 - (c\mu)^2 -2c(\beta,\mu) \over 2c} \) dx \\
&=& \sum_{\mu\in L^*/L} \(\int_{L\tensor\RRR} e^{-\pi i x^2/c} dx \)
\e \( c\mu^2/2 + (\beta,\mu) \) \\
&=& \sqrt{|L/cL|} \sqrt{i}^{\ -\sgn(L)} \sum_{\mu\in L^*/L} \e \(
c\mu^2/2 + (\beta,\mu) \).
\end{eqnarray*}
Putting this all together we get the result.\qed

We can split the sum over $L^*/L$ into a product over local factors, so
we only need to compute these local factors.  Let $A$ be the
discriminant form for a local factor.  Assume that the local factor
corresponds to an odd prime, then

\begin{lemma} The Gauss sum is zero if $\beta \not\in A^c$ and is
$$G_L(1,\beta,c) = {\sqrt{i}^{\ -\sgn(A)} \over \sqrt{|A|}}
\e\(-{c\gamma^2 \over 2}\)\sum_{\mu \in A} \e\( {c\mu^2/2}\)$$
if $\beta = c\gamma$.
\end{lemma}
\proof We have already shown that the Gauss sum is zero unless $\beta
\in A^{c*}$ and if $A$ is a local factor corresponding to an odd prime
then $A^{c*} = A^c$.  So, we can assume that $\beta = c\gamma$ for
some $\gamma \in A$.  In the sum we can now substitute $\mu \mapsto
\mu - \gamma$.  This is a rearrangement and so does not affect the
sum.  Evaluating this gives the lemma. \qed

The sum can now be computed by Milgram's formula to get
$$G_L(1,\beta,c) = \delta\!\[ \beta \in A^c \] \sqrt{i}^{\ \sgn(A^c) -
\sgn(A)} \sqrt{|A_c|}\ \e\(-{c \gamma^2 \over 2}\),$$
where $\gamma$ is some $c$-th root of $\beta$ and
$$\delta\!\[ \beta \in A^c \] = \begin{cases}
1 & \hbox{if $\beta\in A^c$,} \\
0 & \hbox{otherwise.}
\end{cases}$$

If $A$ corresponds to the even prime and it is an even Jordan
component then we get exactly the same as above.  If it corresponds to
an odd component then we get the sum
$${\sqrt{i}^{\ -\sgn(A)} \over \sqrt{|A|}} \sum_{\mu \in A} \e\(
{c\mu^2/2 + (\beta,\mu)}\).$$
If the power of $2$ dividing $c$ is different from that dividing the
order of $A$ then $\beta = 0$, otherwise $\beta \ne 0$.  The first
case can be evaluated as above.  For the second note that $(\beta,
\mu) \equiv c\mu^2/2$ for all $\mu \in A_c$.  And, $A_c = A$ by the
condition on powers of $2$.  It is easy now to see that the sum gives
$|A|$.  So, we always get

\begin{theorem} If $A = L^*/L$ is the discriminant form for $L$ then
$$G_L(1,\beta,c) = \delta\[ \beta \in A^{c*}\]
\sqrt{|A_c|}\ \sqrt{i}^{\ \sgn(A(a)^c) - \sgn(A)}.$$
\end{theorem}
\proof We need to find the relationship between the discriminant forms
for $L$ and $L(a)$.  As $a$ and $c$ are coprime it is clear that the
number of elements of order $c$ will be the same in both.  So, we only
need to worry about the signature of $A^c$.\qed

If we know the genus of the lattice then it is easy to work out the
genus of the scaled lattice $L(a)$ and hence what the discriminant
form is.  This allows us to work out the signature of $A^c$ easily.  A
computer program was implemented in {\tt PARI} to do these
calculations.  It is listed in Appendix \ref{programs}.

Using these Gauss sums we can rapidly evaluate the numbers
$N_{\gamma,n}(m)$ required in the formula for the Eisenstein series
coefficients.  By definition
$$N_{\gamma,n}(m) = \#\{ \lambda \in L/mL : (\lambda-\gamma)^2/2 + n
\equiv 0\ \mod m\}.$$
This can be computed by summing Gauss sums
\begin{lemma}
$$N_{\gamma,n}(m) = {1\over m} \sum_{\lambda\in L/mL} \sum_{k=1}^m
\e\(\[ {(\lambda-\gamma)^2 \over 2} + n\]k/m\).$$
Hence
$$N_{\gamma,n}(m) = {\sqrt{|L/mL|} \over m} \sum_{k=1}^m
G_L(k,-\gamma, m)\ \e\( {nk \over m}\).$$
\end{lemma}
\proof The first formula follows from the fact that sums of powers of
roots of unity is zero unless the root of unity is $1$.\qed

So, we can compute the numbers $N_{\gamma,n}(m)$ by summing over
something $m$ times the size of the discriminant group.  This is
considerably less effort than searching all of $L/mL$ for vectors with
the correct norm modulo $m$.

\section{Singular Weight Forms}

Now that we have a computer program that calculates the Fourier
coefficients of the Eisenstein series we implement the following
algorithm
\begin{enumerate}
\item Fix an (elementary) level $N$;
\item Find all lattices of level $N$, signature $(2,n)$ and signature
  $\ge -24$.  This is easy to do by simply computing what the possible
  genera are;
\item For each possible lattice, computer the Fourier coefficients of
  the Eisenstein series corresponding to the $q^{-1/n}$ terms and
  constant term;
\item See if it is possible to solve the resulting equation with
  constant coefficent the singular weight and singular coefficients
  $1$ or $0$;
\item Output the lattices which pass this test;
\item For each lattice which passes try to find a singular weight
  form.
\end{enumerate}

A sample {\tt PARI} program to do this is listed in Appendix
\ref{programs} (it was modified in the obvious ways for different
levels).

Currently we have only used this program in the case where there are
no odd Jordan $2$--components.  We intend to implement this soon and
expect to find many more lattices with singular weight forms.  Our
results can be found in Appendix \ref{singular_weight}.  Notice
that there is a strong correlation between the singular weight forms
and the $\eta$--products listed in Appendix \ref{eta_products}.

\nb When the calculations in this section were performed we did not
have bounds on the signature as good as those obtained in the previous
chapter.  As such, the {\tt PARI} program tried many primes and low
signatures.

One very interesting thing can be seen from the results in Appendix
\ref{singular_weight}:  all the singular weight forms we find are
related to elements of Conway's group.  The $\eta$--products which
turn up are the ``shapes'' of these elements which may partly explain
why we see so many $\eta$--products.  Some further calculations were
performed (but the results are not listed in this dissertation) which
suggest that this pattern continues even in the non--elementary case.
So, it seems like the singular weight automorphic forms are strongly
related to elements of Conway's group.

We briefly recall the construction of the fake monster Lie algebra in
\cite{borcherds_fake}.  Given a non--singular even lattice $L$ there
is an associated vertex algebra $V_L$.  By the Sugawara construction
there is an action of the Virasoro algebra on $V_L$, where the
Virasoro algebra is spanned by the operators $1$ and $L_i$, $i\in
\ZZZ$ with
$$\[L_i, L_j\] = (i-j)L_{i+j} + {i^3 - i \over 12} \dim(L)
\delta_{i+j}$$
Let $P^i$ be the vectors in $V_L$ which have eigenvalue $i$ under
$L_0$ and are annihilated by $L_n$ for $n<0$.  Then $P^1 / L_1(P^0)$
can be made into a Lie algebra with the bracket and bilinear form
coming from the vertex operator structure of $V_L$.  In the case where
$L$ is $\even_{1,25}$ the no ghost theorem \cite{goddard_thorn}
implies that the bilinear form is positive semidefinite and so we can
take the quotient by its kernel.  This gives the fake monster Lie
algebra.  There are various other construction of vertex algebras and
associated Lie algebras which seem to correspond to other elements of
Conway's group: for example, the element $2A$ should give the ``fake
baby monster Lie algebra''.

There is a construction in the physical literature known as the
``orbifold construction'' which we recall here.  Suppose $V$ is a
vertex algebra with only one simple module and $G$ is a group of
automorphisms acting on $V$.  The fixed points of $V$ under $G$ form
another vertex algebra with only finitely many simple modules.
Conjecturally, it should be possible to endow the collection of these
modules with a vertex algebra structure.  This is the orbifold
construction.  It has been proved rigorously in only a few cases (for
example, when the group of automorphisms is $\ZZZ/2\ZZZ$).

Given the results displayed in Appendix \ref{singular_weight} and
various other computer experiments we conjecture:

\begin{conjecture} All simple BKM algebras with denominator formula a
singular weight automorphic form are obtained from the fake monster
Lie algebra by the orbifolding construction.
\end{conjecture}

A result similar to this would explain the finiteness of the set of
such BKM algebras, the appearance of Conway's group and
``shapes'' of its elements.  Given the already known connections
between Conway's group and genus $0$ congruence subgroups of
$\lieSL_2(\ZZZ)$ this could also explain the strange genus $0$
properties noticed in \cite{borcherds_refl}.

As there are only a finite number of such genus $0$ subgroups
\cite{thompson_genus0} this would explain Nikulin's conjectured
finiteness result for such automorphic forms.  It would also give
large amounts of information about the levels that could occur since
these genus $0$ groups are well studied.  For example, the only prime
levels that could occur would be primes dividing the order of the
Monster,
$$2^{46}\cdot 3^{20}\cdot 5^9\cdot 7^6\cdot 11^2\cdot 13^3\cdot
17\cdot 19\cdot 23\cdot 29\cdot 31\cdot 41\cdot 47\cdot 59\cdot 71,$$
by a well known theorem of Ogg \cite{ogg}.  For these primes all
the groups $\Gamma_0(p)+$ are genus $0$.

Algebras similar to BKM algebras have been studied by Nikulin
\cite{nikulin_lkma}; these algebras, known as {\it Lorentzian
Kac--Moody Algebras}, are more general because they do not have the
restriction to singular weight.  However, many of the interesting
example do seem to be of singular weight so, perhaps, not much is lost
by this restriction.

\section{Arithmetic Mirror Symmetry}

In \cite{mirror_K3} and \cite{mirror_CY} Gritsenko and Nikulin
formulate the {\it arithmetic mirror symmetry conjecture}.  At the
level of lattices this conjecture relates reflective lattices of
signature $(1,n)$ and lattices which occur in BKM algebras in
signature $(2,n)$.

\begin{conjecture}\hskip-.5ex\emph{(Arithmetic Mirror Symmetry)}
\begin{enumerate}
\renewcommand{\labelenumi}{\roman{enumi}.}
\item Let $L$ be a reflective lattice of signature $(2,n)$ and $\Phi$
  a reflective automorphic form for $L$ with primitive root system
  $\Delta(\Phi)$.  Let $c\in L$ be a primitive isotropic element of
  $L$ and $K = c^\perp / \ZZZ c$ be the corresponding hyperbolic
  lattice with the root system
$$\Delta(\Phi)|_K = \Delta(\Phi)\cap c^\perp\ \mod [c]$$
Then the root system $\Delta(\Phi)|_K \subset K$ is elliptic or
parabolic (in particular, the lattice $K$ is reflective) if the set of
roots $\Delta(\Phi)|_K$ is non--empty.
\item Any hyperbolic reflective lattice $K$ with an elliptic or
  parabolic primitive root system $\Delta$ may be obtained from some
  lattice $L$ with signature $(2,n)$ and a reflective automorphic form
  $\Phi$ of $L$ by the construction (i) above.
\end{enumerate}
\end{conjecture}

Roughly speaking this conjectures predicts that automorphic
forms with zeros along divisors corresponding to primitive roots
should give rise to nice Lorentzian reflection groups at the ``cusps''
and all such Lorentzian reflection groups occur this way.

The fact that these two objects (automorphic forms on signature
$(2,n)$ lattice and reflection groups on signature $(1,n)$ lattices)
are both related to modular forms on $\lieMP_2(\ZZZ)$ gives some new
insight into this conjecture.

Let $L$ be a lattice of signature $(2,n)$.  Suppose that $\Phi(\nu)$
is an automorphic form on $\Gr(L)$ with its zeros orthogonal to
primitive roots of $L$.  By results of Bruinier, there should be a
modular form $F(\tau)$ for the lattice $L$ whose singular theta
transformation gives $\Phi(\nu)$.  Let $z$ and $z'$ be elements of
$L^*$ which represent a cusp of $\Gr(L)$.  Define
$$K = L \cap z^\perp \cap z'^\perp \cong (L \cap z^\perp)/\ZZZ z$$
If is possible to restrict the modular form $F(\tau)$ to a form
$F|_K(\tau)$ for the lattice $K$ (see \cite{borcherds_grass} for
details about this construction).  Taking the singular theta
transformation of $F|_K$ we obtain a Weyl vector for $K$ and hence $K$
is a reflective lattice.

Suppose that $K$ is a reflective lattice of signature $(1,n)$.
Results in this dissertation show that there should be a modular form
$F(\tau)$ for the lattice $K$ whose singular theta transformation
gives the Weyl vector.  Define $L = K \oplus \even_{1,1}$.  The
modular form $F(\tau)$ can trivially be regarded as a form for $L$.
Taking the singular theta transformation of $F$ we obtain an
automorphic form $\Phi(\nu)$ on $\Gr(L)$ with its zeros orthogonal
to primitive roots of $L$.

Of course, if we assume that the results in this dissertation and
Bruinier's work hold in greater generality this gives the arithmetic
mirror symmetry conjecture for many more lattices.


\bibliography{bib_thesis}

\appendix
\chapter{Computer Programs}\label{programs}

The following is the {\tt PARI} \cite{PARI98} program used to compute the
coefficients of the Bruinier--Kuss Eisenstein series.

\begin{verbatim}
\\
\\ Procedures
\\

\\ Compute the signature of the Jordan block [p^power, dim]
block_sig(p, power, dim) =
{
  local(sig);

  sig = Mod(0,8);
  if(p>2, sig = 1 - Mod(p,8)^power;\
    sig = sig * abs(dim);\
    if(dim<0 && Mod(power,2)==1, sig = sig+4)\
  );
  if(p==2 && dim<0 && Mod(power,2)==1, sig = Mod(4,8));

  return(sig);
}


\\ Compute the signature of the lattice
lattice_sig() =
{
  local(sig);

  sig = Mod(0,8);
  for(l=1,size-1, sig += block_sig(genus[l,1],\
    genus[l,2], genus[l,3]));

  return(sig);
}


\\ Compute the signature of L scaled by "a" and shrunk by "c"
modified_sig(a,c) =
{
  local(sig, p, power, dim);

  sig = Mod(0,8);
  for(l=1, size-1,\
    p = genus[l,1];\
    power = genus[l,2];\
    dim = genus[l,3];\

    power += valuation(a,p);\
    power -= valuation(c,p);\

    if(power>0,\
      if(Mod(dim,2)==1,\
        dim *= kronecker(remove(p,a*c),p);\
      );\
      sig += block_sig(p, power, dim);\
    );\
  );

  return(sig);
}


\\ Compute the determinant of the lattice
lattice_det() =
{
  local(det);

  det = 1;
  for(l=1,size-1, det *= genus[l,1]^abs(genus[l,2]*genus[l,3]));
  if(Mod(b_minus,2) == 1, det = -det);

  return(det);
}


\\ Compute the size of A_c
Ac_order(c) =
{
  local(ord);

  ord = 1;
  for(l=1,size-1, ord *= gcd(c, genus[l,1]^genus[l,2])\
    ^abs(genus[l,3]));

  return(ord);
}


\\ Remove powers of "a" from "b"
remove(a,b) = b / a^valuation(b,a)


\\ Fill in the 1-components of the genus
fill_genus() =
{
  local(l_sign, l_det, new_size);

  new_size=size;

  fordiv(det, p,\
    if(isprime(p,1),\
      l_sign = kronecker(remove(p,det),p);\
      l_dim = dim;\
      for(l=1,size-1,\
        if(Mod(genus[l,1],p)==0,\
          l_sign *= sign(genus[l,3]);\
          l_dim  -= abs(genus[l,3]);\
      ));\
      genus[new_size,] = [p, 0, l_sign*l_dim];\
      printp("  Inserted Jordan component " genus[new_size,]);\
      new_size++;\
  ));

  size = new_size;
}


\\ Compute the Gauss sum
gauss(a,c) =
{
  local(result, temp, vec_norm);

  temp = gcd(a,c);
  a /= temp;
  c /= temp;
  result = sqrt(temp)^dim;
  result *= sqrt(Ac_order(c));

  temp  = modified_sig(a,c);
  temp -= modified_sig(a,1);
  temp  = lift(temp);

  result *= exp(2*Pi*I*temp/8);

  vec_norm = 1;
  for(l=1,size-1,\
    if(valuation(vec[l,2]/c , genus[l,1])>=0,\
      vec_norm *= exp(-2*Pi*I*a*c*vec[l,1]*\
        lift(Mod(vec[l,2]/c, genus[l,1]^genus[l,2]))^2)\
    , vec_norm=0\
  ));

  result *= vec_norm;

  return(result);
}


\\ Compute the n-th Bernoulli polynomial
bernoulli(n,x) =
{
  local(result);

  result = 0;
  for(l=0,n,\
    result += binomial(n,l) * bernfrac(l) * x^(n-l);\
  );

  return(result);
}


\\ Compute the sigma function
char_sigma(n,k,d) =
{
  local(result);

  result = 0;
  fordiv(n,l,\
    result += kronecker(d,l) * l^k;\
  );

  return(result);
}


\\ Count vectors in L/mL
count_vectors(n,m) =
{
  local(result);

  result = 0;
  for(l=1, m,\
    result += gauss(l,m) * exp(-Pi*I*l*n/m);\
  );

  result *= m^(dim/2 - 1);
\\  result = round(result);

  return(result);
}


\\ Compute the value of the L-function
Lfunction(s,d) =
{
  local(result, temp, conductor, modulus, discriminant);

  conductor = 1;
  modulus = 1;

  result = -(2*Pi)^s / s!;

  fordiv(d,p,\
    if(isprime(p,1),\
      temp = valuation(d,p);\
      if(temp>0 && Mod(temp,2)==1, conductor *= p);\
      if(temp>0, modulus *= p);\
  ));

  discriminant = sign(d)*conductor;

  fordiv(d,p,\
    if(isprime(p,1),\
      result *= 1 - kronecker(discriminant,p)/p^s;\
  ));

  temp = 0;
  for(l=1, conductor,\
    temp += kronecker(discriminant,l) * bernoulli(s, l/conductor);\
  );

  result *= temp;

  temp = 0;
  for(l=1, conductor,\
    temp += kronecker(discriminant,l) * exp(2*Pi*I*l/conductor);\
  );

  temp *= I^(-s) + kronecker(discriminant,-1)*I^s;

  result /= temp;

  return(result);
}



\\ Compute the coefficient of the Eisenstein series
coefficient(n) =
{
  local(result, discriminant, wp, k);

  k = dim/2;
  discriminant = 4*det;
  if(Mod(k,2)==1, discriminant = -discriminant);

  result = (2*Pi)^k * n^(k-1) / (k-1)!;
  result /= sqrt(abs(det));
  result *= char_sigma(n*level^2,1-k,discriminant);
  result /= Lfunction(k,discriminant);

  fordiv(2*det, p,\
    if(isprime(p,1),\
      wp = 1 + 2*valuation(2*n*level,p);\
      result *= count_vectors(2*n, p^wp);\
      result *= p^(wp*(1-2*k));\
  ));

  return(result);
}


\\ Work out the norm of the L*/L vector
vector_norm() =
{
  local(vec_norm);

  vec_norm=0;
  for(l=1,size-1,vec_norm += vec[l,1] * vec[l,2]^2);
  vec_norm = frac(vec_norm);

  return(vec_norm);
}


\\ Work out the level of the L*/L vector
level() =
{
  local(result);

  result = 1;
  for(l=1,size-1,\
    result = lcm(result, genus[l,1]^genus[l,2] /\
    gcd(genus[l,1]^genus[l,2] , vec[l,2]));\
  );

  return(result);
}


\\
\\ Main programme
\\

\\ Significant digits
\p 50;

\\ Number of components in genus
size = 1;

\\ Matrix to hold genus information
genus=matrix(100,3);
genus[size,1] = 7; genus[size,2] = 1; genus[size,3] = +5; size++;

\\ Matrix to hold vector in L*/L
vec=matrix(100,2);
vec[1,1] = 1/3; vec[1,2] = 0;

\\ Signature of the lattice
b_plus = 2; b_minus = 8;

\\ Dimension of the lattice
dim = b_plus+b_minus;

\\ Determinant of the lattice
det = lattice_det;

printp("Genus gives signature : " lattice_sig);
printp("Actual signature is   : " Mod(b_plus-b_minus,8));
printp("Dimension is          : " dim);
printp("Determinant is        : " det);
printp("Half vector norm is   : " vector_norm);
printp("Vector level is       : " level);
print;
printp("Inserting Jordan components...");

fill_genus;
\end{verbatim}

The following is a {\tt PARI} program to search for lattices which
pass the Eisenstein series test for the existence of a singular weight
form.  This program looks for lattices of level $pq$ for $p$ and $q$
distinct primes.  It is obvious how to modify it to search over
different levels.
\begin{verbatim}
\\
\\ User defined functions
\\

\\ Loop over all possible weights
allw() = for(w=1,12,all(w));

\\ Loop over a choice of level pq
all(w) = 
{
  printp("Weight = "w"\n");
  for(l=9,20,\
  for(m=l+1,20,\
    printp(prime(l)" "prime(m));\
    go(w,prime(l),prime(m));\
  ));\
}

\\ Loop over possible discriminant forms
go(wt,pp,qq) =
{
local(dd, s, a, b, c, d, x, y, test, found);

dd = 4+2*wt;
s = -2*wt;

found=0;

\\ Pick x,y such that y/p + x/q = 1/pq
x = lift(Mod(-pp,qq)^-1);
y = lift(Mod(-qq,pp)^-1);

for(l=-dd+1,dd-1,\
for(m=-dd+1,dd-1,\
if((l*m != 0),\
if(lattice(2,2-s,[[pp,1,l],[qq,1,m]],,0)==0,\

  \\ Compute Eisenstein coefficients
  a=coefficient(1);\

  coset([[(pp-1)/pp,1],[(qq-1)/qq,0]]);\
  b=coefficient(1/pp)*count_block_vectors((pp-1)/pp,pp,1,l);\

  coset([[(pp-1)/pp,0],[(qq-1)/qq,1]]);\
  c=coefficient(1/qq)*count_block_vectors((qq-1)/qq,qq,1,m);\

  coset([[y/pp,1],[x/qq,1]]);\
  d=coefficient(1/(pp*qq))*count_block_vectors(y/pp,pp,1,l)*\
    count_block_vectors(x/qq,qq,1,m);\

  \\ Possible singular coefficients
  test=Set([a,b,c,d,a+b,a+c,a+d,b+c,b+d,c+d,a+b+c,a+b+d,a+c+d,\
    b+c+d,a+b+c+d]);

  \\ Try to find choice of singular coefficients giving
  \\ singular weight form; display results
  if(setsearch(test,s)!=0,\
    printp1("II_{2,"(2-s)"}("pp"^"l" "qq"^"m"): ");\
    printp1(a" "b" "c" "d);\
    printp(" **********");\
    found++;\
  );\

))));

return(found);
}
\end{verbatim}

Finally, we have the program used to implement the Selberg trace
formula for the dimensions of the space of obstructions to the
existence of reflective forms.

\begin{verbatim}
\\
\\ Functions
\\

\\ Work out the order of the matrix
ord(x)=
{
  local(n, id);

  id = matid(matsize(x)[1]);

  n = 0;
  until(norml2(x^n - id)<1/1000, n++);

  return(n);
}

\\ Compute the function delta_N(X,g)
delta(n, x, g)=
{
  local(result, j);

  result = 0;
  for(j=1,n-1, result += trace((x^j)*g)/(1-exp(2*Pi*I*j/n)));

  return(result/n);
}

\\ Compute the function delta_\infty(X,g)
delta_inf(x, g, n)=
{
  local(result);

  result = trace(g)/(2*n);
  result += delta(n,x,g);

  return(result);
}

\\ Compute the psi function
Psi(g,k,rho_t,rho_s,n)=
{
  local(result);

  result = ((k-1)/12)*trace(g);
  result += delta_inf(rho_t^-1,g,n);
  result += delta(2,exp(2*Pi*I*k/(2*2))*rho_s,g);
  result += delta(3,exp(2*Pi*I*k/(2*3))*rho_t*rho_s,g);

  return(result);
}

\\ Compute the dimension of the modular forms
dim_M(rho_t, rho_s, k, n)=
{
  local(rho_z, result, j);

  rho_t = conj(rho_t~);
  rho_s = conj(rho_s~);

  rho_z = rho_s^2;

  result = 0;
  for(j=0,3, result += exp(2*Pi*I*j*k/2)*\
                  Psi(rho_z^j,k,rho_t,rho_s,n));

  result = result/4;
  if(abs(result - round(result))>1/1000, print("** ERROR **"));

  result=round(result);

  if(result==0, result=-1);

  return(result);
}

\\ Compute the dimension of the cusp forms
dim_S(rho_t, rho_s, k, n)=
{
  local(result, eisen,j);

  result = dim_M(rho_t,rho_s,k,n);

  eisen=0;
  for(j=1, matsize(rho_t)[1], eisen+=(rho_t[j,j]==1));

  if(result!=-1, result-=eisen);

  return(result);
}

\\ Make a matrix
make(p,s)=
{
  local(result);

  if(p==2 && s==-2, result = [[1,0;0,-1] , (-1/2)*[1,3;1,-1],\
			Mod([0,1],2)]);

  if(p==2 && s==+2, result = [[1,0,0;0,1,0;0,0,-1] ,\
			(+1/2)*[1,2,1;1,0,-1;1,-2,1], Mod([0,0,1],2)]);

  if(p==2 && s==-4, result = [[1,0,0;0,1,0;0,0,-1] ,\
			(-1/4)*[1,5,10;1,-3,2;1,1,-2], Mod([0,0,1],2)]);

  if(p==2 && s==+4, result = [[1,0,0;0,1,0;0,0,-1] ,\
			(+1/4)*[1,9,6;1,1,-2;1,-3,2], Mod([0,0,1],2)]);

  if(p!=2,\
    if(abs(s)==1, result = makep1(p,s));\
    if(abs(s)>=2, result = makep2(p,s));\
  );

  return(result);
}

make2odd(s,t)=
{
  local(result);

  if(abs(s)==1, result = make2odd1(s,t));
  if(abs(s)==2, result = make2odd2(s,t));
  if(abs(s)==3, result = make2odd3(s,t));

  return(result);
}

\\ Make a p^{+/-1} matrix
makep1(p,s)=
{
  local(result,q,rt,rs,rn,j,k);

  if(abs(s)!=1, error("Only call me with signature +/-1"));
  if(p==2, error("Don't call me with the prime 2"));

  q=2;\
  until(kronecker(q,p)==s, q+=2);\

  rt = matrix((p+1)/2,(p+1)/2);\
  rs = matrix((p+1)/2,(p+1)/2,j,k,1);\
  rn = vector((p+1)/2);\

  for(j=1,(p+1)/2, rt[j,j] = exp(Pi*I*(j-1)^2*q/p));\

  for(j=1,(p+1)/2, rn[j] = Mod((j-1)^2*q/2,p));\

  for(j=1,(p+1)/2, \
    for(k=2,(p+1)/2, rs[j,k] = \
            2*real(exp(-2*Pi*I*(j-1)*(k-1)*q/p))));\

  result = [rt,rs/(2*trace(rt)-1),rn];

  return(result);
}

\\ Make a p^{+/-2} matrix 
makep2(p, s)=
{
  local(rt,rs,rn,j,j2,k,k2,size,mx1,mx2,j3,k3);

  if(abs(s)<2, error("Only call me with signature at least +/-2"));

  mx1 = make(p,+1);
  mx2 = make(p,sign(s)*(abs(s)-1));

  size = p;
  if(2*kronecker(-1,p)==s || abs(s)>2, size++);

  rt = matrix(size,size);
  rs = matrix(size,size);
  rn = vector(size);

  for(j=1,p, rt[j,j] = exp(2*Pi*I*(j-1)/p));
  if(size==p+1, rt[p+1,p+1] = 1);

  for(j=1,p, rn[j] = Mod(j-1,p));
  if(size==p+1, rn[p+1] = Mod(0,p));

  for(j=1,matsize(mx1[2])[1], \
    for(k=1,matsize(mx1[2])[2], \
      for(j2=1,matsize(mx2[2])[1], \
        for(k2=1,matsize(mx2[2])[2], \
          j3=lift(mx1[3][j]+mx2[3][j2]);\
          k3=lift(mx1[3][k]+mx2[3][k2]);\
          if(j3==0 && (j!=1 || j2!=1), j3=p);\
          if(k3==0 && (k!=1 || k2!=1), k3=p);\
          rs[j3+1,k3+1] += mx1[2][j,k]*mx2[2][j2,k2];\
  ))));

  for(j=1,size,\
    j2 = rs[1,1]/rs[j,1];\
    for(k=1,size,\
      rs[j,k] *= j2;\
  ));

  result = [rt,rs,rn];\

  return(result);
}

\\ Make a 2^{+/-1}_t matrix
make2odd1(s,t)=
{
  local(rho_t, rho_s, rho_n);

  if(Mod(t,2)!=1, error("Only call me with odd signatures"));
  if(kronecker(t,2) != s, error("Signature is invalid"));

\\ Norm t/4 characteristic vector
  rho_t = [1,0;0,exp(2*Pi*I*t/4)];
  rho_s = [1,1;1,exp(-2*Pi*I*t/2)]/trace(rho_t);
  rho_n = [Mod(0,4), Mod(t,4)];

  return([rho_t, rho_s, rho_n]);
}

\\ Make a 2^{+/-2}_t matrix
make2odd2(s,t)=
{
  local(result);

  result = 0;

  if(Mod(t,2)!=0, error("Only call me with even signatures"));

  if((s==+2 && Mod(t,8)==0) || (s==-2 && Mod(t,8)==4),\
\\ Norm 0 characteristic vector
    result = [[1,0,0,0;0,1,0,0;0,0,I,0;0,0,0,-I],\
              [1,1,1,1;1,1,-1,-1;1,-1,-1,1;1,-1,1,-1]/2,\
              [Mod(0,4), Mod(0,4), Mod(1,4), Mod(3,4)]]);

  if((s==+2 && Mod(t,8)==2) || (s==-2 && Mod(t,8)==6),\
\\ Norm 1/2 characteristic vector
    result = [[1,0,0;0,I,0;0,0,-1],\
              [1,2,1;1,0,-1;1,-2,1]/(2*I),\
              [Mod(0,4), Mod(1,4), Mod(2,4)]]);

  if((s==+2 && Mod(t,8)==6) || (s==-2 && Mod(t,8)==2),\
\\ Norm 3/4 characteristic vector
    result = [[1,0,0;0,-I,0;0,0,-1],\
              [1,2,1;1,0,-1;1,-2,1]/(-2*I),\
              [Mod(0,4), Mod(3,4), Mod(2,4)]]);

  if(result == 0, error("Signature is invalid"));

  return(result);
}

\\ Make a 2^{+/-3}_t matrix
make2odd3(s,t)=
{
  local(result);

  if(Mod(t,2)!=1, error("Only call me with odd signatures"));

\\ Norm 1/4 characteristic vector
  if((s==+3 && Mod(t,8)==1) || (s==-3 && Mod(t,8)==5),\
    result = [[1,0,0,0,0,0;0,1,0,0,0,0;0,0,I,0,0,0;\
               0,0,0,-1,0,0;0,0,0,0,-I,0;0,0,0,0,0,I],\
              [1,2,2,1,1,1;1,0,0,-1,-1,1;1,0,0,-1,1,-1;\
               1,-2,-2,1,1,1;1,-2,2,1,-1,-1;1,2,-2,1,-1,-1]\
              /(2*(1+I)),\
              [Mod(0,4), Mod(0,4), Mod(1,4), Mod(2,4),\
               Mod(3,4), Mod(1,4)]]);

\\ Norm 3/4 characteristic vector
  if((s==+3 && Mod(t,8)==3) || (s==-3 && Mod(t,8)==7),\
    result = [[1,0,0,0;0,I,0,0;0,0,-1,0;0,0,0,-I],\
              [1,3,3,1;1,1,-1,-1;1,-1,-1,1;1,-3,3,-1]/(2*(I-1)),\
              [Mod(0,4), Mod(1,4), Mod(2,4), Mod(3,4)]]);

\\ Norm 1/4 characteristic vector
  if((s==+3 && Mod(t,8)==5) || (s==-3 && Mod(t,8)==1),\
    result = [[1,0,0,0;0,I,0,0;0,0,-1,0;0,0,0,-I],\
              [1,1,3,3;1,-1,3,-3;1,1,-1,-1;1,-1,-1,1]/(-2*(I+1)),\
              [Mod(0,4), Mod(1,4), Mod(2,4), Mod(3,4)]]);

\\ Norm 3/4 characteristic vector
  if((s==+3 && Mod(t,8)==7) || (s==-3 && Mod(t,8)==3),\
    result = [[1,0,0,0,0,0;0,1,0,0,0,0;0,0,I,0,0,0;\
               0,0,0,-1,0,0;0,0,0,0,-I,0;0,0,0,0,0,-I],\
              [1,2,1,1,2,1;1,0,-1,-1,0,1;1,-2,-1,1,2,-1;\
               1,-2,1,1,-2,1;1,0,1,-1,0,-1;1,2,-1,1,-2,-1]\
              /(2*(1-I)),\
              [Mod(0,4), Mod(0,4), Mod(1,4), Mod(2,4),\
               Mod(3,4), Mod(3,4)]]);


  return(result);
}

\\ Tensor product two matrices
tensor(mx1, mx2)=
{
  local(j,j2,k,k2,rs,rt,rn,s1,s2);

  s1 = matsize(mx1[1])[1];
  s2 = matsize(mx2[1])[1];

  rt = matrix(s1*s2,s1*s2);
  rs = matrix(s1*s2,s1*s2);
  rn = vector(s1*s2);

  for(j=1,s1,\
    for(k=1,s1,\
      for(j2=1,s2,\
        for(k2=1,s2,\
          rt[j+(j2-1)*s1, k+(k2-1)*s1] = mx1[1][j,k]*mx2[1][j2,k2];\
          rs[j+(j2-1)*s1, k+(k2-1)*s1] = mx1[2][j,k]*mx2[2][j2,k2];\
          rn[j+(j2-1)*s1] = chinese(mx1[3][j],mx2[3][j2]);\
  ))));

  return([rt,rs,rn]);
}

\\ Display results
test(p,s)=
{
  local(j, result, rho, rho_t, rho_s,eisen,last);

  rho=make(p,s);

  rho_t = rho[1];
  rho_s = rho[2];

  eisen = 0;
  for(j=1,matsize(rho_t)[1], if(rho_t[j,j]==1, eisen++));

  for(j=5,36,\
    result=dim_S(rho_t,rho_s,j/2);\
    if(result!=-1, last=j; print("Signature "4-j"\t Dimension \t"\
       result"\t Approx \t"matsize(rho_t)[1]*(j/2-1)/12-eisen+0.0)));

  j=last+968;
  result=dim_S(rho_t,rho_s,j/2);\
  if(result!=-1, last=j; print("Signature "4-j"\t Dimension \t"\
     result"\t Approx \t"matsize(rho_t)[1]*(j/2-1)/12-eisen+0.0));

  return(rho[3]);
}

\\ Display results
testfast(p,s)=
{
  local(j, result, rho, rho_t, rho_s,eisen,last);

  rho=make(p,s);

  rho_t = rho[1];
  rho_s = rho[2];

  eisen = 0;
  for(j=1,matsize(rho_t)[1], if(rho_t[j,j]==1, eisen++));

  j=5;
  last=0;
  while(last==0,\
    result=dim_S(rho_t,rho_s,j/2);\
    if(result!=-1, last=j; print("Signature "4-j"\t Dimension \t"\
       result"\t Approx \t"matsize(rho_t)[1]*(j/2-1)/12-eisen+0.0));\
    j++;
  );

  forstep(j=last+4,36,4,\
    result=dim_S(rho_t,rho_s,j/2);\
    if(result!=-1, last=j; print("Signature "4-j"\t Dimension \t"\
       result"\t Approx \t"matsize(rho_t)[1]*(j/2-1)/12-eisen+0.0)));

  j=last+968;
  result=dim_S(rho_t,rho_s,j/2);\
  if(result!=-1, last=j; print("Signature "4-j"\t Dimension \t"\
     result"\t Approx \t"matsize(rho_t)[1]*(j/2-1)/12-eisen+0.0));

  return(rho[3]);
}

\\ Display results
test2(rho)=
{
  local(j, result, rho_t, rho_s, last, eisen);

  rho_t = rho[1];
  rho_s = rho[2];

  eisen = 0;
  for(j=1,matsize(rho_t)[1], if(rho_t[j,j]==1, eisen++));

  for(j=5,36,\
    result=dim_S(rho_t,rho_s,j/2);\
    if(result!=-1, last=j; print("Signature "4-j"\t Dimension \t"\
       result"\t Approx \t"matsize(rho_t)[1]*(j/2-1)/12-eisen+0.0)));

  j=last+968;
  result=dim_S(rho_t,rho_s,j/2);\
  if(result!=-1, last=j; print("Signature "4-j"\t Dimension \t"\
     result"\t Approx \t"matsize(rho_t)[1]*(j/2-1)/12-eisen+0.0));

  return(rho[3]);
}

\\ Display results
test2fast(rho,n)=
{
  local(j, result, rho_t, rho_s, last, eisen, sig);

  rho_t = rho[1];
  rho_s = rho[2];

  for(j=1,matsize(rho[3])[2], print(lift(rho[3][j])/\
                               component(rho[3][j],1) - 1));

  sig = Mod(round(arg(rho_s[1,1])*-4/Pi),8);
  print("Signature: "sig"\n");

  eisen = 0;
  for(j=1,matsize(rho_t)[1], if(rho_t[j,j]==1, eisen++));

  sig=Mod(lift(sig),4);

  forstep(j=5+lift(-sig-1),28,4,\
    result=dim_S(rho_t,rho_s,j/2,n);\
    if(result!=-1, last=j; print("Signature "4-j"\t"\
       lift(Mod(4-j,8))"\t Dimension \t"result"\t Approx \t"\
       matsize(rho_t)[1]*(j/2-1)/12-eisen+0.0))\
  );

  j=last+968;
  result=dim_S(rho_t,rho_s,j/2,n);\
  if(result!=-1, last=j; print("Signature "4-j\
     "\t Dimension \t"result"\t Approx \t"\
     matsize(rho_t)[1]*(j/2-1)/12-eisen+0.0));
}

\p 4
\end{verbatim}

\chapter{Eta Products}\label{eta_products}

We list the results of a {\tt PARI} program to find $\eta$--products
with negative weight and singularities at cusps corresponding to roots
of an associated lattice.  The results are given for elementary
lattices only which corresponds to the level being square--free or
twice a square--free number.  In the case where there are odd
$2$--Jordan blocks, which is exactly when the level is divisible by
$4$, the $\eta$--products can be allowed to have double poles at
certain cusps; these forms are marked with a dagger ($\dagger$).

\begin{table}[ht]
$$
\begin{array}{|c||c|c|c|}
\hline
\hbox{\bf Level} & \hbox{\bf Shape} & \hbox{\bf Poles (1/n)} & \hbox{\bf
Weight} \\
\hline
\hline
1 & 1^{-24} & 1 & -12 \\
\hline
2 & 1^{-8} 2^{-8} & 1,2 & -8 \\
2 & 1^{-16} 2^{8} & 1 & -4 \\
2 & 1^{8} 2^{-16} & 2 & -4 \\
\hline
3 & 1^{-6} 3^{-6} & 1,3 & -6 \\
3 & 1^{-9} 3^{3} & 1 & -3 \\
3 & 1^{3} 3^{-9} & 3 & -3 \\
\hline
\end{array}
$$
\caption{$\eta$--products with reflective singularities.}
\end{table}

\addtocounter{table}{-1}
\begin{table}[p]
$$
\begin{array}{|c||c|c|c|}
\hline
\hbox{\bf Level} & \hbox{\bf Shape} & \hbox{\bf Poles (1/n)} & \hbox{\bf
Weight} \\
\hline
\hline
4 & 1^{-8} 2^{-8} & 1,2,4^\dagger & -8 \\
4 & 1^{-12} 2^{2} 4^{-4} & 1,2,4^\dagger & -7 \\
4 & 1^{-14} 2^{7} 4^{-6} & 1,2,4^\dagger & -13/2 \\
4 & 2^{-12} & 1,2,4 & -6 \\
4 & 1^{-8} 2^{-12} 4^{8} & 1,2^\dagger & -6 \\
4 & 1^{-16} 2^{12} 4^{-8} & 1,4^\dagger & -6 \\
4 & 1^{-18} 2^{17} 4^{-10} & 1,4^\dagger & -11/2 \\
4 & 1^{-4} 2^{-2} 4^{-4} & 1,2,4 & -5 \\
4 & 1^{-12} 2^{-2} 4^{4} & 1,2^\dagger & -5 \\
4 & 1^{-20} 2^{22} 4^{-12} & 1,4^\dagger & -5 \\
4 & 1^{-6} 2^{3} 4^{-6} & 1,2,4 & -9/2 \\
4 & 1^{-14} 2^{3} 4^{2} & 1,2^\dagger & -9/2 \\
4 & 1^{-22} 2^{27} 4^{-14} & 1,4^\dagger & -9/2 \\
4 & 1^{-8} 2^{8} 4^{-8} & 1,4 & -4 \\
4 & 1^{-24} 2^{32} 4^{-16} & 1,4^\dagger & -4 \\
4 & 1^{-8} 2^{-16} 4^{16} & 1,2 & -4 \\
4 & 2^{-16} 4^{8} & 1,2 & -4 \\
4 & 1^{8} 2^{-16} & 2,4 & -4 \\
4 & 1^{-16} 2^{8} & 1^\dagger & -4 \\
4 & 1^{-10} 2^{13} 4^{-10} & 1,4 & -7/2 \\
4 & 1^{-26} 2^{37} 4^{-18} & 1,4^\dagger & -7/2 \\
4 & 1^{-18} 2^{13} 4^{-2} & 1^\dagger & -7/2 \\
4 & 1^{-4} 2^{-6} 4^{4} & 1,2 & -3 \\
4 & 1^{4} 2^{-6} 4^{-4} & 2,4 & -3 \\
4 & 1^{-12} 2^{-6} 4^{12} & 1,2^\dagger & -3 \\
\hline
\end{array}
$$
\caption{$\eta$--products with reflective singularities. (\emph{ctd.})}
\end{table}

\addtocounter{table}{-1}
\begin{table}[p]
$$
\begin{array}{|c||c|c|c|}
\hline
\hbox{\bf Level} & \hbox{\bf Shape} & \hbox{\bf Poles (1/n)} & \hbox{\bf
Weight} \\
\hline
\hline
4 & 1^{-12} 2^{18} 4^{-12} & 1,4 & -3 \\
4 & 1^{-28} 2^{42} 4^{-20} & 1,4^\dagger & -3 \\
4 & 1^{-20} 2^{18} 4^{-4} & 1^\dagger & -3 \\
4 & 1^{-6} 2^{-1} 4^{2} & 1,2 & -5/2 \\
4 & 1^{2} 2^{-1} 4^{-6} & 2,4 & -5/2 \\
4 & 1^{-14} 2^{23} 4^{-14} & 1,4 & -5/2 \\
4 & 1^{-14} 2^{-1} 4^{10} & 1,2^\dagger & -5/2 \\
4 & 1^{-30} 2^{47} 4^{-22} & 1,4^\dagger & -5/2 \\
4 & 1^{-22} 2^{23} 4^{-6} & 1^\dagger & -5/2 \\
4 & 1^{-8} 2^{-20} 4^{24} & 1,2^\dagger & -2 \\
4 & 2^{-20} 4^{16} & 1,2 & -2 \\
4 & 1^{16} 2^{-20} & 2,4 & -2 \\
4 & 1^{-32} 2^{52} 4^{-24} & 1,4^\dagger & -2 \\
4 & 1^{-16} 2^{28} 4^{-16} & 1,4 & -2 \\
4 & 1^{8} 2^{-20} 4^{8} & 2 & -2 \\
4 & 1^{-16} 2^{4} 4^{8} & 1^\dagger & -2 \\
4 & 1^{-24} 2^{28} 4^{-8} & 1^\dagger & -2 \\
4 & 1^{-8} 2^{4} & 1 & -2 \\
4 & 2^{4} 4^{-8} & 4 & -2 \\
4 & 1^{-18} 2^{33} 4^{-18} & 1,4 & -3/2 \\
4 & 1^{-34} 2^{57} 4^{-26} & 1,4^\dagger & -3/2 \\
4 & 1^{-10} 2^{9} 4^{-2} & 1 & -3/2 \\
4 & 1^{-2} 2^{9} 4^{-10} & 4 & -3/2 \\
4 & 1^{-18} 2^{9} 4^{6} & 1^\dagger & -3/2 \\
4 & 1^{-26} 2^{33} 4^{-10} & 1^\dagger & -3/2 \\
\hline
\end{array}
$$
\caption{$\eta$--products with reflective singularities. (\emph{ctd.})}
\end{table}

\addtocounter{table}{-1}
\begin{table}[p]
$$
\begin{array}{|c||c|c|c|}
\hline
\hbox{\bf Level} & \hbox{\bf Shape} & \hbox{\bf Poles (1/n)} & \hbox{\bf
Weight} \\
\hline
\hline
4 & 1^{-4} 2^{-10} 4^{12} & 1,2 & -1 \\
4 & 1^{12} 2^{-10} 4^{-4} & 2,4 & -1 \\
4 & 1^{-12} 2^{-10} 4^{20} & 1,2^\dagger & -1 \\
4 & 1^{-36} 2^{62} 4^{-28} & 1,4^\dagger & -1 \\
4 & 1^{-20} 2^{38} 4^{-20} & 1,4 & -1 \\
4 & 1^{4} 2^{-10} 4^{4} & 2 & -1 \\
4 & 1^{-4} 2^{14} 4^{-12} & 4 & -1 \\
4 & 1^{-12} 2^{14} 4^{-4} & 1 & -1 \\
4 & 1^{-20} 2^{14} 4^{4} & 1^\dagger & -1 \\
4 & 1^{-28} 2^{38} 4^{-12} & 1^\dagger & -1 \\
4 & 1^{-6} 2^{-5} 4^{10} & 1,2 & -1/2 \\
4 & 1^{10} 2^{-5} 4^{-6} & 2,4 & -1/2 \\
4 & 1^{-14} 2^{-5} 4^{18} & 1,2^\dagger & -1/2 \\
4 & 1^{-38} 2^{67} 4^{-30} & 1,4^\dagger & -1/2 \\
4 & 1^{-22} 2^{43} 4^{-22} & 1,4 & -1/2 \\
4 & 1^{2} 2^{-5} 4^{2} & 2 & -1/2 \\
4 & 1^{-22} 2^{19} 4^{2} & 1^\dagger & -1/2 \\
4 & 1^{-14} 2^{19} 4^{-6} & 1 & -1/2 \\
4 & 1^{-6} 2^{19} 4^{-14} & 4 & -1/2 \\
4 & 1^{-30} 2^{43} 4^{-14} & 1^\dagger & -1/2 \\
\hline
5 & 1^{-4} 5^{-4} & 1,5 & -4 \\
5 & 1^{-5} 5^{1} & 1 & -2 \\
5 & 1^{1} 5^{-5} & 5 & -2 \\
\hline
6 & 1^{-2} 2^{-2} 3^{-2} 6^{-2} & 1,2,3,6 & -4 \\
6 & 1^{-1} 2^{-4} 3^{-5} 6^{4} & 1,2,3 & -3 \\
\hline
\end{array}
$$
\caption{$\eta$--products with reflective singularities. (\emph{ctd.})}
\end{table}

\addtocounter{table}{-1}
\begin{table}[p]
$$
\begin{array}{|c||c|c|c|}
\hline
\hbox{\bf Level} & \hbox{\bf Shape} & \hbox{\bf Poles (1/n)} & \hbox{\bf
Weight} \\
\hline
6 & 1^{-4} 2^{-1} 3^{4} 6^{-5} & 1,2,6 & -3 \\
6 & 1^{-5} 2^{4} 3^{-1} 6^{-4} & 1,3,6 & -3 \\
6 & 1^{4} 2^{-5} 3^{-4} 6^{-1} & 2,3,6 & -3 \\
6 & 2^{-6} 3^{-8} 6^{10} & 1,2,3 & -2 \\
6 & 1^{-6} 3^{10} 6^{-8} & 1,2,6 & -2 \\
6 & 1^{-8} 2^{10} 6^{-6} & 1,3,6 & -2 \\
6 & 1^{10} 2^{-8} 3^{-6} & 2,3,6 & -2 \\
6 & 1^{-3} 2^{-3} 3^{1} 6^{1} & 1,2 & -2 \\
6 & 1^{-4} 2^{2} 3^{-4} 6^{2} & 1,3 & -2 \\
6 & 1^{-7} 2^{5} 3^{5} 6^{-7} & 1,6 & -2 \\
6 & 1^{5} 2^{-7} 3^{-7} 6^{5} & 2,3 & -2 \\
6 & 1^{2} 2^{-4} 3^{2} 6^{-4} & 2,6 & -2 \\
6 & 1^{1} 2^{1} 3^{-3} 6^{-3} & 3,6 & -2 \\
6 & 1^{1} 2^{-8} 3^{-11} 6^{16} & 1,2,3 & -1 \\
6 & 1^{-8} 2^{1} 3^{16} 6^{-11} & 1,2,6 & -1 \\
6 & 1^{-11} 2^{16} 3^{1} 6^{-8} & 1,3,6 & -1 \\
6 & 1^{16} 2^{-11} 3^{-8} 6^{1} & 2,3,6 & -1 \\
6 & 1^{-2} 2^{-5} 3^{-2} 6^{7} & 1,2 & -1 \\
6 & 1^{-5} 2^{-2} 3^{7} 6^{-2} & 1,2 & -1 \\
6 & 1^{-3} 3^{-7} 6^{8} & 1,3 & -1 \\
6 & 1^{-7} 2^{8} 3^{-3} & 1,3 & -1 \\
6 & 1^{-9} 2^{6} 3^{11} 6^{-10} & 1,6 & -1 \\
6 & 1^{-10} 2^{11} 3^{6} 6^{-9} & 1,6 & -1 \\
6 & 1^{6} 2^{-9} 3^{-10} 6^{11} & 2,3 & -1 \\
6 & 1^{11} 2^{-10} 3^{-9} 6^{6} & 2,3 & -1 \\
\hline
\end{array}
$$
\caption{$\eta$--products with reflective singularities. (\emph{ctd.})}
\end{table}

\addtocounter{table}{-1}
\begin{table}[p]
$$
\begin{array}{|c||c|c|c|}
\hline
\hbox{\bf Level} & \hbox{\bf Shape} & \hbox{\bf Poles (1/n)} & \hbox{\bf
Weight} \\
\hline
\hline
6 & 1^{-2} 2^{7} 3^{-2} 6^{-5} & 3,6 & -1 \\
6 & 1^{7} 2^{-2} 3^{-5} 6^{-2} & 3,6 & -1 \\
6 & 2^{-3} 3^{8} 6^{-7} & 2,6 & -1 \\
6 & 1^{8} 2^{-7} 6^{-3} & 2,6 & -1 \\
6 & 1^{-6} 2^{3} 3^{2} 6^{-1} & 1 & -1 \\
6 & 1^{3} 2^{-6} 3^{-1} 6^{2} & 2 & -1 \\
6 & 1^{2} 2^{-1} 3^{-6} 6^{3} & 3 & -1 \\
6 & 1^{-1} 2^{2} 3^{3} 6^{-6} & 6 & -1 \\
\hline
7 & 1^{-3} 7^{-3} & 1,7 & -3 \\
\hline
10 & 1^{-1} 2^{-2} 5^{-3} 10^{2} & 1,2,5 & -2 \\
10 & 1^{-2} 2^{-1} 5^{2} 10^{-3} & 1,2,10 & -2 \\
10 & 1^{-3} 2^{2} 5^{-1} 10^{-2} & 1,5,10 & -2 \\
10 & 1^{2} 2^{-3} 5^{-2} 10^{-1} & 2,5,10 & -2 \\
\hline
11 & 1^{-2} 11^{-2} & 1,11 & -2 \\
\hline
12 & \# \sim 1600 & & \\
\hline
14 & 1^{-1} 2^{-1} 7^{-1} 14^{-1} & 1,2,7,14 & -2 \\
14 & 1^{-2} 2^{1} 7^{-2} 14^{1} & 1,7 & -1 \\
14 & 1^{1} 2^{-2} 7^{1} 14^{-2} & 2,14 & -1 \\
\hline
15 & 1^{-1} 3^{-1} 5^{-1} 15^{-1} & 1,3,5,15 & -2 \\
15 & 1^{-2} 3^{1} 5^{1} 15^{-2} & 1,15 & -1 \\
15 & 1^{1} 3^{-2} 5^{-2} 15^{1} & 3,5 & -1 \\
\hline
20 & 2^{-2} 10^{-2} & 1,2,4,5,10,20 & -2 \\
20 & 1^{-2} 2^{-1} 5^{2} 10^{-3} & 1,2,4,10,20^\dagger & -2 \\
20 & 1^{2} 2^{-3} 5^{-2} 10^{-1} & 2,4,5,10,20^\dagger & -2 \\
20 & 1^{-1} 2^{-2} 5^{-3} 10^{2} & 1,2,4,5^\dagger & -2 \\
\hline
\end{array}
$$
\caption{$\eta$--products with reflective singularities. (\emph{ctd.})}
\end{table}

\addtocounter{table}{-1}
\begin{table}[p]
$$
\begin{array}{|c||c|c|c|}
\hline
\hbox{\bf Level} & \hbox{\bf Shape} & \hbox{\bf Poles (1/n)} & \hbox{\bf
Weight} \\
\hline
\hline
20 & 1^{-1} 2^{-3} 4^{2} 5^{-3} 10^{3} 20^{-2} & 1,2,5,20^\dagger & -2 \\
20 & 1^{-3} 2^{3} 4^{-2} 5^{-1} 10^{-3} 20^{2} & 1,4,5,10^\dagger & -2 \\
20 & 1^{-3} 2^{2} 5^{-1} 10^{-2} & 1,5,10,20^\dagger & -2 \\
20 & 1^{-2} 2^{3} 4^{-2} 5^{-2} 10^{3} 20^{-2} & 1,4,5,20 & -1 \\
20 & 1^{1} 2^{-1} 4^{-1} 5^{-1} 10^{-1} 20^{1} & 2,4,5,10 & -1 \\
20 & 1^{-1} 2^{-1} 4^{1} 5^{1} 10^{-1} 20^{-1} & 1,2,10,20 & -1 \\
20 & 1^{-1} 4^{-1} 5^{1} 10^{-2} 20^{1} & 1,2,4,10 & -1 \\
20 & 1^{1} 2^{-2} 4^{1} 5^{-1} 20^{-1} & 2,5,10,20 & -1 \\
20 & 2^{-2} 5^{-4} 10^{8} 20^{-4} & 1,2,4,5,20 & -1 \\
20 & 1^{-4} 2^{8} 4^{-4} 10^{-2} & 1,4,5,10,20 & -1 \\
20 & 1^{-1} 2^{1} 4^{-2} 5^{-1} 10^{-2} 20^{4} & 1,2,4,5,10 & -1/2 \\
20 & 1^{-2} 2^{1} 4^{-1} 5^{4} 10^{-2} 20^{-1} & 1,2,4,10,20 & -1/2 \\
20 & 1^{-1} 2^{-2} 4^{4} 5^{-1} 10^{1} 20^{-2} & 1,2,5,10,20 & -1/2 \\
20 & 1^{4} 2^{-2} 4^{-1} 5^{-2} 10^{1} 20^{-1} & 2,4,5,10,20 & -1/2 \\
20 & 2^{-2} 5^{-6} 10^{13} 20^{-6} & 1,2,4,5,20 & -1/2 \\
20 & 1^{-6} 2^{13} 4^{-6} 10^{-2} & 1,4,5,10,20 & -1/2 \\
20 & 1^{-1} 5^{-1} 10^{-1} 20^{2} & 1,2,5,10 & -1/2 \\
20 & 1^{-1} 2^{-1} 4^{2} 5^{-1} & 1,2,5,10 & -1/2 \\
20 & 4^{-1} 5^{2} 10^{-1} 20^{-1} & 2,4,10,20 & -1/2 \\
20 & 1^{2} 2^{-1} 4^{-1} 20^{-1} & 2,4,10,20 & -1/2 \\
20 & 1^{-2} 2^{3} 4^{-2} 5^{-4} 10^{8} 20^{-4} & 1,4,5,20 & -1/2 \\
20 & 1^{-4} 2^{8} 4^{-4} 5^{-2} 10^{3} 20^{-2} & 1,4,5,20 & -1/2 \\
20 & 1^{-1} 4^{-1} 5^{-1} 10^{3} 20^{-1} & 1,2,4 & -1/2 \\
20 & 1^{-1} 2^{3} 4^{-1} 5^{-1} 20^{-1} & 5,10,20 & -1/2 \\
20 & 1^{-1} 2^{-1} 4^{1} 5^{-1} 10^{4} 20^{-3} & 1,2,20 & -1/2 \\
\hline
\end{array}
$$
\caption{$\eta$--products with reflective singularities. (\emph{ctd.})}
\end{table}

\addtocounter{table}{-1}
\begin{table}[p]
$$
\begin{array}{|c||c|c|c|}
\hline
\hbox{\bf Level} & \hbox{\bf Shape} & \hbox{\bf Poles (1/n)} & \hbox{\bf
Weight} \\
\hline
\hline
20 & 1^{-3} 2^{4} 4^{-1} 5^{1} 10^{-1} 20^{-1} & 1,10,20 & -1/2 \\
20 & 1^{1} 2^{-1} 4^{-1} 5^{-3} 10^{4} 20^{-1} & 2,4,5 & -1/2 \\
20 & 1^{-1} 2^{4} 4^{-3} 5^{-1} 10^{-1} 20^{1} & 4,5,10 & -1/2 \\
20 & 1^{1} 2^{-2} 4^{1} 5^{-3} 10^{5} 20^{-3} & 2,5,20 & -1/2 \\
20 & 1^{-3} 2^{5} 4^{-3} 5^{1} 10^{-2} 20^{1} & 1,4,10 & -1/2 \\
20 & \vdots & \vdots{}^\dagger & \vdots \\
\hline
23 & 1^{-1} 23^{-1} & 1,23 & -1 \\
\hline
28 & 1^{-1} 2^{-1} 7^{-1} 14^{-1} & 1,2,4,7,14,28^\dagger & -2 \\
28 & 1^{-1} 2^{-1} 7^{-5} 14^{9} 28^{-4} & 1,2,4,7,28^\dagger & -1 \\
28 & 1^{-5} 2^{9} 4^{-4} 7^{-1} 14^{-1} & 1,4,7,14,28^\dagger & -1 \\
28 & 1^{-1} 2^{1} 4^{-1} 7^{-1} 14^{1} 28^{-1} & 1,4,7,28 & -1 \\
28 & 1^{-1} 2^{-2} 4^{2} 7^{-1} 14^{-2} 28^{2} & 1,2,7,14^\dagger & -1 \\
28 & 1^{-3} 2^{4} 4^{-2} 7^{-3} 14^{4} 28^{-2} & 1,4,7,28^\dagger & -1 \\
28 & 2^{-2} 4^{1} 14^{-2} 28^{1} & 1,2,7,14 & -1 \\
28 & 1^{1} 2^{-2} 7^{1} 14^{-2} & 2,4,14,28 & -1 \\
28 & 1^{-2} 2^{1} 7^{-2} 14^{1} & 1,7^\dagger & -1 \\
28 & 1^{-1} 2^{1} 4^{-1} 7^{-3} 14^{6} 28^{-3} & 1,4,7,28 & -1/2 \\
28 & 1^{-3} 2^{6} 4^{-3} 7^{-1} 14^{1} 28^{-1} & 1,4,7,28 & -1/2 \\
28 & 1^{-3} 2^{4} 4^{-2} 7^{-5} 14^{9} 28^{-4} & 1,4,7,28^\dagger & -1/2 \\
28 & 1^{-5} 2^{9} 4^{-4} 7^{-3} 14^{4} 28^{-2} & 1,4,7,28^\dagger & -1/2 \\
28 & 1^{-2} 2^{1} 4^{-1} 28^{1} & 1,2,4^\dagger & -1/2 \\
28 & 4^{1} 7^{-2} 14^{1} 28^{-1} & 7,14,28^\dagger & -1/2 \\
28 & 1^{-2} 2^{1} 7^{-4} 14^{6} 28^{-2} & 1,7^\dagger & -1/2 \\
28 & 1^{-4} 2^{6} 4^{-2} 7^{-2} 14^{1} & 1,7^\dagger & -1/2 \\
\hline
30 & 2^{-1} 3^{-1} 5^{-1} 6^{1} 10^{1} 30^{-1} & 1,2,3,5,15,30 & -1 \\
\hline
\end{array}
$$
\caption{$\eta$--products with reflective singularities. (\emph{ctd.})}
\end{table}

\addtocounter{table}{-1}
\begin{table}[p]
$$
\begin{array}{|c||c|c|c|}
\hline
\hbox{\bf Level} & \hbox{\bf Shape} & \hbox{\bf Poles (1/n)} & \hbox{\bf
Weight} \\
\hline
\hline
30 & 1^{-1} 3^{1} 5^{1} 6^{-1} 10^{-1} 15^{-1} & 1,2,6,10,15,30 & -1 \\
30 & 1^{-1} 2^{1} 6^{-1} 10^{-1} 15^{-1} 30^{1} & 1,3,5,6,10,15 & -1 \\
30 & 1^{1} 2^{-1} 3^{-1} 5^{-1} 15^{1} 30^{-1} & 2,3,5,6,10,30 & -1 \\
\hline
44 & 2^{-1} 22^{-1} & 1,2,4,11,22,44 & -1 \\
\hline
\end{array}
$$
\caption{$\eta$--products with reflective singularities. (\emph{ctd.})}
\end{table}

\chapter{Obstruction Spaces for Reflective Forms}
\label{obstructions}

We list some of the results of a {\tt PARI} program to calculate the
dimension of certain obstruction spaces.  These spaces were defined in
\cite{borcherds_gkz} and a formula to compute their dimensions was
presented in \cite{borcherds_refl}.

If a discriminant form is preceded by a star ($\star$) this means that
one of the components of the discriminant form is regular.

\begin{sidewaystable}[p]
$$
\begin{array}{|c|c|c|c|c|c|c|c|c|c|c|c|c|c|c|c|c|c|c|c|c|c|c|c|c|}
\hline
\hbox{\bf Discriminant} & 1 & 2 & 3 & 4 & 5 & 6 & 7 & 8 & 9 & 10 & 11
& 12 & 13 & 14 & 15 & 16 & 17 & 18 & 19 & 20 & 21 & 22 & 23 & 24 \\
\hline\hline
{}^\star 1 &
- & - & - & 0 & - & - & - & 0 & - & - & - & 0 &
- & - & - & 0 & - & - & - & 1 & - & - & - & 0 \\
\hline
2^{-2} &
- & - & - & 0 & - & - & - & 0 & - & - & - & 1 &
- & - & - & 1 & - & - & - & 1 & - & - & - & 2 \\
{}^\star 2^{+2} &
- & - & - & 0 & - & - & - & 0 & - & - & - & 1 &
- & - & - & 1 & - & - & - & 2 & - & - & - & 2 \\
\hline
3^{-1} &
- & 0 & - & - & - & 0 & - & - & - & 0 & - & - &
- & 1 & - & - & - & 1 & - & - & - & 1 & - & - \\
3^{+1} &
- & 0 & - & - & - & 0 & - & - & - & 1 & - & - &
- & 1 & - & - & - & 1 & - & - & - & 2 & - & - \\
3^{+2} &
- & - & - & 0 & - & - & - & 1 & - & - & - & 1 &
- & - & - & 2 & - & - & - & 2 & - & - & - & 3 \\
{}^\star 3^{-2} &
- & - & - & 0 & - & - & - & 1 & - & - & - & 1 &
- & - & - & 2 & - & - & - & 3 & - & - & - & 3 \\
\hline
2_1^{+1} &
- & - & 0 & - & - & - & 0 & - & - & - & 0 & - &
- & - & 1 & - & - & - & 1 & - & - & - & 1 & - \\
2_7^{+1} &
0 & - & - & - & 0 & - & - & - & 1 & - & - & - &
1 & - & - & - & 1 & - & - & - & 2 & - & - & - \\
2_0^{+2} &
- & - & - & 0 & - & - & - & 1 & - & - & - & 1 &
- & - & - & 2 & - & - & - & 3 & - & - & - & 3 \\
2_2^{+2} &
- & 0 & - & - & - & 0 & - & - & - & 1 & - & - &
- & 1 & - & - & - & 2 & - & - & - & 2 & - & - \\
2_6^{+2} &
- & 0 & - & - & - & 1 & - & - & - & 1 & - & - &
- & 2 & - & - & - & 2 & - & - & - & 3 & - & - \\
2_3^{+3} &
0 & - & - & - & 1 & - & - & - & 1 & - & - & - &
2 & - & - & - & 3 & - & - & - & 3 & - & - & - \\
2_5^{+3} &
- & - & 0 & - & - & - & 1 & - & - & - & 2 & - &
- & - & 2 & - & - & - & 3 & - & - & - & 4 & - \\
{}^\star 2_1^{+3} &
- & - & 0 & - & - & - & 1 & - & - & - & 2 & - &
- & - & 3 & - & - & - & 4 & - & - & - & 5 & - \\
{}^\star 2_7^{+3} &
0 & - & - & - & 1 & - & - & - & 2 & - & - & - &
3 & - & - & - & 4 & - & - & - & 5 & - & - & - \\
\hline
\end{array}
$$
\caption{Dimensions of obstruction spaces.}
\end{sidewaystable}

\addtocounter{table}{-1}
\begin{sidewaystable}[p]
$$
\begin{array}{|c|c|c|c|c|c|c|c|c|c|c|c|c|c|c|c|c|c|c|c|c|c|c|c|c|}
\hline
\hbox{\bf Discriminant} & 1 & 2 & 3 & 4 & 5 & 6 & 7 & 8 & 9 & 10 & 11
& 12 & 13 & 14 & 15 & 16 & 17 & 18 & 19 & 20 & 21 & 22 & 23 & 24 \\
\hline\hline
5^{-1} &
- & - & - & 0 & - & - & - & 1 & - & - & - & 1 &
- & - & - & 2 & - & - & - & 2 & - & - & - & 3 \\
5^{+1} &
- & - & - & 0 & - & - & - & 1 & - & - & - & 1 &
- & - & - & 2 & - & - & - & 2 & - & - & - & 3 \\
5^{-2} &
- & - & - & 1 & - & - & - & 1 & - & - & - & 3 &
- & - & - & 3 & - & - & - & 4 & - & - & - & 5 \\
{}^\star 5^{+2} &
- & - & - & 1 & - & - & - & 1 & - & - & - & 3 &
- & - & - & 3 & - & - & - & 5 & - & - & - & 5 \\
\hline
2^{-2} 3^{+1} &
- & 0 & - & - & - & 1 & - & - & - & 1 & - & - &
- & 2 & - & - & - & 3 & - & - & - & 3 & - & - \\
2^{-2} 3^{-1} &
- & 0 & - & - & - & 1 & - & - & - & 2 & - & - &
- & 2 & - & - & - & 3 & - & - & - & 4 & - & - \\
2^{-2} 3^{+2} &
- & - & - & 1 & - & - & - & 2 & - & - & - & 3 &
- & - & - & 4 & - & - & - & 5 & - & - & - & 6 \\
{}^\star 2^{+2} 3^{+1} &
- & 0 & - & - & - & 1 & - & - & - & 2 & - & - &
- & 3 & - & - & - & 4 & - & - & - & 5 & - & - \\
{}^\star 2^{+2} 3^{-1} &
- & 0 & - & - & - & 1 & - & - & - & 2 & - & - &
- & 3 & - & - & - & 4 & - & - & - & 5 & - & - \\
{}^\star 2^{-2} 3^{-2} &
- & - & - & 1 & - & - & - & 2 & - & - & - & 4 &
- & - & - & 5 & - & - & - & 6 & - & - & - & 8 \\
{}^\star 2^{+2} 3^{+2} &
- & - & - & 1 & - & - & - & 3 & - & - & - & 4 &
- & - & - & 6 & - & - & - & 7 & - & - & - & 9 \\
{}^\star 2^{+2} 3^{-2} &
- & - & - & 1 & - & - & - & 3 & - & - & - & 5 &
- & - & - & 7 & - & - & - & 9 & - & - & - & 11 \\
\hline
7^{+1} &
- & 0 & - & - & - & 0 & - & - & - & 1 & - & - &
- & 2 & - & - & - & 2 & - & - & - & 3 & - & - \\
7^{-1} &
- & 1 & - & - & - & 1 & - & - & - & 2 & - & - &
- & 3 & - & - & - & 3 & - & - & - & 4 & - & - \\
7^{+2} &
- & - & - & 1 & - & - & - & 3 & - & - & - & 3 &
- & - & - & 5 & - & - & - & 6 & - & - & - & 7 \\
{}^\star 7^{-2} &
- & - & - & 1 & - & - & - & 3 & - & - & - & 3 &
- & - & - & 5 & - & - & - & 7 & - & - & - & 7 \\
\hline
\end{array}
$$
\caption{Dimensions of obstruction spaces. (\emph{ctd.})}
\end{sidewaystable}

\addtocounter{table}{-1}
\begin{sidewaystable}[p]
$$
\begin{array}{|c|c|c|c|c|c|c|c|c|c|c|c|c|c|c|c|c|c|c|c|c|c|c|c|c|}
\hline
\hbox{\bf Discriminant} & 1 & 2 & 3 & 4 & 5 & 6 & 7 & 8 & 9 & 10 & 11
& 12 & 13 & 14 & 15 & 16 & 17 & 18 & 19 & 20 & 21 & 22 & 23 & 24 \\
\hline\hline
2^{-2} 5^{+1} &
- & - & - & 1 & - & - & - & 2 & - & - & - & 3 &
- & - & - & 4 & - & - & - & 5 & - & - & - & 6 \\
2^{-2} 5^{-1} &
- & - & - & 1 & - & - & - & 2 & - & - & - & 3 &
- & - & - & 4 & - & - & - & 5 & - & - & - & 6 \\
{}^\star 2^{+2} 5^{+1} &
- & - & - & 1 & - & - & - & 3 & - & - & - & 4 &
- & - & - & 6 & - & - & - & 7 & - & - & - & 9 \\
{}^\star 2^{+2} 5^{-1} &
- & - & - & 1 & - & - & - & 3 & - & - & - & 4 &
- & - & - & 6 & - & - & - & 7 & - & - & - & 9 \\
2^{-2} 5^{-2} &
- & - & - & 2 & - & - & - & 4 & - & - & - & 5 &
- & - & - & 7 & - & - & - & 9 & - & - & - & 10 \\
{}^\star 2^{-2} 5^{+2} &
- & - & - & 2 & - & - & - & 4 & - & - & - & 6 &
- & - & - & 8 & - & - & - & 10 & - & - & - & 12 \\
{}^\star 2^{+2} 5^{-2} &
- & - & - & 3 & - & - & - & 5 & - & - & - & 8 &
- & - & - & 10 & - & - & - & 13 & - & - & - & 15 \\
{}^\star 2^{+2} 5^{+2} &
- & - & - & 3 & - & - & - & 5 & - & - & - & 9 &
- & - & - & 11 & - & - & - & 15 & - & - & - & 17 \\
\hline
11^{-1} &
- & 0 & - & - & - & 1 & - & - & - & 2 & - & - &
- & 3 & - & - & - & 4 & - & - & - & 5 & - & - \\
11^{+1} &
- & 1 & - & - & - & 2 & - & - & - & 3 & - & - &
- & 4 & - & - & - & 5 & - & - & - & 6 & - & - \\
11^{+2} &
- & - & - & 2 & - & - & - & 4 & - & - & - & 6 &
- & - & - & 8 & - & - & - & 9 & - & - & - & 12 \\
{}^\star 11^{-2} &
- & - & - & 2 & - & - & - & 4 & - & - & - & 6 &
- & - & - & 8 & - & - & - & 10 & - & - & - & 12 \\
\hline
13^{-1} &
- & - & - & 1 & - & - & - & 3 & - & - & - & 3 &
- & - & - & 5 & - & - & - & 6 & - & - & - & 7 \\
13^{+1} &
- & - & - & 1 & - & - & - & 3 & - & - & - & 3 &
- & - & - & 5 & - & - & - & 6 & - & - & - & 7 \\
13^{-2} &
- & - & - & 3 & - & - & - & 5 & - & - & - & 7 &
- & - & - & 9 & - & - & - & 12 & - & - & - & 13 \\
{}^\star 13^{+2} &
- & - & - & 3 & - & - & - & 5 & - & - & - & 7 &
- & - & - & 9 & - & - & - & 13 & - & - & - & 13 \\
\hline
\end{array}
$$
\caption{Dimensions of obstruction spaces. (\emph{ctd.})}
\end{sidewaystable}

\addtocounter{table}{-1}
\begin{sidewaystable}[p]
$$
\begin{array}{|c|c|c|c|c|c|c|c|c|c|c|c|c|c|c|c|c|c|c|c|c|c|c|c|c|}
\hline
\hbox{\bf Discriminant} & 1 & 2 & 3 & 4 & 5 & 6 & 7 & 8 & 9 & 10 & 11
& 12 & 13 & 14 & 15 & 16 & 17 & 18 & 19 & 20 & 21 & 22 & 23 & 24 \\
\hline\hline
2^{-2} 7^{+1} &
- & 0 & - & - & - & 2 & - & - & - & 3 & - & - &
- & 4 & - & - & - & 6 & - & - & - & 7 & - & - \\
2^{-2} 7^{-1} &
- & 1 & - & - & - & 3 & - & - & - & 4 & - & - &
- & 5 & - & - & - & 7 & - & - & - & 8 & - & - \\
{}^\star 2^{+2} 7^{+1} &
- & 0 & - & - & - & 2 & - & - & - & 4 & - & - &
- & 6 & - & - & - & 8 & - & - & - & 10 & - & - \\
{}^\star 2^{+2} 7^{-1} &
- & 2 & - & - & - & 4 & - & - & - & 6 & - & - &
- & 8 & - & - & - & 10 & - & - & - & 12 & - & - \\
2^{-2} 7^{+2} &
- & - & - & 3 & - & - & - & 5 & - & - & - & 8 &
- & - & - & 10 & - & - & - & 12 & - & - & - & 15 \\
{}^\star 2^{-2} 7^{-2} &
- & - & - & 3 & - & - & - & 5 & - & - & - & 9 &
- & - & - & 11 & - & - & - & 13 & - & - & - & 17 \\
{}^\star 2^{+2} 7^{+2} &
- & - & - & 4 & - & - & - & 8 & - & - & - & 11 &
- & - & - & 15 & - & - & - & 18 & - & - & - & 22 \\
{}^\star 2^{+2} 7^{-2} &
- & - & - & 4 & - & - & - & 8 & - & - & - & 12 &
- & - & - & 16 & - & - & - & 20 & - & - & - & 24 \\
\hline
3^{+1} 5^{+1} &
- & 0 & - & - & - & 1 & - & - & - & 2 & - & - &
- & 3 & - & - & - & 4 & - & - & - & 5 & - & - \\
3^{-1} 5^{+1} &
- & 1 & - & - & - & 2 & - & - & - & 3 & - & - &
- & 4 & - & - & - & 5 & - & - & - & 6 & - & - \\
3^{+2} 5^{+1} &
- & - & - & 2 & - & - & - & 3 & - & - & - & 5 &
- & - & - & 6 & - & - & - & 8 & - & - & - & 9 \\
3^{-2} 5^{+1} &
- & - & - & 2 & - & - & - & 4 & - & - & - & 6 &
- & - & - & 8 & - & - & - & 10 & - & - & - & 12 \\
3^{+1} 5^{-1} &
- & 1 & - & - & - & 2 & - & - & - & 3 & - & - &
- & 4 & - & - & - & 5 & - & - & - & 6 & - & - \\
3^{-1} 5^{-1} &
- & 0 & - & - & - & 1 & - & - & - & 2 & - & - &
- & 3 & - & - & - & 4 & - & - & - & 5 & - & - \\
3^{+2} 5^{-1} &
- & - & - & 2 & - & - & - & 3 & - & - & - & 5 &
- & - & - & 6 & - & - & - & 8 & - & - & - & 9 \\
\hline
\end{array}
$$
\caption{Dimensions of obstruction spaces. (\emph{ctd.})}
\end{sidewaystable}

\addtocounter{table}{-1}
\begin{sidewaystable}[p]
$$
\begin{array}{|c|c|c|c|c|c|c|c|c|c|c|c|c|c|c|c|c|c|c|c|c|c|c|c|c|}
\hline
\hbox{\bf Discriminant} & 1 & 2 & 3 & 4 & 5 & 6 & 7 & 8 & 9 & 10 & 11
& 12 & 13 & 14 & 15 & 16 & 17 & 18 & 19 & 20 & 21 & 22 & 23 & 24 \\
\hline\hline
3^{-2} 5^{-1} &
- & - & - & 2 & - & - & - & 4 & - & - & - & 6 &
- & - & - & 8 & - & - & - & 10 & - & - & - & 12 \\
3^{+1} 5^{+2} &
- & 1 & - & - & - & 3 & - & - & - & 5 & - & - &
- & 7 & - & - & - & 9 & - & - & - & 11 & - & - \\
3^{-1} 5^{+2} &
- & 1 & - & - & - & 3 & - & - & - & 5 & - & - &
- & 7 & - & - & - & 9 & - & - & - & 11 & - & - \\
3^{+2} 5^{+2} &
- & - & - & 3 & - & - & - & 7 & - & - & - & 9 &
- & - & - & 13 & - & - & - & 15 & - & - & - & 19 \\
3^{-2} 5^{+2} &
- & - & - & 4 & - & - & - & 8 & - & - & - & 12 &
- & - & - & 16 & - & - & - & 20 & - & - & - & 24 \\
3^{+1} 5^{-2} &
- & 1 & - & - & - & 3 & - & - & - & 4 & - & - &
- & 6 & - & - & - & 8 & - & - & - & 9 & - & - \\
3^{-1} 5^{-2} &
- & 1 & - & - & - & 3 & - & - & - & 5 & - & - &
- & 6 & - & - & - & 8 & - & - & - & 10 & - & - \\
3^{+2} 5^{-2} &
- & - & - & 3 & - & - & - & 6 & - & - & - & 8 &
- & - & - & 11 & - & - & - & 13 & - & - & - & 16 \\
3^{-2} 5^{-2} &
- & - & - & 4 & - & - & - & 7 & - & - & - & 11 &
- & - & - & 14 & - & - & - & 17 & - & - & - & 21 \\
\hline
\end{array}
$$
\caption{Dimensions of obstruction spaces. (\emph{ctd.})}
\end{sidewaystable}

\chapter{Singular Weight Forms}\label{singular_weight}

We list the results of a {\tt PARI} program (and hand calculation) to
find singular weight forms suitable for constructing generalized
Kac--Moody algebras.  Currently this is only compiled for elementary
lattices with no odd $2$--Jordan blocks.  All forms which occur seem
to be $\eta$--products (listed in the Shape column) and all seem to
correspond to regular lattices (which seriously limits the number that
can occur).  The lattice for which this gives a singular weight form
is given in the Genus column.  Group refers to the invariance group of
the $\eta$--product.  $Co_0$ gives the class in Conway's group which
corresponds to the Shape (see \cite{conway_norton} and
\cite{kondo_leech}).

$$
\begin{array}{|c||c|c|c|c||c|}
\hline
\hbox{\bf Level} & \hbox{\bf Shape} & \hbox{\bf Genus} & \hbox{\bf
Group} & \hbox{\bf $Co_0$} & \hbox{\bf Notes}\\
\hline
\hline
1 & 1^{-24} & \even_{2,26} & \Gamma_0(1) & 1A+ & \hbox{Monster $M$} \\
\hline
2 & 1^{-8} 2^{-8} & \even_{2,18}(2^{+10}) & \Gamma_0(2)+ & 2A+ &
\hbox{Baby Monster $B$} \\
2 & 1^{-16} 2^{8} & \even_{2,10}(2^{+2}) & \Gamma_0(2) & (2A-) &
\hbox{Conway $Co_1$} \\
2 & 1^{8} 2^{-16} & \even_{2,10}(2^{+10}) & \Gamma_0(2) & 2A- & \\
\hline
3 & 1^{-6} 3^{-6} & \even_{2,14}(3^{-8}) & \Gamma_0(3)+ & 3B+ &
\hbox{Fischer $Fi_{24}$} \\
3 & 1^{-9} 3^{3} & \even_{2,8}(3^{+3}) & \Gamma_0(3) & (3C+) &
\hbox{Suzuki $Suz$} \\
3 & 1^{3} 3^{-9} & \even_{2,8}(3^{+7}) & \Gamma_0(3) & 3C+ & \\
\hline
\end{array}
$$

$$
\begin{array}{|c||c|c|c|c||c|}
\hline
\hbox{\bf Level} & \hbox{\bf Shape} & \hbox{\bf Genus} & \hbox{\bf
Group} & \hbox{\bf $Co_0$} & \hbox{\bf Notes}\\
\hline
\hline
5 & 1^{-4} 5^{-4} & \even_{2,10}(5^{+6}) & \Gamma_0(5)+ & 5B+ &
\hbox{$HN$} \\
5 & 1^{-5} 5^{1} & \even_{2,6}(5^{+3}) & \Gamma_0(5) & (5C+) &
\hbox{$HJ$} \\
5 & 1^{1} 5^{-5} & \even_{2,6}(5^{+5}) & \Gamma_0(5) & 5C+ & \\
\hline
6 & 1^{-2} 2^{-2} 3^{-2} 6^{-2} & \even_{2,10}(2^{+6}3^{-6}) & 
\Gamma_0(6)+ & 6E+ & \hbox{Fischer $Fi_{22}$} \\
6 & 1^{-1} 2^{-4} 3^{-5} 6^{4} & \even_{2,8}(2^{+2} 3^{+3}) &
\Gamma_0(6) & (6C+) & \\
6 & 1^{-4} 2^{-1} 3^{4} 6^{-5} & \even_{2,8}(2^{-8} 3^{-3}) &
\Gamma_0(6) & 6C+ & \\
6 & 1^{-5} 2^{4} 3^{-1} 6^{-4} & \even_{2,8}(2^{+2} 3^{+7}) &
\Gamma_0(6) & 6D- & \\
6 & 1^{4} 2^{-5} 3^{-4} 6^{-1} & \even_{2,8}(2^{-8} 3^{-7}) &
\Gamma_0(6) & 6C- & \\
6 & 1^{-3} 2^{-3} 3^{1} 6^{1} & \even_{2,6}(2^{-4} 3^{-2}) &
\Gamma_0(6)+2 & (6F+) & \\
6 & 1^{-4} 2^{2} 3^{-4} 6^{2} & \even_{2,6}(2^{+2} 3^{-4}) &
\Gamma_0(6)+3 & (6E-) & \\
6 & 1^{2} 2^{-4} 3^{2} 6^{-4} & \even_{2,6}(2^{+6} 3^{-4}) &
\Gamma_0(6)+3 & 6E- & \\
6 & 1^{1} 2^{1} 3^{-3} 6^{-3} & \even_{2,6}(2^{-4} 3^{-6}) &
\Gamma_0(6)+2 & 6F+ & \\
6 & 1^{-6} 2^{3} 3^{2} 6^{-1} & \even_{2,4}(2^{+2} 3^{-3}) &
\Gamma_0(6) & (6F-) & \\
6 & 1^{3} 2^{-6} 3^{-1} 6^{2} & \even_{2,4}(2^{-4} 3^{+3}) &
\Gamma_0(6) & (6F-) & \\
6 & 1^{2} 2^{-1} 3^{-6} 6^{3} & \even_{2,4}(2^{+2} 3^{-3}) &
\Gamma_0(6) & (6F-) & \\
6 & 1^{-1} 2^{2} 3^{3} 6^{-6} & \even_{2,4}(2^{-4} 3^{+3}) &
\Gamma_0(6) & 6F- & \\
\hline
7 & 1^{-3} 7^{-3} & \even_{2,8}(7^{+5}) & \Gamma_0(7)+ &
7B+ & \hbox{Held $He$} \\
\hline
10 & 1^{-1} 2^{-2} 5^{-3} 10^{2} & \even_{2,6}(2^{+2} 5^{+3}) &
\Gamma_0(10) & (10D+) & \\
10 & 1^{-2} 2^{-1} 5^{2} 10^{-3} & \even_{2,6}(2^{-6} 5^{-3}) &
\Gamma_0(10) & 10D+ & \\
10 & 1^{-3} 2^{2} 5^{-1} 10^{-2} & \even_{2,6}(2^{+2} 5^{+5}) &
\Gamma_0(10) & 10E- & \\
10 & 1^{2} 2^{-3} 5^{-2} 10^{-1} & \even_{2,6}(2^{-6} 5^{-5}) &
\Gamma_0(10) & 10D- & \\
\hline
11 & 1^{-2} 11^{-2} & \even_{2,6}(11^{-4}) & \Gamma_0(11)+ &
11A+ & \hbox{Mathieu $M_{12}$} \\
\hline
14 & 1^{-1} 2^{-1} 7^{-1} 14^{-1} & \even_{2,6}(2^{+4} 7^{-4}) &
\Gamma_0(14)+ & 14B+ & \\
14 & 1^{-2} 2^{1} 7^{-2} 14^{1} & \even_{2,4}(2^{+2} 7^{+3}) &
\Gamma_0(14)+7 & (14B-) & \\
14 & 1^{1} 2^{-2} 7^{1} 14^{-2} & \even_{2,4}(2^{+4} 7^{+3}) &
\Gamma_0(14)+7 & 14B- & \\
\hline
15 & 1^{-1} 3^{-1} 5^{-1} 15^{-1} & \even_{2,6}(3^{+4} 5^{-4}) &
\Gamma_0(15)+ & 15D+ & \\
15 & 1^{1} 3^{-2} 5^{-2} 15^{1} & \even_{2,4}(3^{+3} 5^{+3}) &
\Gamma_0(15)+15 & (15E+) & \\
15 & 1^{-2} 3^{1} 5^{1} 15^{-2} & \even_{2,4}(3^{-3} 5^{-3}) &
\Gamma_0(15)+15 & 15E+ & \\
\hline
23 & 1^{-1} 23^{-1} & \even_{2,4}(23^{+3}) & \Gamma_0(23)+ &
23A/23B & \\
\hline
30 & 1^{-1} 2^{1} 6^{-1} 10^{-1} 15^{-1} 30^{1} & 
\even_{2,4}(2^{+2} 3^{+3} 5^{+3}) & \Gamma_0(30)+15 & (30E-) & \\
30 & 2^{-1} 3^{-1} 5^{-1} 6^{1} 10^{1} 30^{-1} &
\even_{2,4}(2^{+2} 3^{-3} 5^{-3}) & \Gamma_0(30)+15 & 30E- & \\
30 & 1^{-1} 3^{1} 5^{1} 6^{-1} 10^{-1} 15^{-1} &
\even_{2,4}(2^{+4} 3^{+3} 5^{+3}) & \Gamma_0(30)+15 & 30D+ & \\
30 & 1^{1} 2^{-1} 3^{-1} 5^{-1} 15^{1} 30^{-1} &
\even_{2,4}(2^{+4} 3^{-3} 5^{-3}) & \Gamma_0(30)+15 & 30D- & \\
\hline
\end{array}
$$

\end{document}